\documentclass[a4paper,12pt]{amsart}

\usepackage[driverfallback=dvipdfm]{hyperref}
\usepackage[all]{xy}
\usepackage{tikz}
\usepackage{amsaddr}
\usepackage{amsthm}

\newtheorem{thm}{Theorem}[section]
\newtheorem*{thm*}{Theorem}
\newtheorem{cor}[thm]{Corollary}
\newtheorem{lem}[thm]{Lemma}
\newtheorem{prop}[thm]{Proposition}

\theoremstyle{definition}
\newtheorem{defn}[thm]{Definition}

\newtheorem{nt}[thm]{Notation}
\newtheorem{rem}[thm]{Remark}
\newtheorem{ex}[thm]{Example}

\newtheorem*{theoremaux}{Theorem \theoremauxnum}
\gdef\theoremauxnum{1}

\newenvironment{theoremff}[2][]{%
  \def\theoremauxnum{\ref{#2}}
  \begin{theoremaux}[#1]
}{%
  \end{theoremaux}
}

\def \p{{\mathbb P}}          

\def\N{{\mathbb N}}           
\def\R{{\mathbb R}}           
\def\O{{\mathcal O}}          
\def\G{{\mathcal G}}          
\def\e{\varepsilon}           

\def\Aut{\operatorname{Aut}}
\def\Out{\operatorname{Out}}
\def\Lip{\operatorname{Lip}}
\def\Hor{\operatorname{Hor}}
\def\PL{\operatorname{PL}}
\def\opt{\operatorname{opt}}
\def\wopt{\operatorname{weakopt}}
\def\vol{\operatorname{vol}}
\def\Min{\operatorname{Min}}
\def\LocMin{\operatorname{LocMin}}
\def\TT{\operatorname{TT}}
\def\TTo{\operatorname{TT_0}}
\def\rank{\operatorname{rank}}
\def\core{\operatorname{core}}

\def\simfk{\langle \sim_{f^k}\rangle}
\def\wt{\widetilde}

\title[Level set of automorphisms are connected]{On the connectivity of level sets of automorphisms of free groups, with applications to
  decision problems}
\author{Stefano Francaviglia}
\address{Dipartimento di Matematica of the University of
Bologna}
\email{stefano.francaviglia@unibo.it}
\author{Armando Martino}
\address{Mathematical Sciences, University of Southampton }
\email{A.Martino@soton.ac.uk}
\begin{document}


\begin{abstract}
We show that the level sets of automorphisms of free groups with respect to the Lipschitz
metric are connected as subsets of Culler-Vogtmann space. In fact we prove our result in a 
 more general setting of deformation spaces. As applications, we give metric solutions of the
 conjugacy problem for irreducible automorphisms  and
 the detection of reducibility. We additionally prove technical results that may be of independent
interest --- such as the fact that the set of displacements is well ordered.
\end{abstract}
\maketitle
\tableofcontents

\section{Introduction}

We consider $F_n$ the free group of rank $n$, usually with a basis $B$ (a free generating set). We are interested in the automorphism group, $\Aut(F_n)$ and the Outer automorphism group, which is defined as $\Out(F_n)= \Aut(F_n)/\operatorname{Inn}(F_n)$. 

In recent years there has been a great deal of attention given to the Lipschitz metric on $CV_n$, Culler-Vogtmann space, see \cite{MR2862155}, \cite{MR2863547}, \cite{BestvinaBers} for instance. It has been considered even more generally in \cite{MR3342683}.

The main goal of the paper, Theorem~\ref{tconnected}, is to prove a result about the
connectedness of the level sets of the displacement function $\lambda_\phi$ of $\phi \in \Out(F_n)$. That is, one considers the Lipschitz
metric then one can take the infimum of all displacements of points for $\phi$ in $CV_n$. This
infimum may not be realised in general. However, we show that for any $\epsilon > 0$, the set
of all points of $CV_n$ displaced at most $\epsilon$ more than this infimum is
connected. Formally\footnote{Theorem~\ref{tconnected} holds in a general setting, the result
in $CV_n$ is for $\Gamma=F_n$. Here $\overline{\O(\Gamma)}^\infty$ denotes the simplicial
bordification of the outer space. Precise definitions are given through the paper.},

\begin{theoremff}[Level sets are connected]{tconnected}
Let $\phi \in\Out(\Gamma)$. For any $\e > 0 $ the set
$$\{X\in\O(\Gamma):\lambda_\phi(X)\leq \lambda(\phi) + \e \}$$ is connected in $\O(\Gamma)$ by simplicial paths.
The set
$$\{X\in\overline{\O(\Gamma)}^\infty:\lambda_\phi(X)=\lambda(\phi)\}$$
is connected by simplicial paths in $\overline{\O(\Gamma)}^\infty$.
\end{theoremff}

This result is of independent interest, and surprisingly strong. We also show how to deduce
algorithmic results from this geometric one. We note that these algorithmic results are already
known, but are mainly a demonstration of the power of the result.  In general, if one wants to
check a property $P$, detectable by simplicial maps somewhere in outer space, an
algorithm of the type ``go to neighbour simplex and check $P$'' clearly stops if it finds a
points having $P$, but it has no a stopping procedure. Connectivity of level sets provides
 stopping criteria.

Specifically, we solve the conjugacy problem for irreducible automorphisms and prove that it is
determining whether an automorphism is irreducible or not is decidable. The conjugacy problem
for irreducible automorphisms has already been solved by \cite{MR1396778} and
\cite{MR2395795}. Deciding irreducibility of automorphisms has been proven by \cite{MR3194747}
and improved in \cite{kapovichpolyred}. While our solution of the former has some similarities
to that in \cite{MR1396778} (namely, peak reduction), our approach is distinctly geometric and
provides a uniform framework for dealing with this type of problems for general deformation spaces. Moreover, our connectivity result is general and does not assume irreducibility, so there is hope of pushing these techniques further even though one generally needs irreducibility in order to avoid singularities (which for us means entering the thin part of $CV_n$). 

\medskip

 In proving Theorem~\ref{tconnected} we
obtain a collection of results that may be of independent interest. Namely:
\begin{itemize}
  \item In Section~\ref{section_conv} we give a detailed analysis of the convexity properties of the displacement function.
  \item We prove that the global simplex-displacement spectrum of $\Aut(F_n)$ is well-ordered. (Theorem~\ref{conj}.)
  \item Generalising a result of~\cite{FM13}, we show that local minima of the displacement are
    global minima 
    (Lemma~\ref{Lemmatt4}). This allows, together with Theorem~\ref{conj} to implement an
    efficient gradient method for finding train tracks. 
  \item We study the behaviour of the displacement at bordification points, providing a
    characterisation of those points where the displacement does not jump. (Corollary~\ref{corlalx}.)
  \item We show that train tracks at infinity minimise the displacement. (Therem~\ref{corttE}.)   
  \item Given an automorphism $\phi$, we show that any invariant free factor is visible in
    a train track map. (Corollaries~\ref{corred} and~\ref{strongcorred})
  \item  We also wish to mention Theorem~\ref{thmZnew}: a technical result, with explicit
    estimates, that can be phrased as ``Folds of illegal turns in a simplex may be closely read in
    close simplices''.
\end{itemize}
The main results of the paper are proved by induction on the rank, and Theorem~\ref{tconnected} is
assumed inductively true in many points.
\medskip

The paper is extremely technical, even though the ideas are fundamentally straightforward. In order to motivate the detailed discussion, we provide here the two algorithms for solving conjugacy in the irreducible case and for detecting irreducibility. We present these algorithms as naively as possible, in order to make them more accessible. That is, one could understand and implement them without any knowledge of the Lipschitz metric, Culler-Vogtmann space or train track maps. As such we have made no attempt to streamline the algorithms in any way; they are brute force searches in an exponential space. 

However, we would stress that our point of view is fundamentally that these procedures would be better run as path searches in Culler-Vogtmann space, enumerating optimal $PL$-maps and calculating displacements via candidates. That abundance of terminology would make the algorithms much harder to describe, so we instead translate everything to a more manageable setting; bases of $F_n$ and generating sets for $\Out(F_n)$. However, the technical point of view is more helpful in developing an intuition of the processes and is likely the way to vastly improve the algorithmic complexity. 

\medskip

Let us know describe our algorithms, whose correctness is proved at the end of the paper. First, we recall some terminology. In order to work algorithmically with $\Out(F_n)$ we need a generating set. The best known of these is the set of Nielsen generators, but it is more convenient for us to work with the following:

\begin{defn}[CMT Automorphisms, \cite{MR721773} and \cite{MR748994}]
A {\em CMT} automorphism of $F_n$ is one that is induced by a change of maximal tree. More precisely, let $X$ be a graph with fundamental group of rank $n$, and let $R$ be the rose of rank $n$ (the graph with one vertex and $n$ edges). Let $T, T'$ be two maximal trees of $X$, and let $\rho_T, \rho_{T'}$ be the corresponding projections from $X$ to $R$. Then the (outer) automorphism induced by changing the maximal tree from $T$ to $T'$ is the (homotopy class of the) map $\rho_{T'} {\rho_T}^{-1}$, where the inverse denotes a homotopy inverse. 

The set of CMT maps includes all Whitehead automorphisms, (see \cite{MR721773}, Theorem~5.5 and \cite{MR1922276}) and is a finite set which generates $\Out(F_n)$. 

For convenience, we will include all graph automorphisms of $R$, including inversions of generators,  in the set of CMT automorphisms. 
\end{defn}

Next we need a notion of size of an automorphism, which will provide a termination criterion for our algorithms. 

\begin{defn}
Let $\phi \in \Out(F_n)$, and let $B$ be a basis of $F_n$. Define $|| \phi ||_B$ to be $\sup_{1
  \neq g \in F_n} \frac{|| \phi g||_B}{||g||_B}$, where $||g||_B$ denotes the cyclic reduced length of
$g$ with respect to $B$. This supremum is a maximum and is realised by an element of cyclic
length $\leq 2$.  
\end{defn}

\begin{rem}
Note that for any constant, $C$, there are only finitely many $\phi \in \Out(F_n)$ such that $||\phi||_B \leq C$.
\end{rem}

Our first application is then as follows. (See Section~\ref{appl} for the proof.)

\begin{thm} \label{conjirred}
The following is an algorithm to determine whether two irreducible automorphisms are conjugate. 

\medskip

\noindent
Let $\phi, \psi$ be two irreducible outer automorphisms of $F_n$, and $B$ a basis of $F_n$.
\begin{itemize}
\item Choose any $\mu > \max \{ ||\phi||_B, ||\psi||_B \}$. 
\item Inductively construct a finite set, $S=S_{\phi, \mu}$, as follows (which depends on both $\phi$ and $\mu$): 
\begin{itemize}
\item Start with $S_0 = \{ \phi \}$. 
\item Set $K= n(3n-3) \mu^{3n-1}$. 
\item Inductively put $S_{i+1}$ to be all possible automorphisms $ \zeta \phi_i \zeta^{-1}$, where $\phi_i$ is any element of $S_i$, $\zeta$ is any CMT automorphism, subject to the constraint that $||   \zeta \phi_i \zeta^{-1} ||_B \leq K$. (We include the identity as a CMT automorphism so that $S_{i-1} \subseteq S_i$). 
\item End this process when $S_i = S_{i+1}$, and let this final set be $S$. 
\end{itemize} 
\item Then $\psi$ is conjugate to $\phi$ if and only if $\psi \in S$. 
\end{itemize}
\end{thm}

\medskip

Of course, one would like to also be able to decide when an automorphism is irreducible when it is given by images of a basis, for instance. In order to do so, we recall the definition of irreducibility. 

\begin{defn}[see \cite{BestvinaHandel}]
An (outer) automorphism, $\psi$ of $F_n$ is called {\em reducible} if there are free factors, $F_1, \ldots, F_k, F_{\infty}$ such that $F_n = F_1 * \ldots F_k * F_{\infty}$ and each $\psi(F_i)$ is conjugate to $F_{i+1}$ (subscripts taken modulo $k$). If $k=1$ we further require that $F_{\infty} \neq 1$. (In general $\psi(F_{\infty})$ is not conjugate to $F_{\infty}$. Otherwise $\phi$ is called irreducible. 

Equivalently, $\psi$ is reducible if it is represented by a homotopy equivalence, $f$, on a core graph, $X$, such that $X$ has a proper, homotopically non-trivial subgraph, $X_0$, such that $f(X_0) = X_0$. (Being represented by $f$ means that there is an isomorphism, $\tau: F_n \to \pi_1(X)$ such that $\psi = \tau^{-1} f_* \tau$).  
\end{defn}

We add the following, which constitutes an obvious way that one can detect irreducibility by inspection.  

\begin{defn}
Consider $F_n$ with basis $B$ and let $\psi$ be an outer automorphism of $F_n$. We say that $\psi$ is visibly reducible with respect to $B$, or simply visibly reducible, if there exist disjoint subsets $B_1, \ldots, B_k$ of $B$ such that $\psi(\langle B_i \rangle)$ is conjugate to $\langle B_{i+1} \rangle$ (with subscripts taken modulo $k$). If $k=1$ we also require that $B_1 \neq B$. 

More generally, we say that a homotopy equivalence on the rose is visibly reducible if it is visibly reducible with respect to the basis given by the edges of the rose. 
\end{defn}

This is, in fact, easy to check by classical methods due to Stallings, \cite{Sta}.

\begin{lem}
If $\psi$ is visibly reducible, it is reducible. Moreover, there is an algorithm to determine if $\psi$ is visibly reducible with respect to $B$. 
\end{lem}
\proof
The first statement is clear, since each subset of a basis generates a free factor, and disjoint subsets generate complementary free factors. Since there are only finitely many subsets to check, we simply need to determine if the conditions that $\psi(\langle B_i \rangle)$ is conjugate to $\langle B_{i+1} \rangle$ hold. But this can readily be checked since two subgroups of a free group are conjugate if and only if the core of their Stallings graphs are equal, \cite{Sta}. 
\qed

We can now describe our second algorithm. (See Section~\ref{appl} for the proof.)

\begin{thm} \label{detectirred}
The following is an algorithm to determine whether or not an outer automorphism of $F_n$ is irreducible.

\medskip

\noindent
Let $\phi$ be an automorphism of $F_n$, and $B$ a basis of $F_n$. Construct $S=S_{\phi}$ as above. Namely, 
\begin{itemize}
\item Choose any $\mu > ||\phi||_B$. 
\item Inductively construct the finite set, $S=S_{\phi, \mu}$: 
\begin{itemize}
\item Start with $S_0 = \{ \phi \}$. 
\item Set $K= n(3n-3) \mu^{3n-1}$. 
\item Inductively put $S_{i+1}$ to be all possible automorphisms $ \zeta \phi_i \zeta^{-1}$, where $\phi_i$ is any element of $S_i$, $\zeta$ is any CMT automorphism, subject to the constraint that $||   \zeta \phi_i \zeta^{-1} ||_B \leq K$. (We include the identity as a CMT automorphism so that $S_{i-1} \subseteq S_i$). 
\item End this process when $S_i = S_{i+1}$, and let this final set be $S$. 
\end{itemize} 
\item Let $S^+$ be the set of all possible automorphisms $ \zeta \phi_i \zeta^{-1}$, where $\phi_i$ is any element of $S$, $\zeta$ is any CMT automorphism, with no other constraint. 
\item If some $\psi \in S^+$ is visibly reducible with respect to $B$, then $\phi$ is reducible. Otherwise, $\phi$ is irreducible. 
\end{itemize}
\end{thm}

\section{Preliminaries}
\subsection{Motivation for new definitions}
First, we want to motivate the definitions that we are going to
introduce. This is because they are a little different and at times more complicated
than those usually present in literature. Our aim is to study automorphisms of free groups
which are possibly reducible. If $\Gamma$ is a marked graph with $\pi_1(\Gamma)=F$ a free
group, and $\phi\in\Aut(F)$, then $\phi$ can be represented by a simplicial map (sending vertices to vertices and edges to edge paths) $f:\Gamma\to\Gamma$. That is, $f$ represents $\phi$ if there is an isomorphism $\tau: F_n \to \pi_1(\Gamma)$ such that $\phi = \tau^{-1} f_* \tau$. 

If $\phi$ is reducible it may happen that there is a collection of
disjoint connected sub-graphs $\Gamma_1,\dots,\Gamma_k$ such that $f$ permutes the $\Gamma_i$'s. In
order to study the properties of $\phi$ it may help to collapse such an invariant collection.
If we want to keep track of all the relevant information, we will be faced with the study of some
particular kind of moduli space. Namely, moduli spaces of actions on trees with possibly
non-trivial vertex stabilizers (when we collapse the $\Gamma_i$'s) and product of such spaces
(when we consider the restriction to $\phi$ to the $\Gamma_i$'s.)

Since the notation that we are going to introduce may be cumbersome, we will often abuse it, making no distinction between an element of outer spaces and its (projective) class, or
between $G$-trees and $G$-graphs.
\subsection{General definitions and notations}

Let $G=G_1*\dots* G_p* F_n$ be a free product of groups, where $F_n$ denotes the free group of
rank $n$ (we allow $n$ to be zero, in that case we omit $F_n$). We do not assume the $G_i$'s are
indecomposable. Throughout the paper, $G$ will be a free group. In particular, each
$G_i$ will be a free factor of $G$. Thus, there is no uniqueness of this free product decomposition, since $G$ has many different
splittings as a free product. We use the notation $\mathcal G: G=G_1*\dots*G_p*F_n$ to indicate
a splitting of $G$.
We briefly recall the definition of the outer space $\O(G)$
of $G$ corresponding to the splitting $\mathcal G$, referring to
\cite{FM13,GuirardelLevitt} for a detailed discussion of definitions and general properties of
$\O(G)$.

\begin{defn}[Outer space]
  The (projectivized) outer space of $G$, relative to the splitting $\mathcal G:G=G_1*\dots*G_p*F_n$, consists of (projective) classes of minimal simplicial metric  $G$-trees $X$ such that:
  \begin{itemize}
  \item For every $G_i$ there is exactly one orbit of vertices whose stabilizer is conjugate to
    $G_i$. Such vertices are called {\em non-free}. Remaining vertices have trivial
    stabilizer and are called {\em free} vertices.
  \item $X$ has no redundant vertex (i.e. free and two-valent).
  \item $X$ has trivial edge stabilizers.
  \end{itemize}
We use the notation $\O(G)$ to indicate the outer space of $G$ and, if we want to emphasize the
splitting, we write $\O(G;\mathcal G)$ (and $\p\O(G)$ and $\p\O(G;\mathcal G)$ for
projectivized ones). We stress here that when the distinction between $\O(G)$ and $\p\O(G)$ is
not crucial, we will often make no distinction between $\O(G)$ and $\p\O(G)$.
\end{defn}

\begin{rem}
The equivalence relation that defines $\p\O(G)$ is the following: $X$ and $Y$ are equivalent if
there is an homothety (isometry plus a rescaling by a positive number) $X\to Y$ conjugating the
actions of $G$ on $X$ and $Y$.  If there is no ambiguity, we will make no distinction between
a $G$-tree $X$, its class in $\O(G)$, and its projective class in $\p\O(G)$.
\end{rem}

\begin{rem}
  If $G=G_1$, then $\O(G)$ consists of a single element: a point stabilized by $G_1$, and in
  this case the equivalence relation is trivial.
\end{rem}

\begin{defn}
A splitting $\mathcal S: G=H_1*\dots...*H_q*F_r$ is a {\em sub-splitting} of
$G=G_1*\dots*G_p*F_n$ if each $H_i$ decomposes as $$H_i=G_{i_1}*\dots G_{i_{l}}*F_{s}$$ where $F_s$
is a free factor of $F_n$ and $i_1,\dots,i_l\in\{1,\dots,p\}$. Sometime we will make use of the
notation $\O(G;H_1*\dots*H_q*F_r)$ to mean $\O(G;\mathcal S)$.
\end{defn}

\begin{rem}
If $T\in\O(G)$, the quotient $X=G\backslash T$ is a finite metric graph of groups with trivial
edge-groups, together with a marking that identifies $G$ with $\pi_i(X)$ and maps the $G_i$'s
to the vertex-groups. We will refer to such graphs as $G$-graphs (or $(G,\mathcal G)$-graph if
we need to specify the splitting). On the other hand, given a metric $G$-graph $Y$, its universal
cover $\wt Y$ is a $G$-tree in the unprojectivized outer space $\O(G)$. Here $\pi_1(X)$ means
the fundamental group of $X$ as graph of groups. (The fundamental group of $X$ as a topological
space is just $F_n$.)
\end{rem}

\begin{nt}
  If there is no ambiguity we will make no distinction between $G$-tress and $G$-graphs. In case
  of necessity we will use the {\em tilde}-notation: $X$ for a $G$-graph and $\wt X$ for a
  $G$-tree, meaning that $X=G\backslash \wt X$. As usual, if $x\in X$ then $\wt x$ will denote
  a lift of $x$ in $\wt X$. The same for subsets: if $A\subset X$ then $\wt A\subset \wt X$ is
  one of its lifts.
\end{nt}

\begin{defn}[Immersed loops]
   A path $\gamma$ in a $G$-graph $X$ is called {\em immersed} if it is has a lift $\wt \gamma$
   in $\wt X$ which is embedded. (Note that
$\gamma$ could not be topologically immersed.)
\end{defn}

Let $X$ be a $G$-graph and let $\Gamma=\sqcup_i\Gamma_i$ be a sub-graph of $X$ whose connected
components $\Gamma_i$ have non-trivial fundamental groups (as graphs of groups).
Then $\Gamma$ induces a sub-splitting $\mathcal S$ of $\mathcal G$ where the factor-groups $H_j$ are either
\begin{itemize}
\item the fundamental groups $\pi_1(\Gamma_i)$, or
\item the vertex-groups of non-free vertices in $X\setminus \Gamma$.
\end{itemize}

\begin{nt}\label{not:Ogamma}
  We will use the notation
 $$\O(X):=\O(G)\qquad \O(X/\Gamma):=\O(G;\mathcal S)\qquad \O(\Gamma):=\Pi_i\O(\pi_1(\Gamma_i))$$
\end{nt}

The above notation leads to the following general definition.
\begin{defn}
  Let $\Gamma=\sqcup_{i=1}^k\Gamma_i$ be a finite disjoint union of finite connected graphs of
  groups $\Gamma_i$ with trivial edge-groups and non-trivial fundamental group (as graphs of
  groups). Let
  $H_i=\pi_1(\Gamma_i)$, equipped with the splitting given by the vertex-groups of $\Gamma_i$
  (hence $\Gamma_i$ is an $H_i$-graph). We define $\O(\Gamma)$ as the product of the
  $\O(H_i)$'s
$$\O(\Gamma)=\Pi_{i=1}^k\O(H_i)=\Pi_{i=1}^k\O(\Gamma_i).$$
We tacitly identify $X=(X_1,\dots,X_k)\in \O(\Gamma)$ with the labelled disjoint union
$X=\sqcup_i X_i$. An element of $\O(\Gamma)$ will be also called {\bf $\Gamma$-graph} (or
{\bf $\Gamma$-tree} if we work with universal covers).
\end{defn}

Here we need to be more precise about projectivization.
There is a natural action of $\R^+$ on $\O(\Gamma)$ given by scaling each component by the same
amount. The quotient of $\O(\Gamma)$ by such action is the projective outer space of $\Gamma$
and it is denoted by $\p\O(\Gamma)$.

\begin{nt}\label{not:gamma}
In what follows we use the following convention:
\begin{itemize}
\item $G$ will always mean a group with  a splitting $\mathcal G:G=G_1*\dots*G_p*F_n$;
\item  $\Gamma=\sqcup \Gamma_i$ will always mean that $\Gamma$ is a finite disjoint
union of finite graphs of groups $\Gamma_i$, each with trivial edge-groups and non-trivial fundamental group $H_i=\pi_i(\Gamma_i)$, each
$H_i$ being equipped with the splitting given by the vertex-groups.
\end{itemize}
\end{nt}

\begin{defn}\label{pr4_rank}
  The rank of the splitting $G=G_1*\dots*G_p*F_n$ is $n+p$. The rank of a graph of groups $X$ is
  the rank of the splitting induced on $\pi_1(X)$, finally if $\Gamma=\sqcup \Gamma_i$ we set $$\rank(\Gamma)=\sum_i\rank(\Gamma_i).$$
\end{defn}
By definition, the rank is a natural number greater or equal to one.
Note the the rank of a graph of groups $X$ is nothing but the rank of its fundamental group as
a topological space plus the number of non-free vertices.

We will also consider moduli spaces with marked points.
\begin{nt}
  The moduli space of $G$-trees with $k$ labelled points $p_1,\dots,p_k$ (not necessarily
  distinct) is denoted by $\O(G,k)$
  or $\O(G;\mathcal G,k)$. If $\Gamma=\sqcup_{i=1}^s \Gamma_i$, given $k_1,\dots,k_s$ we
  set $$\O(\Gamma,k_1,\dots,k_s)=\Pi_i\O(\Gamma_i,k_i).$$
If $X$ is a $\Gamma$-graph and $A\subset X$ is a subgraph whose components have non-trivial
fundamental group, we define $\O(X/A)$ and $\O(A)$ as in Notation~\ref{not:Ogamma}.
\end{nt}

\subsection{Simplicial structure}

\begin{defn}[Open simplices]
Given a $G$-tree $X$,  the open simplex $\Delta_X$ is the set of $G$-trees
equivariantly homeomorphic to $X$ . The Euclidean topology on $\Delta_X$ is given by assigning a
$G$-invariant positive lent $L_X(e)$ to each edge $e$ of $X$. Therefore, if $X$ has $k$ orbit of
edges, then $\Delta_X$ is isomorphic to the
standard open $(k-1)$-simplex if we work in $\p\O(G)$, and to the positive cone over it if we
work on $\O(G)$. Given two elements $X,Y$ in the same simplex
$\Delta\subset \O(G)$ we define the {\bf
  Euclidean} sup-distance $d_\Delta^{Euclid}(X,Y)$ ($d_\Delta(X,Y)$ for short)
$$d_\Delta^{Euclid}(X,Y)=d_\Delta(X,Y)=\max_{e\text{ edge}}|L_X(e)-L_Y(e)|.$$
\end{defn}

Such definitions extend to the case of $\Gamma=\sqcup_i\Gamma_i$.

\begin{defn}
    If $X=(X_1,\dots,X_k)\in\O(\Gamma)$, the simplex $\Delta_X$ is the set of $\Gamma$-trees
    equivariantly homeomorphic to $X$ (component by component).
The Euclidean topology and  distance on $\Delta_X$ are defined by $$d_\Delta(X,Y)=\sup_i d_{\Delta_{X_i}}(X_i,Y_i).$$
\end{defn}

We notice that the simplicial structure of $\p\O(\Gamma)$ is not the product of the structures of
$\p\O(\pi_1(\Gamma_i))$.

\begin{rem}
  If $X\in\O(G)$, then $\O(X)=\O(G)$. In other words, $\O(G)$ is a particular case of
  $\O(\Gamma)$ whit $\Gamma$ connected. In the following we will therefore develop the
  theory of $\O(\Gamma)$ and that of $\O(G)$ at once.
\end{rem}

\begin{defn}[Faces and closed simplices]\label{def:face}
Let $X$ be a $G$-graph (resp. a $\Gamma$-graph) and let $\Delta=\Delta_X$ be the corresponding
open simplex. Let $F\subset X$ be a forest whose trees each contains at most one non-free
vertex. The collapse of $F$ in $X$ produces a new $G$-graph (resp. $\Gamma$-graph), whence a
simplex $\Delta_F$. Such a simplex is called a {\em face} of $\Delta$.

  The {\em closed} simplex $\overline{\Delta}$ is defined by
$$\overline{\Delta}=\Delta\cup\{\text{all the faces of $\Delta$}\}.$$
\end{defn}

\subsection{Simplicial bordification}

  There are two natural topologies on $\O(G)$ (resp. $\O(\Gamma)$), the simplicial one and the equivariant Gromov
  topology, which are in general different. Here we will
  mainly use the simplicial topology. We notice that if $\Delta$ is an open simplex,
  the simplex $\overline\Delta$ is not the standard simplicial closure of $\Delta$, because  not all
  its simplicial faces are {\em faces} according to Definition~\ref{def:face}. This is because
  some simplicial face of $\Delta$ are not in $\O(G)$ (resp. $\O(\Gamma)$) as defined. Such faces are somehow ``at
  infinity'' and describe limit points of sequences in $\O(G)$ (resp. $\O(\Gamma)$).
We give now precise definitions
  to handle such limit points.

  We will sometimes refer to the faces of $\Delta$, as defined in Definition~\ref{def:face} as {\em finitary faces} of $\Delta$.
\begin{defn}
  Given an open simplex $\Delta$ in $\O(\Gamma)$, its {\em boundary at the finite} is the set of its
  proper faces: $$\partial_\O\Delta=\partial_\O\overline\Delta=\overline\Delta\setminus\Delta.$$
\end{defn}

\begin{defn}
  A {\em core-graph} is a connected graph of groups whose leaves (univalent vertices) have
  non-trivial vertex-group. Given a graph $X$ we define $\core(X)$ to be the maximal core
  sub-graph of $X$. (If the vertex groups are all trivial, so that $X$ is simply a graph, then a core graph has no valence one vertices).
\end{defn}
Note that $\core(X)$ is obtained by recursively cutting edges ending at leaves.

\medskip

Let $X$ be a $\Gamma$-graph and $\Delta=\Delta_X$. Let $A$ be a proper subgraph of $X$ having
at least a component which is not a tree with at most one non-free vertex. Let $Y$ be the graph
of groups obtained by collapsing each component of $A$ to a point (different components to
different points). Then, $Y\in \O(X/A)$. The corresponding simplex $\Delta_Y$ is a simplicial
face of $\Delta_X$ obtained by setting to zero the edge-lengths of $A$.
\begin{defn}
  A face $\Delta_Y$ obtained as just described is called a {\em face at infinity} of
  $\overline\Delta_X$. If in addition we have that all components of $A$ are core-graphs, then
  we say that $\Delta_Y$ is a face at infinity of $\Delta_X$.

We define the {\em boundaries at infinity} by
$$\partial_\infty\Delta=\{\text{faces at infinity of } \Delta\}$$
$$\partial_\infty\overline\Delta=\{\text{faces at infinity of } \overline\Delta\},$$

and the {\em closure at infinity by}
$$\overline\Delta^\infty=\overline\Delta\cup\partial_\infty\overline\Delta.$$
\end{defn}

If we denote by $\partial \Delta$ the simplicial boundary of $\Delta$, we have
$$\partial \Delta=\partial_\infty\overline \Delta\cup\partial_\O\overline\Delta$$
and $$\partial_\infty \overline\Delta=\bigcup_{F=\text{face of }\Delta}\partial_\infty F$$
(where the union is over all faces of $\Delta$, $\Delta$ included.) Moreover, the simplicial
closure of $\Delta$ is just $\overline\Delta^\infty$.

\begin{defn}
  We define the boundary at infinity and the simplicial bordification of $\O(\Gamma)$ as
$$\partial_\infty\O(\Gamma)=\bigcup_{\Delta\text{ simplex}}\partial_\infty\Delta
\qquad \text{and}\qquad \overline{\O(\Gamma)}=\overline{\O(\Gamma)}^\infty=\O(\Gamma)\cup\partial_\infty\O(\Gamma).$$
\end{defn}

\subsection{Horospheres and regeneration}\label{sechor}
\begin{defn}
  Given $X\in\partial_\infty\O(\Gamma)$, the horosphere $\Hor(\Delta_X)$ of $\Delta_X$ in 
$\O(\Gamma)$ is the union of simplices $\Delta$ such that $X\in\partial_\infty \Delta$.
If $X\in\O(\Gamma)$ we  set  $\Hor(\Delta_X)=\Delta_X$. 

The horosphere $\Hor(X)$ of $X$ in $\O(\Gamma)$ is the set formed by points $Y\in\Hor(\Delta_X)$ such that
$L_Y(e)=L_X(e)$ for any edge $e$ of $X$. (In particular, if $X\in\O(\Gamma)$ we
  have  $\Hor(X)=X$.)
\end{defn}
Thus, $Y$ is in the horosphere of $X$ if $X$ is obtained from $Y$ by collapsing a proper family
of core sub-graphs.  On the other hand, $\Hor(X)$ can be regenerated from $X$ as follows.

Suppose $X\in\partial_\infty\O(\Gamma)$. Thus there is a $\Gamma$-graph $Y$ and a sub-graph
$A=\sqcup_iA_i\subset Y$ whose components $A_i$ are core-graphs, and such that $X=Y/A$. Let
$v_i$ be the non-free vertex of $X$ corresponding to $A_i$. In order to recover a generic point
$Z\in\Hor(X)$, we need to replace each $v_i$ with an element $V_i\in\O(A_i)$. Moreover,
in order to define the marking on $Z$, we need to know where to attach to $V_i$ the edges of $X$
incident to $v_i$, and this choice has to be done in the universal covers $\wt{V_i}$. No more
is needed. Therefore, if $k_i$ denote the valence of the vertex $v_i$ in $X$, we have
$$\Hor(X)=\Pi_i\O(A_i,k_i).$$
(Note that some $k_i$ could be zero, e.g.  if $A_i$ is a connected component of $Y$.)
There is a natural projection $\Hor(X)\to\O(A)$ which forgets the marking. We will be mainly
interested in cases when we collapse $A$ uniformly, for that reason we will use the
projection $$\pi:\Hor(X)\to\mathbb P\O(A)$$
where $\Hor(X)$ is intended to be not projectivized.

Note that if $[P]\in\mathbb P\O(A)$, then $\pi^{-1}(P)$ is connected because it is just
$\Pi_{i}(A_i^{k_i})$. Since $\O(A)$ is connected, then $\Hor(X)$ is connected.

Finally, we notice that a graph of groups $X$ can be considered as a point at infinity of
different spaces. If we need
to specify in which space we work we write $\Hor_\Gamma(X)$ or $\Hor_G(X)$.

\subsection{The groups $\Aut(\Gamma)$ and $\Out(\Gamma)$}

\begin{defn}
  Let $G$ be endowed with the splitting $\G:G=G_1*\dots*G_p*F_n$.
The group of automorphisms of $G$
  that preserve the set of conjugacy classes of the $G_i$'s is
  denoted by $\operatorname{Aut}(G;\G)$. We set
  $\operatorname{Out}(G;\G)=\operatorname{Aut}(G;\G)/\operatorname{Inn}(G)$\footnote{Clearly $\operatorname{Inn}(G)\subset\Aut(G;\G)$.}.
\end{defn}


The group $\Aut(G,\G)$ acts on $\O(G)$ by
changing the marking (i.e. the action), and $\operatorname{Inn}(G)$ acts trivially. Hence
$\operatorname{Out}(G;\G)$ acts on
$\O(G)$. If $X\in\O(G)$ and
$\phi\in\operatorname{Out}(G;\G)$ then $\phi X$ is the same metric tree as $X$, but the action
is $(g,x)\to \phi(g)x$. The action is simplicial and continuous w.r.t. both simplicial and
equivariant Gromov topologies.

We now extend the definition of $\Aut(G,\G)$ to the case of $\Gamma=\sqcup_i\Gamma_i$. Let
$\mathfrak S_k$ denotes the group of permutations of $k$ elements.
\begin{defn}
  Let $G$ and $H$ be two isomorphic groups endowed with splitting $\G:G=G_1*\dots G_p*F_n$ and
$\mathcal H:H=H_1*\dots H_p*F_n$. The set of isomorphisms from $G$ to $H$ that maps each $G_i$
to a conjugate of one of the $H_i$ is denoted by $\operatorname{Isom}(G,H)$. If we need to
specify the splittings we write $\operatorname{Isom}(G,H;\G,\mathcal H)$.

\end{defn}

\begin{defn} For $\Gamma=\sqcup_{i=1}^k\Gamma_i$ as in Notation~\ref{not:gamma},
  we set $$\Aut(\Gamma)=\{\phi=(\sigma,\phi_1,\dots,\phi_k):\ \sigma\in\mathfrak S_k\text{
and } \phi_i\in\operatorname{Isom}(H_i,H_{\sigma_i})\}.$$
\end{defn}

The composition of $\Aut(\Gamma)$ is component by component defined as follows. Given
$\phi=(\sigma,\phi_1,\dots,\phi_k)$ and $\psi=(\tau,\psi_1,\dots,\psi_k)$ we have
$$\psi\phi=(\tau\sigma,\psi_{\sigma(1)}\phi_1,\dots,\psi_{\sigma(k)}\phi_k)$$
\begin{rem}
  Not all permutations appear. For instance, if the groups $H_i$ are mutually not isomorphic,
  then the only possible $\sigma$ is the identity.
\end{rem}

\begin{defn}
  We set:
$$\operatorname{Inn}(\Gamma)=\{(\sigma,\phi_1,\dots,\phi_k)\in\Aut(\Gamma): \sigma=id,
\phi_i\in\operatorname{Inn}(H_i)\}$$
$$\Out(\Gamma)=\Aut(\Gamma)/\operatorname{Inn}(\Gamma).$$
\end{defn}

\begin{ex}
If $X\in\O(G)$ and $f:X\to X$ is a homotopy equivalence which leaves invariant a subgraph $A$,
then $f|_A$ induces and element of $\Aut(A)$, and its free homotopy class an element of $\Out(A)$.
\end{ex}

The group $\Out(\Gamma)$ acts on $\O(\Gamma)$ as follows. If $X=(X_1,\dots,X_k)\in\O(\Gamma)$,
then each $X_i$ is an $H_i$-tree. If $(\sigma,\phi_1,\dots,\phi_k)\in\Aut(\Gamma)$ then
$X_{\sigma(i)}$ becomes an $H_i$-tree via the pre-composition of $\phi_i:H_i\to H_{\sigma(i)}$
with the $H_{\sigma(i)}$-action. We denote such an $H_{i}$-tree by
$\phi_iX_{\sigma(i)}$. With that notation we have
$\phi(X_i,\dots,X_n)=(\phi_{1}X_{\sigma(1)},\dots,\phi_{k}X_{\sigma(k)}).$
(We remark that despite the left-positional notation, this is a right-action.)

\section{$\PL$-maps, gate structures, and optimal maps.}
In this section we describe the theory of maps between graphs (or trees) representing points in outer
spaces. We will treat in parallel the ``connected'' case $\O(G)$ and the general case
$\O(\Gamma)$, where $G$ and $\Gamma$ are as in Notation~\ref{not:gamma}.

\subsection{$\PL$-maps} Now we will mainly work with trees.

\begin{defn}[$\O$-maps in $\O(G)$]
Let $X,Y\in\O(G)$. A map $f:X\to Y$ is called an $\O$-map if it is Lipschitz-continuous and
$G$-equivariant. The Lipschitz constant of $f$ is denoted by $\Lip(f)$.
\end{defn}

We recall that we tacitly identify $X=(X_1,\dots,X_k)\in\O(\Gamma)$ with the labelled disjoint
union $\sqcup_i X_i$. Hence, if $X,Y\in\O(\Gamma)$, a continuous map $f:X\to Y$ is a collection
of continuous maps $f_i;X_i\to Y_{j}$ for some $j=j(i)$.

\begin{defn}[$\O$-maps in $\O(\Gamma)$]
Let $X=(X_1,\dots,X_k)$ and $Y=(Y_1,\dots,Y_k)$ be two elements of $\O(\Gamma)$.
A map $f=(f_1,\dots,f_k):X\to Y$ is called an $\O$-map if for each $i$ the map
 $f_i$ is an $\O$-map from $X_i$ to $Y_i$.
\end{defn}

\begin{defn}[$\PL$-maps]
Let $X,Y$ be two metric trees.
A Lipschitz-continuous  map $f:X\to Y$ is a
$\PL$-map if it has constant speed on edges, that is to say, for any edge $e$
of $X$ there is a non-negative number $\lambda_e(f)$ such that for any $a,b\in e$ we have
$d_Y(f(a),f(b))=\lambda_e(f)d_X(a,b)$. If $X,Y\in\O(G)$ then we require any $\PL$-map to be an
$\O$-map. A $\PL$-map between elements of $\O(\Gamma)$ is an $\O$-map whose components are
$\PL$. (If $X,Y$ are $\Gamma$-graph, we understand that $f:X\to Y$ is a $\PL$-map if its lift to the universal
covers is $\PL$:)

\end{defn}

\begin{rem}\label{rem:24} $\O$-map always exists and the images of non-free vertices is
  determined a priori by  equivariance
  (see~\cite{FM13}).
  For any $\O$-map $f$ there is a unique $\PL$-map, denoted by $\PL(f)$, that conincides with $f$
  on vertices. We have $\Lip(\PL(f))\leq \Lip(f)$.
\end{rem}

\begin{defn}[$\lambda_{\max}$ and tension graph]
  Let $f:X\to Y$ be a $\PL$-map. We set $$\lambda(f)=\lambda_{\max}(f)=\max_{e}\lambda_e(f)=\Lip(f).$$  We
  define the {\em tension  graph} of $f$ by
$$X_{\max}(f)=\{e \text{ edge of } X : \lambda_e(f)=\lambda_{\max}\}.$$
If there are no ambiguities we set $\lambda_{\max}=\lambda_{\max}(f)$ and $X_{\max}=X_{\max}(f)$.
\end{defn}

\begin{defn}[Stretching factors]
For $X,Y\in\O(\Gamma)$ we define
$$\Lambda(X,Y)=\min_{f:X\to Y\ \O\text{-map}}\Lip(f)$$
\end{defn}

The theory of stretching factors is well-developed in the connected case (i.e. for $\O(G)$),
but one can readily see that connectedness of trees plays no role, and the theory extends without
modifications to the non-connected case. In fact,
 $\Lambda$ is well-defined, (see~\cite{FM11,FM13}for details) and it satisfies the
 multiplicative triangular inequality: $$\Lambda(X,Z)\leq\Lambda(X,Y)\Lambda(Y,Z)$$
It can be used to define a non-symmetric metric $d_R(X,Y)=\log(\Lambda(X,Y))$ and its
symmetrized version $d_R(X,Y)+d_R(Y,X)$ (see~\cite{FM11,FM12,FM13} for details) which
induces the Gromov topology. The group $\operatorname{Out}(\Gamma)$ acts by isometries on
$\O(\Gamma)$.

Moreover, there is an effective way to compute $\Lambda$, via the so-called
``sausage-lemma'' (see~\cite[Lemma 3.14]{FM11},\cite[Lemma 2.16]{FM12} for the classical case,
and~\cite[Theorem 9.10]{FM13} for the case of trees with non-trivial vertex-groups). We briefly
recall here how it works.

Let
$X,Y\in\O(\Gamma)$ (now seen as graphs).
Any non-elliptic element $\gamma\in\pi_1(\Gamma)$ (i.e. an element not in a
vertex-group) is represented by an immersed loop $\gamma_X$ in $X$ and one $\gamma_Y$ in
$Y$. The loop $\gamma_X$ (or its lift to $\wt X$) is usually called {\bf axis} of $\gamma$ in $X$
(or in $\wt X$).
They have lengths $L_X(\gamma_X)$ and $L_Y(\gamma_Y)$ that correspond to the minimal
translation length of the element $\gamma$ acting on $X$ and $Y$. (So
$L_X(\gamma_X)=L_X(\gamma)$ and $L_Y(\gamma_Y))=L_Y(\gamma)$.)  We can define the stretching
factor $L_Y(\gamma)/L_X(\gamma)$. Then $\Lambda(X,Y)$ is the minimum of the stretching factors
of all non-elliptic elements.

\begin{thm}[Saussage Lemma~{\cite[Theorem 9.10]{FM13}}]\label{sausagelemma}
  Let $X,Y,\in\O(\Gamma)$. The stretching factor $\Lambda(X,Y)$ is realized by an element
  $\gamma$ whose axis $\gamma_X$ has one of the following forms:
  \begin{itemize}
  \item Embedded simple loop $O$;
  \item embedded ``infinity''-loop $\infty$;
  \item embedded barbel $O$--- $O$;
  \item singly degenerate barbel $\bullet$---$O$;
  \item doubly degenerate barbel $\bullet$---$\bullet$.
  \end{itemize}
(the $\bullet$ stands for a non-free vertex.) Such loops are usually named ``candidates''.
\end{thm}

\begin{rem}
  The stretching factor $\Lambda(X,Y)$ is defined on $\O(\Gamma)$ and not in
  $\p\O(\Gamma)$. However, we will mainly interested in computing factors of type
  $\Lambda(X,\phi X)$ (for $\phi\in\Out(\Gamma)$) and that factor is scale invariant.
\end{rem}

\begin{defn}[Gate structures] Let $X$ be any graph.
  A {\em gate} structure on $X$ is an equivalence relation on germs of
  edges at vertices of $X$. Equivalence classes of germs are called {\em gates}.
A {\em train-track} structure on $X$ is a gate structure having at least two gates at every
vertex. A {\em turn} is a pair of germs of edges incident to the same vertex. A turn is {\em
  illegal} if the two germs are in the same gate, it is {\em legal} otherwise. An immersed path
in $X$ is legal if it has only legal turns.

If $X=(X_1,\dots,X_k)\in\O(\Gamma)$ we require the equivalence relation to be $H_i$-invariant
on each $X_i$.
\end{defn}

Any $\PL$-map induces a gate structure as follows.
\begin{defn}[Gate structure induced by $f$]\label{pr4_gate} Given $X,Y\in\O(\Gamma)$ and a $\PL$-map $f:X\to Y$,
  the gate structure induced by $f$, denoted by $$\sim_f$$
  is defined by declaring equivalent two germs that have the
  same non-degenerate $f$-image.
\end{defn}

\begin{rem}[See~\cite{FM13}]\label{rem:2gated}
Given $X,Y\in\O(\Gamma)$ and $f:X\to Y$ a $\PL$-map. If $v$ is a non-free vertex of $X$ and $e$
if an edge incident to $v$ which is not collapsed
by $f$, then $e$ and $ge$ are in different gates for any $id\neq g\in\operatorname{Stab}(v)$.
\end{rem}

\begin{defn}[Optimal maps] Given $X,Y\in\O(\Gamma)$,
  a map $f:X\to Y$ is {\em weakly  optimal} if it is $\PL$ and $\lambda(f)=\Lambda(X,Y)$.

 A map $f:X\to Y$ is {\em optimal} if the restriction of the gate
  structure induced by $f$, to the tension graph,  is a train track structure (in other words,
  if the vertices of
  $X_{\max}$ are at least two-gated in $X_{\max}$).
\end{defn}

\begin{rem} Optimal maps always exist and are weakly optimal.
  A map between two $\Gamma$-trees is weakly optimal if and only there is an periodic immersed
  legal line in the tension graph (i.e. a legal immersed loop in the quotient graph).
\end{rem}

In general optimal maps are neither unique nor do they form a discrete set, even if $X_{\max}=X$,
as the following example shows. (If $X_{\max}\neq X$ then one can use freedom given by the lengths of
edges not in $X_{\max}$ to produce examples.)

\begin{ex}[A continuous family of optimal maps with $X_{\max}=X$]

Consider $G=F_2$. Let $X$ be a graph with three edges $e_1,e_2,e_3$ and two free vertices
$P,Q$, as in Figure~\ref{fig:ex3.14}. Set the length of $e_2$ to be $2$, name $x$ the length of
$e_1$, and $1+\delta$ that of $e_3$. The parameters $x,\delta$ will be determined below. For
any $t\in [0,1]$ consider the point $P_t$ at distance $1+t$ from $P$ along $e_2$, and the
point $Q_t$  at distance $1-t$ from $P$ along $e_3$. $P_t$ divides $e_2$ in oriented segments
$a_t,c_t$. $Q_t$ divides $e_3$ into $b_t,d_t$.
\begin{figure}[htbp]
  \centering
  \begin{tikzpicture}[x=5ex,y=5ex]
    \draw (0,0) circle [radius =1.5];
    \draw (0,1.5) arc (30:330:3) ;
    \draw [arrows=-{latex}] (0,1.5) arc (30:50:3) ; 
    \draw [arrows=-{latex}] (0,1.5) arc (90:130:1.5) ; 
    \draw [arrows=-{latex}] (0,1.5) arc (90:50:1.5) ; 
    \filldraw (0,1.5) circle(2pt) node[above right] {$P$};
    \filldraw (0,-1.5) circle(2pt) node[below right] {$Q$};
    \filldraw (-1.5,0) circle(2pt) node[left] {$P_t$};
    \filldraw (1.5,0) circle(2pt) node[right] {$Q_t$};
    \draw (-1.1,1.1) node[below right] {$a_t$};
    \draw (1.1,1.1) node[below left] {$b_t$};
    \draw (-1.1,-1.1) node[above right] {$c_t$};
    \draw (1.1,-1.1) node[above left] {$d_t$};
    \draw (-5,2) node[left] {$e_1$};
    \draw (3.5,-3.5) node[above right] {\parbox{20ex}{$e_2=a_tc_t$ \\ $e_3=b_t d_t$ \\ $f_t(e_1)=\bar{a_t}e_1\bar
        c_t\bar a_t b_t$\\ $f_t(e_2)=c_t \bar d_t$\\ $f_t(e_3)=c_t \bar d_t \bar b_t a_t c_t \bar d_t$\\
        Length$(a_t)=1+t$\\ Length$(b_t)=1-t$\\ Length$(c_t)=1-t$\\Length$(d_t)=\delta+t$\\
        Length$(e_1)=x$\\Length$(e_2)=2$\\Length$(e_3)=1+\delta$}};
    \draw[red, dashed] (-1.5,0) to [out=90,in=190] (-.3,1.6);
    \draw[red, dashed] (-.3,1.6) arc (40:330:2.7);
    \draw[red, dashed] (0,-1.5) arc (-95:-345:1.4);
    \draw[blue] (-1.5,0) arc(-182:72:1.6);
    \draw[blue] (1.5,0) arc(-10:-290:1.3);
  \end{tikzpicture}
  \caption{A continuous family of optimal maps with $X_{\max} = X$. The red dashed line is
    $f(e_1)$ and the blue line is $f(e_3)$ ($f(e_2)$ is not depicted).}
  \label{fig:ex3.14}
\end{figure}
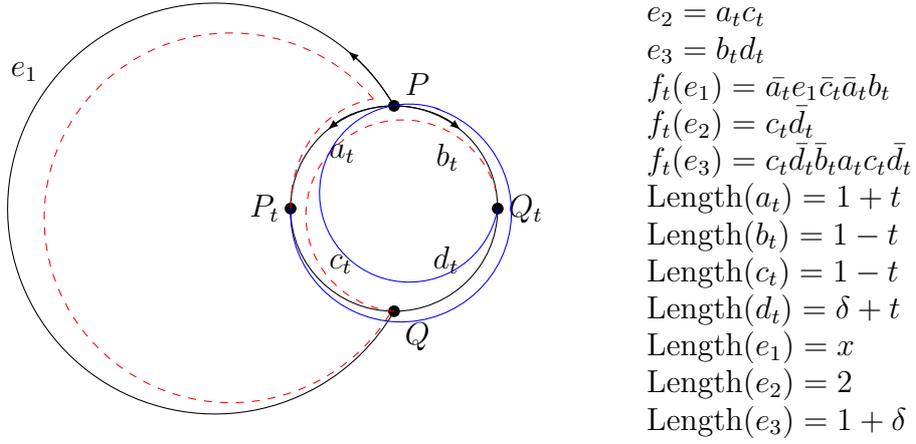 
Consider the $\PL$-map $f:X\to X$ defined as in the Figure, sending $P$ to $P_t$ and $Q$ to $Q_t$. If we collapse $e_3$, and we
homotop $P_t$ to $P$ along $a$, this corresponds to the automorphism $e_1\mapsto
e_1\overline{e_2}, e_2\mapsto \overline{e_2}$. 

A direct calculation shows that if we set $\delta=1+2\sqrt2$ and $x=2\sqrt2$, the map $f_t$ is
optimal for any $t$ and all the three edges are stretched by the same amount, as follows. 

\medskip

The edges
 $e_1$ and $e_2$ are in different gates at $P$ and $e_1$ and $e_3$ are in different gates at
 $Q$. In order to check that $f_t$ is optimal it suffices to check that every edge is stretched
 the same. 
$$\lambda_{e_1}(f_t)=\frac{x+4}{x}\qquad \lambda_{e_2}(f_t)=\frac{1+\delta}{2}\qquad\lambda_{e_3}(f_t)=\frac{4+2\delta}{1+\delta}.$$
In particular they do not depend on $t$. If we set $x=2\sqrt2$ and $\delta=1+2\sqrt2$ we get
$$\lambda_{e_1}(f_t)=\frac{2\sqrt2+4}{2\sqrt2}\qquad \lambda_{e_2}(f_t)=\frac{2+2\sqrt2}{2}\qquad\lambda_{e_3}(f_t)=\frac{6+4\sqrt2}{2+2\sqrt2}$$
and all of them are $1+\sqrt2$. \qed
\end{ex}

However, given a $\PL$-map, we can choose an optimal map which is in some sense the closest
possible. Given two $\O$-maps $f,g:X\to Y$ we define $$d_\infty(f,g)=\max_{x\in X} d_Y(f(x),g(x)).$$

For $X\in\O(G)$ we define its (co-)volume $\vol(X)$ as the sum of lengths of edges in
$G\backslash X$. If
$X=(X_1,\dots,X_k)\in\O(\Gamma)$ we set $\vol(X)=\sum_i\vol(X_i)$.

\begin{thm}[Optimization]\label{Lemma_opt}
Let $X,Y\in\O(\Gamma)$ and let $f:X\to Y$ be a $\PL$-map. There is a map\footnote{We describe
  an algorithm to fine the map $\wopt(f)$, but the algorithm will depend on some choice, hence
  the map $\wopt(f)$ may be not unique in general.} $\wopt(f):X\to Y$
which is weakly optimal and such that
$$d_\infty(f,\wopt(f))\leq \vol(X)(\lambda(f)-\Lambda(X,Y))$$
Moreover, for any $\varepsilon >0$ there is an optimal map $g:X\to Y$ such that
$d_\infty(g,\wopt(f))<\varepsilon$.
\end{thm}

\proof By arguing component by component, we may assume without loss of generality that
$\Gamma$ is connected, hence that we are in the case $X,Y\in\O(G)$.
For this proof it will be convenient to work with graphs rather than with trees (so
$X=G\backslash \wt X$, whit $\wt X$ a $G$-tree). By Remark~\ref{rem:2gated} a non-free vertex
will never be considered one-gated (because it is never one-gated in $\wt X$).

Let us concentrate on the first claim.

Let $\lambda=\Lambda(X,Y)$. Since $\PL$-maps are uniquely
determined by their value on vertices, we need only to define $\wopt(f)$ (and $g$) on vertices of
$X$. By Remark~\ref{rem:24} the image of non-free vertices is fixed.
We define $\PL$-maps $f_t$ for $t\in[0,\lambda_f-\lambda]$ by moving the images of all
one-gated vertices of $X_{\max}(f_t)$, in the direction given by the gate, so that $$\frac{d}{dt}\lambda(f_t)=-1.$$

Let us be more precise on this point. We define a flow which is piecewise linear, depending on
the geometry of the tension graph at time $t$.
The key remark to have in mind is that if an edge is not in
$X_{\max}(f)$, then it remains in the complement of the tension graph for small perturbations
of $f$. Therefore, we can restrict our attention to the tension graph.

Suppose we are at time $t$. We recursively define sets of vertices and edges as follows:

\begin{itemize}
\item $V_0$ is the set vertices of $X_{\max}(f_t)$ which are one-gated in
  $X_{\max}(f_t)$;
\item $E_0$ is the set of edges of $X_{\max}(f_t)$ incident to vertices in $V_0$. We agree that
  such edges contain the vertices in $V_0$ but not others. (If and edge has both vertices in
  $V_0$ then it contains both, otherwise it contains only one of its vertices.)
\end{itemize}
Having defined $V_0,\dots,V_i$ and $E_0,\dots,E_i$ , we define $V_{i+1}$ and $E_{i+1}$ as follows:
\begin{itemize}
\item $V_{i+1}$ is the set of
  one-gated vertices of $X_{\max}(f_t)\setminus\cup_{j=0}^iE_j$;
\item $E_{i+1}$ is the set of edges of $X_{\max}(f_t)\setminus\cup_{i=0}^iE_i$ incident to vertices in
  $V_{i+1}$. (As above such edges contain vertices in $V_{i+1}$ but not others.)
\end{itemize}
We notice  that since $X$ is a finite $G$-graph, we have only finitely many sets
$V_i$, say $V_0,\dots,V_k$.
\begin{lem}
  If $f_t$ is not weakly optimal, then $X_{\max}(f_t)\setminus\cup_{i=0}^kE_i$ is a (possibly empty)
collection of vertices, that we name terminal vertices.
\end{lem}
\proof
Note that no vertex in
$X_{\max}(f_t)\setminus\cup_{i=0}^kE_i$ can be one-gated, hence any vertex in
$X_{\max}(f_t)\setminus\cup_{i=0}^kE_i$ is either isolated or has at least two gates in
$X_{\max}(f_t)\setminus\cup_{i=0}^kE_i$. Thus if there is an edge $e$ in $X_{\max}(f_t)\setminus\cup_{i=0}^kE_i$,
the component of $X_{\max}(f_t)\setminus\cup_{i=0}^kE_i$ containing $e$ must also contain an
immersed legal loop and so $f_t$ is weakly optimal.\qed

\medskip

By convention we denote the set of terminal
vertices  by $V_\infty$.

\begin{rem}\label{rem:vivj}
Any $e\in E_i$ has by definition at least one endpoint in
$V_i$, and  the
other endpoint is in some $V_j$ with $j\geq i$.
\end{rem}
Our flow is defined by moving the images $f_t(v)$ of vertices in $X_{\max}(f_t)$. We need to define a
direction and a speed $s(v)\geq 0$ for any $f_t(v)$.

For $i<\infty$ each vertex in $V_{i}$ has a preferred gate: the one that survives in
$X_{\max}(f_t)\setminus\cup_{j=0}^{i-1}E_j$. That gate gives us the direction in which we move $f_t(v)$.

Thus a vertex in $V_0$ is one-gated, and hence we define the flow so as to reduce the Lipschitz constant for every edge in $E_0$ (shrinking the image of each $E_0$ edge). Similarly, every vertex in $V_1$ is one gated in $X_{\max}(f_t)\setminus E_0$, so we define the flow to reduce the Lipschitz constants of edges in $E_1$ and so on. 

Now we define the speeds.

\begin{lem}
  There exists speeds $s(v)\geq 0$ such that if we move the
  images of any $v$ at speed $s(v)$ in the direction of its preferred gate, then
  for any edge $e\in X_{\max}(t)$ $$\frac{d}{dt}\lambda_e(f_t)\leq -1.$$ Moreover,
  for any $i$, and for any $v\in V_i$, either $s(v)=0$ or there is an edge $e\in E_i$ incident at
  $v$ such that $$\frac{d}{dt}\lambda_e(f_t)=-1.$$
\end{lem}
\proof
We start by giving a total order of the vertices of $X_{\max}(f_t)$ in such a way
that vertices  in $V_i$ are bigger than those in $V_j$ whenever $i>j$. We define the speeds
recursively.

The speed of terminal vertices is set to zero. Let $v$ be a
vertex of $X_{\max}(f_t)$ and suppose that we already defined the speed $s(w)$ for all
$w>v$.

The vertex $v$ belongs to some set $V_i$. For any edge $e\in E_i$ emanating from $v$ let
$u_e$ be the other endpoint of of $e$, and define a sign $\sigma_e(u_e)=\pm1$ as follows:
$\sigma_e(u_e)=-1$ if the germ of $e$ at $u_e$ is in the preferred gate of $u_e$, and
$\sigma_e(u_e)=1$ otherwise. (So, for example, $\sigma_e(u_e)=1$ if $u_e$ is terminal, and
$\sigma_e(u_e)=-1$ if $v=u_e$, or if $u_e\in V_i$.)

With this notation, if we move $f(v)$ and $f(u_e)$ in the direction given by their gates at
speeds $s(v)$ and $\nu$ respectively, then the derivative of $\lambda_e(f_t)$ is given by
$$-\frac{s(v)-\sigma_e(u_e)\nu}{L_X(e)}$$

If $u_e>v$ we already defined its speed. We set
$$s(v)=\max\{0,\max_{u_e>v}\{L_X(e)+\sigma_e(u_e)s(u_e)\},\max_{u_e=v}\frac{L_X(e)}{2}\}$$

where the maxima are taken over all edges $e\in E_i$ emanating from $v$. Note that there may
exist some such edge with $u_e<v$. By Remark~\ref{rem:vivj} in this case $u_e\in V_i$ (same $i$ as $v$), $\sigma_e(u_e)=-1$ and the
derivative of $\lambda_e$ will be settled later, when defining the speed of $u_e$.

With the speeds defined in that way we are sure that for any edge $e$ we have $d/dt
\lambda_e(f_t)\leq -1$ and, if $s(v)\neq 0$, then the edges that realize the above maximum
satisfy $d/dt \lambda_e(f_t)=-1$.
\qed

\medskip
Locally, when we start moving, the tension graph may lose some edges. However, the above lemma
assures that any vertex $v$ with $s(v)\neq 0$ is incident to an edge $e$ which is maximally
stretched and $d/dt\lambda_e=-1$. Hence such an edge remains in the tension graph when we start
moving. Since $d/dt\lambda_e\leq -1$ for any edge, it follows that when we start moving, the
tension graph stabilizes. So our flow is well defined in $[t,t+\epsilon]$ for some
$\epsilon>0$. If at a time $t_1>t$ some edge that was not previously in $X_{\max}(f_t)$ becomes
maximally stretched we recompute speeds and we start again. A priori we may have to recompute
speeds infinitely many times $t<t_1<t_2<\dots$ but the control on $d/dt\lambda(f_t)$
assures that $\sup t_i=T\leq \lambda_f-\lambda$. Since speeds are bounded the flow has a limit
for $t\to T$ an then we can restart from $T$. Therefore the set of times
$s\in[0,\lambda_f-\lambda]$ for which the flow  is well-defined for $t\in[0,s]$ is closed and
open and thus is the whole $[0,\lambda_f-\lambda]$.

With these speeds, we have $d/dt(\lambda(f_t))=-1$. Therefore
for $t=\lambda(f)-\lambda$, and not before, we have $\lambda(f_t)=\lambda$ hence
$f_t$ is weakly optimal. We define
$$\wopt(f)=f_{\lambda(f)-\lambda}.$$

Now we estimate $d_\infty(f,f_t)$. The $d_\infty$-distance between $\PL$-maps is bounded by the $d_\infty$-distance of their restriction to vertices.

We first estimate the speed at which the images of
vertices move. Let $S$ be the maximum speed of vertices,
i.e. $S=\max_v|s(v)|$. Let $v$ be a fastest vertex. Since it moves, it belongs
to $V_s$ for some $s<\infty$.
Let $v=v_1,v_2\dots,v_m$ be a maximal sequence of vertices such that:
\begin{enumerate}
\item $s(v_i)>0$ for $i<m$;
\item there is an edge $e_i$ between $v_i$ and $v_{i+1}$ such that $e_i\in E_a$ if $v_i\in V_a$;
\item $\sigma_{e_i}(v_{i+1})=1$ for $i+1<m$;
\item  $d/dt(\lambda_{e_i}(f_t))=-1$.
\end{enumerate}
By the above lemma, we have that either $s(v_m)=0$ or $\sigma_{e_{m-1}}(v_m)=-1$. Moreover, by
$(2)-(3)$ and Remark~\ref{rem:vivj} we have that $v_i<v_{i+1}$ and therefore the edges $e_i$
are all distinct.

Let $\gamma$ be the path obtained by concatenating the $e_i$'s. By $(2)-(3)$,
$\gamma$ is a legal path in the tension graph. So let

$$L=\sum_i L_X(e_i)=L_X(\gamma)\qquad L_t=\sum_i L_Y(f_t(e_i))=L_Y(f_t(\gamma)).$$

Since the $e_i$'s are in the tension graph and by condition $(4)$ we have $$L_t=\lambda(f_t) L\qquad \frac{d}{dt}L_t=-L$$

On the other hand $-\frac{d}{dt}L_t\geq S$ because by $(3)$ the contributions of the speeds of $v_i$
does not count for $i=2,\dots,m-1$ and $f(v_m)$ either stay  or moves towards $f(v_1)$ . It follows
that
$$S\leq L\leq \vol(X).$$

It follows that for any vertex $v$ we have
$$d_Y(f(v),f_t(v))\leq\int_0^t\left|\frac{d}{ds}f_s(w)\right|ds\leq\int_0^tS=tS\leq t\vol(X)$$
hence
$$d_{\infty}(\wopt(f),f))=d_\infty(f_{\lambda(f)-\lambda},f)\leq (\lambda(f)-\lambda)\vol(X).$$

We prove now the last claim. If $\wopt(f)$ is optimal then we are done. Otherwise, there is
some one-gated vertex in $X_{\max}$. We start moving the one-gated vertices as described above,
for an arbitrarily small amount. Let $g$ be the map obtained, clearly we can make
$d_\infty(g,\wopt(f))$ arbitrarily small.
Since $\wopt(f)$ is optimal, we must have $\lambda(g)=\lambda(\wopt(f))$. It follows that there is a
core sub graph of $X_{\max}$ which survives the moving. In particular, every vertex of
$X_{\max}(g)$ is at least two-gated, hence $g$ is optimal.
\qed

\medskip

\begin{defn}
We denote by $\opt(f)$ any optimal map obtained from $\wopt(f)$ as described in the proof of Theorem~\ref{Lemma_opt}.
\end{defn}

\medskip

%
%

Let $X,Y\in\O(\Gamma)$ and let $f:X\to Y$ be an optimal map. Let $v$ be a vertex of $X$ having an
$f$-illegal turn $\tau=(e_1,e_2)$. Since $f(e_1)$ and $f(e_2)$ share an initial segment, we can
identify an initial segment of $e_1$ and $e_2$. We obtain a new element $X'\in \O(\Gamma)$, with an
induced optimal map, still denoted by $f$, $f:X'\to Y$. This is a particular case of Stallings
folding (\cite{Sta}). We refer to~\cite{FM13} for further details.

\begin{defn}
  We call such an operation {\bf simple fold directed by $f$}.
\end{defn}

\medskip We finish this section by proving the existence of optimal maps with an additional
property, that will be used in the sequel.

\begin{defn}
  Let $X,Y\in\O(\Gamma)$. An optimal map $f:X\to Y$ is {\em minimal} if its tension graph
  consists of the union of axes of maximally stretched elements it contains.
  In other words, if any edge $e\in X_{\max}$ is contained in the axis of some element in
  $\pi_1(X_{\max})$ which is maximally stretched by $f$.
\end{defn}
Note that not all optimal maps are minimal, as the following shows.

\begin{ex}
Let $X$ be the graph consisting of two barbels joined by an edge, as in
Figure~\ref{fig:ex3.22}. All edges have length one except the two lower loops that have length two.
\begin{figure}[htbp]
  \centering
  \begin{tikzpicture}[x=1ex,y=1ex]
   \foreach \x in {0,25} \draw (\x,0) circle (4);
   \foreach \x in {0,29} \draw (\x-2,-14) circle (6);
   \draw (4,0) -- (21,0);
   \draw (4,-14) -- (21,-14);
   \draw (12.5,0) -- (12.5,-14);
   \foreach \x in {(-4,4), (29,4), (8,2), (17,2), (11,-5), (8,-16), (17,-16)} \draw \x node
   {$1$};
   \foreach \x in {(-9,-19), (34,-19)} \draw \x node {$2$};
   \foreach \x in {(4,0),(21,0),(-8,-14),(33,-14)} \filldraw \x circle(2pt);
   \draw (2.5,0) node {$x$} (22.5,0) node {$y$} (-11,-14) node {$f(x)$} (36,-14) node {$f(y)$};
   \draw[red] (-8,-14) arc (-180:-10:5.5) -- (12.5,-14);   
  \end{tikzpicture}
  \caption{A non-minimal optimal map. The dots $f(x)$ and $f(y)$ are not vertices, all other
    crossings are. The red line is the image of the left ``bar-edge'' of the top barbell. }
  \label{fig:ex3.22}
\end{figure}
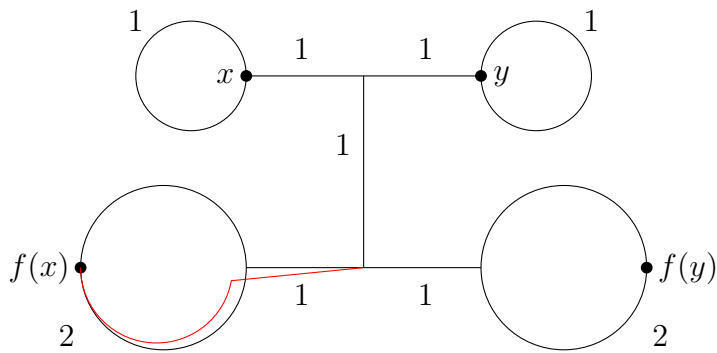

Let $f:X\to X$ be the $\PL$-map that exchanges the the top and bottom barbells (preserving left and right) and maps $x$ to the middle point
of the lower left loop, and $y$ to the middle point of the lower right loop (see the figure).

The restriction of $f$ to the lower barbell is $1$-Lipschitz (each loop is shrunk and the bar is the same length as its image), while the stretching factor of
all top edges is two. Hence the tension graph $X_{\max}$ is the top barbel. The map is optimal
because all vertices of $X_{\max}$ are two gated, but the ``bar-edges'' of the top barbel are
not in the axis of any maximally stretched loop. This is because the only legal loops in
$X_{\max}$ are the two lateral loops of the barbell. Clearly this map can be homotoped to a map
with smaller tension graph. As the next theorem shows this is always the case for non-minimal
optimal maps.
\qed
\end{ex}

\begin{thm}\label{thmminopt}
  Let $X,Y\in\O(\Gamma)$ and let $f:X\to Y$ be an optimal map. If $f$ locally minimizes the
  tension graph amongst all optimal maps $X\to Y$, then $f$ is minimal. Moreover, given $g:X\to
  Y$ optimal, for any $\e>0$ there is a minimal optimal map $f:X\to Y$ with
  $d_\infty(g,f)<\e$.
\end{thm}
\proof The first claim clearly implies the second, because the tension graph is
combinatorially finite, hence the set of possible tension graphs is finite and we can always
locally minimize it.

We will prove the contrapositive, that if $f$ is not minimal then we can decrease the tension
graph by perturbations as small as we want. The spirit is similar to that of the proof of Theorem~\ref{Lemma_opt}.
Again we will work with graphs. At the level of
graphs, the non-minimality of $f$ translates to the fact that there is an edge
$\alpha$ in the tension graph which is not part of any legal loop in $X_{\max}$.

Let $x$ be the terminal vertex of the edge $\alpha$. We say that a path starting at $x$ is $\alpha$-legal, if it is a legal path in the tension graph, whose initial edge, $e$, is not in the same gate as  $\overline{\alpha}$. We say a loop at $x$ is $\alpha$-legal if, considered as paths, both the loop and its inverse are $\alpha$-legal.

If the terminal vertex of $\alpha$ admits an $\alpha$-legal loop and the initial point of $\alpha$ also admits an  $\overline{\alpha}$-legal loop, then we can form the concatenation of these loops with $\alpha$ to get a legal loop in the tension graph crossing $\alpha$ and contradicting our hypothesis. Hence we may assume that the endpoint (rather than the initial point) of $\alpha$ is a vertex, $x$, which admits no $\alpha$-legal loops.

We will show
that it is possible to move the image of $x$ a small amount (and possibly some other vertices) so that we obtain an optimal map with smaller
tension graph.  Let $\e$ small enough so that if an edge is not in $X_{\max}$, than it remains
outside the tension graph for any perturbation of $f$ by less than $\e$.

From now on, we restrict ourselves to the tension graph. We say that a vertex $v$ is legally seen from
$x$ if there is an $\alpha$-legal path $\gamma$ from $x$ to $v$. Note that in this case $v$ is free. Indeed, otherwise the path $\gamma$ followed by its
inverse is in fact  an $\alpha$-legal loop (it has a legal lift to $\wt X$ defined by using the action of the stabilizer
of $v$). Since $v$ is free, we can move $f(v)$. Also observe that the initial point of $\alpha$ is not $\alpha$-legally seen from $x$, since otherwise we would get a legal loop in the tension graph containing $\alpha$. 

We want to chose a direction to move the images of vertices $\alpha$-legally seen from $x$. First, the direction we choose for $f(x)$ is given by the gate of $\alpha$. That is, we move $f(x)$ so as to reduced the length of $\alpha$. For any vertex, $v$, $\alpha$-legally seen from $x$, via a path $\gamma$, we move $f(v)$ backwards via the last gate of $\gamma$. That is, we move $f(v)$ so as to retrace $\gamma$. Note that this direction depends only on $v$ and not on the choice of $\gamma$. This is because, were there to be another $\alpha$-legal path from $x$ to $v$, $\gamma'$, then the concatenation $\gamma \overline{\gamma'}$ would define an $\alpha$-legal loop at $x$ unless the terminal edges of $\gamma$, $\gamma'$ lie in the same gate. Hence directions are well defined. 

We move the images of all vertices by $\e$ in the direction given above. Consider an edge, $\beta$ (not equal to $\alpha$ or its inverse) in the tension graph. If neither vertex of $\beta$ is $\alpha$-legally seen from $x$, then the image of $\beta$ is unchanged and it remains in the tension graph. Otherwise, suppose that the initial vertex of $\beta$ is $\alpha$-legally seen from $x$, via a path $\gamma$, whose terminal edge is $\overline{\eta}$. If $\eta $ and $\beta$ are in different gates, then the terminal vertex of $\beta$ is also $\alpha$-legally seen from $x$ and both vertices are moved the same amount, such that the length of the image of $\beta$ remains unchanged. If, conversely, $\eta $ and $\beta$ are in the same gate then either the length of the image of $\beta$ is reduced (if the terminal vertex is not $\alpha$-legally seen) or it remains unchanged (if it is). Moreover, since the initial vertex of $\alpha$ is not $\alpha$-legally seen, the length of the image of $\alpha$ must strictly decrease. In particular, $\alpha$ itself is no longer in the tension graph. 

On the other hand, since the tension graph has no one-gated vertices, there is at least one $\alpha$-legal
path emanating from $x$, an so some part of the tension graph survives.
Since $f$ is optimal, our assumption on $\e$ implies that the
new map is optimal  and  it  has a tension graph strictly smaller than $f$.\qed

\section{Displacement function and train track maps for automorphisms}
This section is devoted to the study of train track maps from a metric point of view. The
spirit is that of~\cite{BestvinaBers,FM13}. For the rest of the section we fix $G$ and
$\Gamma=\sqcup_i\Gamma_i$ as in Notation~\ref{not:gamma}. We recall that the study of $\O(G)$
is a particular case of $\O(\Gamma)$ when $\Gamma$ is connected.

We recall the main facts proved in~\cite{FM13} for
irreducible elements of $\Aut(G)$, and we generalize such facts to the case of $\Aut(\Gamma)$,
including reducible automorphisms. Connectedness does not really play a crucial role, and most
of the proves of~\cite{FM13} work exactly in the same way. The key here is the passage from
irreducible to reducible automorphisms.

For the rest of the section, if not specified otherwise, $\phi=(\sigma,\phi_1,\dots,\phi_k)$ will be an element of $\Aut(\Gamma)$. By abuse of
notation, we will make no
distinction between  $\phi$ and the element of $\Out(\Gamma)$ it represents. We let act the symmetric group
$\mathfrak \Sigma_k$ on $\O(\Gamma)$ by permuting the components: $\sigma(X_1,\dots,X_k)=(X_{\sigma(1)},\dots,X_{\sigma(k)})$.

As usual, if there is no ambiguity we will make no distinction between $\Gamma$-graphs and
$\Gamma$-trees. If necessary we well use $\wt X$ to refer to the $\Gamma$-tree corresponding to
a $\Gamma$-graph $X$. (So $\wt X$ will be the universal covering of $X$:)

\begin{defn}
  Let $X\in\O(\Gamma)$. We say that a ($\PL$) map $f=(f_1,\dots,f_k):X\to X$ {\em represents $\phi$} if $f_i$
  maps $X_i$ to $X_{\sigma(i)}$, and $\sigma\circ f:X \to \phi X$  is an $\O$-map
  (resp. $\PL$-map). We say that $f$ is optimal if $\sigma\circ f$ is optimal.
\end{defn}

If $\sigma \circ f:X\to \phi X$ is an optimal map, then any fold directed by $f$
gives a new point $X'$ as
well as a new map, still denoted by $f$, such that $f$ represents $\phi$ and $\sigma\circ
f:X'\to\phi X'$ is an optimal map. (This follows from Theorem~\ref{Lemma_opt}, see~\cite{FM13}
for more details.)

\begin{defn}[Displacements]\label{defdispl}
  For any automorphism $\phi\in\Out(\Gamma)$ we define the function
$$\lambda_\phi:\O(\Gamma)\to\R\qquad \lambda_\phi(X)=\Lambda(X,\phi X)$$
If $\Delta$ is a simplex of $\O(\Gamma)$ we define
$$\lambda_\phi(\Delta)=\inf_{X\in\Delta}\lambda_\phi(X)$$
If there is no ambiguity we write simply $\lambda$ instead of $\lambda_\phi$.
Finally, we set
$$\lambda(\phi)=\inf_{X\in\O(\Gamma)}\lambda_\phi(X)$$
\end{defn}

\begin{defn}[Minimally displaced points]
For any automorphism $\phi$ we define sets:
$$\Min(\phi)=\{X\in\O(\Gamma):\lambda(X)=\lambda(\phi)\}$$
$$\LocMin(\phi)=\{X\in\O(\Gamma):\exists U\ni X\text{ open s.t. } \forall Y\in U\  \lambda(X)\leq\lambda(Y)\}$$
\end{defn}

\begin{rem}
  A fold directed by a weakly optimal map does not increase $\lambda$. In particular, $\Min(\phi)$ is
  invariant by folds directed by weakly optimal maps.
\end{rem}

\begin{defn}[Reducibility]
  An automorphism $\phi$ is called {\em reducible} if there is a $\Gamma$-graph $X\in
  \O(\Gamma)$ and $f:X\to X$
  representing $\phi$ having a proper $f$-invariant subgraph $Y\subset X$ such that
  at least a component of $Y$ is not a tree with at most one non-free vertex.

Equivalently, $\phi$ is reducible if the above $\Gamma$ contains the axis of a hyperbolic element. 

  $\phi$ is {\em irreducible} if it is not reducible.
\end{defn}
\begin{rem}
In the  connected case, if $G=F_n$ then this definition coincides with the usual definition of irreducibility.  For irreducible
automorphisms we have $\Min(\phi)\neq\emptyset$. (See~\cite{FM13} for more
details.)

\end{rem}

\begin{defn}[Thin and thick reducible automorphisms]
  A reducible automorphism $\phi$ is called {\em thick} if $\Min(\phi)\neq \emptyset$, and it
  is called {\em thin} otherwise.
\end{defn}

\begin{defn}[Thin and thick simplices]
  Given a (reducible) automorphism $\phi$, a simplex $\Delta$ of $\O(\Gamma)$ is called {\em
    $\phi$-thick} (or simply {\em thick} for for short) if
  $$\inf_\Delta\lambda_\phi \text{ is realized at a
  point of }\overline{\Delta}.$$
  Otherwise $\Delta$   is called {\em $\phi$-thin}. (Recall that $\overline \Delta$ means the
  finitary closure of $\Delta$.)
\end{defn}

\begin{rem}
If $\phi$ is irreducible, then any simplex is $\phi$-thick. (See for
instance~\cite[Section~$8$]{FM13}. See also Proposition~\ref{propfinvariance} below.)  In~\cite{BestvinaBers,FM13} irreducible and thick reducible automorphisms
are called hyperbolic.
\end{rem}

\begin{defn}[Train track between trees]
  Let $\sim$ be a gate structure on a (not necessarily connected) tree $X$. A map
  $f:X\to X$ is a {\bf train track map w.r.t. $\sim$} if
  \begin{enumerate}
  \item $\sim$ is a train track structure (i.e. vertices have at least two gates);
  \item $f$ maps edges to legal paths (in particular, $f$ does not collapse edges);
  \item for any vertex $v$, if $f(v)$ is a vertex, then $f$ maps inequivalent germs at $v$ to
    inequivalent germs at $f(v)$.
  \end{enumerate}
\end{defn}

We already defined the gate structure $\sim_f$ induced by a $\PL$-map (Definition~\ref{pr4_gate}).
\begin{defn}
  Let $X$ be a (not necessarily connected) tree, and let $f:X\to X$ be a maps whose components
  are $\PL$. We define the gate structure $\langle \sim_{f^k}\rangle$ as the equivalence relation on germs generated
  by all $\sim_{f^k}, \ k\in \N$.
\end{defn}

\begin{lem}\label{Lemmatt1}
  Let $\phi\in\Aut(\Gamma)$, $X\in\O(\Gamma)$ and $\sim$ be a gate structure on $X$. Let
  $f:X\to X$ be a $PL$-map representing $\phi$. If
  $f:X\to X$ is a train track map w.r.t. $\sim$, then $\sim\supseteq\simfk$. In particular
  if $f$ is a train track map w.r.t. some $\sim$ then it is a train track map w.r.t $\simfk$.
\end{lem}
See~\cite[Section~8]{FM13} for a proof (where it is proved  in the connected case, but
connectedness plays no role).

Now we give a definition of train track map
representing an automorphism. Our definition is given at once for both reducible and
irreducible automorphisms, and in the irreducible case is the standard one. This definition
well-behaves with respect to the displacement function in the reducible case.

\begin{defn}[Optimal train track maps for automorphisms]\label{defttm}
  Let $\phi \in\Out(\Gamma)$. Let $X$ be a $\Gamma$-graph in $\O(\Gamma)$ and let $f:X\to X$ be a $\PL$-map representing
  $\phi$. Then we say that $f$ is a
   \begin{itemize}
  \item {\bf strict train track map} if there is an $f$-invariant sub-graph
  $Y\subseteq X_{\max}(f)$ such that $f|_Y$ is a train track map w.r.t. $\sim f$.
  \item {\bf train track map} if there is an $f$-invariant sub-graph
  $Y\subseteq X_{\max}(f)$ such that $f|_Y$ is a train track map w.r.t. $\simfk$.
  \end{itemize}
\end{defn}

Here some remarks are needed. The theory of train tracks maps, introduced
in~\cite{BestvinaHandel}, does not have a completely standard terminology, especially for reducible
automorphisms. We want to describe the main properties of train track maps, comparing
topological and metric viewpoints. Usually, {\em topological} train track maps are defined
without requiring that the $f$-invariant sub-graph is in the tension graph.\footnote{Our present
definition of train track map coincides with the notion of {\em optimal} train track map given
in~\cite{FM13} for irreducible automorphisms in the connected case.}

In the case $\phi$ is irreducible there is no much difference. Indeed if $f:X\to X$ is a topological
train track map representing $\phi$, then one can rescale the edge-lengths of $X$ so that $f$
is a train track map for Definition~\ref{defttm}.  And the same holds true if
$f$ has no proper invariant sub-graphs. This is because train track maps does not collapse edges,
hence edge-lengths can be adjusted so that every edge is stretched the same. In particular, the
following two results are proved in~\cite{FM13} for irreducible automorphisms and $\Gamma$
connected. The proves for generic automorphisms are basically the same (details are left to the reader).

\begin{lem}\label{Lemmatt2}
  Let $\phi\in\Aut(\Gamma)$, $X\in\O(\Gamma)$, and $f:\wt X\to \wt X$ be a $\PL$-map representing
  $\phi$. Then $f$ is train track if and only if there is an immersed periodic line $L$ in $\wt
  X_{\max}$  such that $f^k(L)\subseteq
  \wt X_{\max}$  and $f^k|_L$ is injective for all $k\in\N$. In particular if $f$ is train track then
  \begin{enumerate}
  \item $f^k$ is train track;
  \item $\Lip(f)=\Lambda(X,\phi X)$ (hence $f$ is weakly optimal);
  \item $\Lip(f)^k=\Lip(f^k)=\Lambda(X,\phi^kX)$.
  \end{enumerate}
\end{lem}

\begin{cor}
    Let $\phi\in\Aut(\Gamma)$, $X\in\O(\Gamma)$, and $f:\wt X\to \wt X$ be a map
    representing
  $\phi$. Suppose that there is an embedded periodic line $L$ in $\wt X$ such that
  $f^k|_L$ is injective for all $k\in\N$.  Suppose moreover that $\cup_kf^{k}(L)=\wt X$.
  Then there is $X'$ obtained by rescaling edge-lengths of $X$ so that $\PL(f):\wt X'\to \wt X'$ is a
  train track map.
\end{cor}

In general, if $\cup_kf^k(L)$ is just an $f$-invariant subtree $Y$ of $\wt X$, we can adjust edge
lengths so that every edge of $Y$ is stretched the same, but we cannot guarantee a priori that
$Y\subset X_{\max}$. However, the {\em interesting} case is when $\cup_kf^k(L)=\wt X$. 

\begin{defn}[Train track sets]
For any automorphism $\phi\in\Aut(\Gamma)$ we define sets:
$$\TT(\phi)=\{X\in\O(\Gamma):\exists f:X\to X \text{ train track}\}$$
$$\TTo(\phi)=\{X\in\O(\Gamma):\exists f:X\to X \text{ strict train track}\}$$
If we need to specify the map we write $(X,f)\in \TT(\phi)$ or $(X,f)\in\TTo(\phi)$.\footnote{We remark that, since in the irreducible case our present definition of train track
  map corresponds to that of {\em optimal} train track map of~\cite{FM13}, the two definitions of
  $\TT$ and $\TTo$ coincide with those given in~\cite{FM13}.}

\end{defn}

\begin{thm}\label{Theoremtt}
Let $\phi=(\sigma,\phi_1,\dots,\phi_k)\in\Aut(\Gamma)$. Then $$\overline{\TTo(\phi)}=\TT(\phi)=\Min(\phi)=\LocMin(\phi)$$
where the closure is made with respect to the simplicial topology.
\end{thm}
\proof If $\phi$ is irreducible and $\Gamma$ connected, the proof is given in~\cite{FM13} and goes trough the
following steps:
\begin{enumerate}
\item $\TTo(\phi)\subset \TT(\phi)\subseteq\Min(\phi)$.
\item If $X$ locally minimizes $\lambda_\phi$ in $\Delta_X$ then $X_{\max}=X$.
\item $\TTo(\phi)$ is dense in $\LocMin(\phi)$.
\item $\TT(\phi)$ is closed.
\end{enumerate}
We now adapt the proof so that it works also for  $\phi$ reducible and general $\Gamma$. Clearly $\Min(\phi)\subseteq \LocMin(\phi)$. By Lemma~\ref{Lemmatt1} $\TTo(\phi)\subseteq \TT(\phi)$. We see
now that $\TT(\phi)\subseteq\Min(\phi)$. If $X\in\TT(\phi)$ and $\lambda(X)>\lambda(\phi)$ then
there is $Y\in\O(\Gamma)$ such that $\lambda(Y)<\lambda(X)$. By Lemma~\ref{Lemmatt2}
$\Lambda(X,\phi^kX)=\lambda(X)^k$ but then
$$\lambda(X)^k=\Lambda(X,\phi^kX)\leq\Lambda(X,Y)\Lambda(Y,\phi^kY)\Lambda(\phi^kY,\phi^kX)$$
$$=\Lambda(X,Y)\Lambda(Y,\phi^kY)\Lambda(\phi^kY,\phi^kX)
\leq\Lambda(X,Y)\Lambda(Y,X)\lambda(Y)^k$$
thus $(\frac{\lambda(X)}{\lambda(Y)})^k$ is bounded for any $k$, which is impossible if
$\frac{\lambda(X)}{\lambda(Y)}>1$.

Thus we have
$$\TTo(\phi)\subseteq\TT(\phi)\subseteq\Min(\phi)\subseteq\LocMin(\phi).$$
\begin{lem}\label{Lemmatt3}
  Suppose $(X,f)$ locally minimizes $\lambda$ in $\Delta_X$. Then there is $Y\subseteq
  X_{\max}$ which is $f$-invariant.
\end{lem}
\proof For every open neighbourhood $U$ of $X$ we choose a point $X^U$ such that
\begin{itemize}
\item it locally minimizes $\lambda$ (a priori $X^U$ can be $X$)
\item it locally minimizes the tension graph with respect to the inclusion.
\end{itemize}
We still denote by $f$ the optimal map $f:X^U\to X^U$ obtained by optimizing $f$ w.r.t. the metric
of $X^U$ (see Theorem~\ref{Lemma_opt}). If $f(X^U_{\max})$ contains an edge $e$
which is
not in the tension graph, then by shrinking a little such edge, either we reduce $\lambda$ ---
which is impossible --- or we reduce the tension graph --- which is impossible too ---. Thus
$X^U_{\max}$ is $f$-invariant. By choosing a family of nested neighbourhoods $U_i$  we provide a
sequence $X^{U_i}\to X$ having an invariant subgraph in the tension graph. At the limit we get
an invariant subgraph of the tension graph of $X$.\qed

\begin{lem}\label{Lemmatt4}
  $\LocMin(\phi)\subseteq \overline{\TTo(\phi)}$. More precisely, let $X\in\O(\Gamma)$ and fix
  $f:X\to X$ an optimal map representing $\phi$. Suppose $X$ has an open
  neighbourhood $U$ such that for any $Y\in U$ obtained from $X$ by a sequence of simple folds
  directed by $f$, we have $\lambda(X)\leq\lambda(Y)$. Then there
  is $Y_n\in U$, all contained in the same simplex, with $Y_n\to X$ and $Y_n\in \TTo$.
\end{lem}
\proof The proof is basically the same as in~\cite{FM13}.
We recall that for $Y$ obtained from $X$ by folds directed by $f$, we
still denote $f$ the induced optimal map. First we remark that if $Y$ is obtained from $X$ by
folds directed by $f$ then $\lambda(Y)\leq\lambda(X)$ and by minimality of $X$ we have $\lambda(Y)=\lambda(X)$.
We consider the gate structure induced by $f$.
We call a vertex of $Y_{\max}$ {\em foldable} if it
has at least two elements of $Y_{\max}$ in the same gate.

Locally, by using as small as we want folds
in $X_{\max}$, directed by $f$, we find $Y\in U$ such that
\begin{enumerate}
\item $\lambda(Y)=\lambda(X)$;
\item $Y$ maximizes the dimension of $\Delta_Y$ among points reachable from $X$ via folds
  directed by $f$;
\item $Y$ minimizes $Y_{\max}$ among points of $\Delta_Y$ satisfying $(1)$ and $(2)$;
\item $Y$ maximizes the number of foldable vertices of $Y_{\max}$ among points satisfying $(1),(2),(3)$.
\end{enumerate}
Let $A\subseteq Y_{\max}$ be an $f$-invariant subgraph given by Lemma~\ref{Lemmatt3}. We
claim that $f|_A$ is a strict train track map.
Indeed, otherwise there is either an edge $e$ or a legal turn $\tau$ in $A$ having
       illegal image. Let $v$ be the vertex of $\tau$.
       \begin{itemize}
       \item If $f(e)$ contains an illegal turn $\eta$ then by folding it a little, we would
         reduce the tension graph, contradicting $(3)$. (Note that $\eta\subset Y_{\max}$
         because $A\subseteq Y_{\max}$ is $f$-invariant, thus by folding $\eta$ we do not
         change simplex of $\O(\Gamma)$ because of $(2)$.)
       \item If $f(\tau)$ is an illegal turn $\eta$ then we fold it a little.  Either
         $Y_{\max}$ becomes one-gated at $v$, and in this case the optimization process reduces
         the tension graph, contradicting $(3)$, or $v$ were not foldable at $Y$ and became
         foldable, thus contradicting $(4)$.
       \end{itemize}
Finally, note that given such an $Y$, the sequence $Y_n$ can be chosen in $\Delta_Y$.
\qed

In particular, since $\Min(\phi)$ is clearly closed, we now have:
$$\LocMin(\phi)\subseteq
\overline{\TTo(\phi)}\subseteq\overline{\TT(\phi)}\subseteq\overline{\Min(\phi)}=\Min(\phi)\subseteq\LocMin(\phi)$$
hence all inclusions are equalities.
\begin{lem}
$\overline{\TT(\phi)}=\TT(\phi)$.
\end{lem}
\proof Let $X\in\overline{\TT(\phi)}=\Min(\phi)$. Let $f:\wt X\to \wt X$ be an optimal map representing
$\phi$. By Lemma~\ref{Lemmatt4} there is $X_n\to X$ and $f_n\to f$ so that
$(X_n,f_n)\in\TTo(\phi)$. By Lemma~\ref{Lemmatt2} there is an immersed periodic line $L_n$ in
$(\wt X_n)_{\max}$ such that $f_n^k(L_n)\subset (\wt X_n)_{\max} $ is embedded for all $k\in\
N$. Since the points $X_n$ belong to the same simplex, we can suppose that all the $L_n$ are in
fact the same line $L$. Since $f_n\to f$ and the maps are all $\PL$, $L\subset \wt X_{\max}$ and  $f^k(L)\subset \wt X_{\max}$.
Moreover, if $f^k$ where not injective on $L$ for some $k$, then we
could find $\varepsilon >0$ and point $p,q$ with $d_X(p,q)=\varepsilon$ and
$f^k(p)=f^k(q)$. Now the fact that $f_n\to f$ would contradict the fact that $f_n^k|_L$ is a
homothety of ratio $\lambda(\phi)$. Thus $f^k|_L$ is embedded for any $k$, $f$ is a train
track map and so $X\in\TT(\phi)$.\qed

This completes the proof of Theorem~\ref{Theoremtt}.\qed

\begin{cor}
  Let $\phi \in\Out(\Gamma)$. If $\LocMin(\phi)\neq\emptyset$ then $\phi$ is either
  irreducible or thick reducible.
\end{cor}

We end this section by proving a lemma which is basically a rephrasing of
Lemma~\ref{Lemmatt4} with a language which will be more useful in the final part of the paper.

\begin{defn}\label{exitp}
    Let $\phi \in\Out(\Gamma)$ a point $X\in\O(\Gamma)$ is called an {\em exit point} of
    $\Delta_X$ if for any neighbourhood $U$ of $X$ in $\O(\Gamma)$ there is a point $X_E\in U$
    finite sequence of points
    $X=X_0,X_1,\dots,X_m=X_E$ in $U$, each one obtained by a
    simple fold directed by an optimal map representing $\phi$ such that $\Delta_{X_i}$ is face
    at the finite of $\Delta_{X_{i+1}}$, such that $\Delta_{X}$ is a proper face of
    $\Delta_{X_E}$, and     such that $$\lambda_\phi(X_E)<\lambda_\phi(X)$$ (strict inequality).
\end{defn}

\begin{lem}\label{LemmaX}
  Let $\phi \in\Out(\Gamma)$ and $X\in\O(\Gamma)$ such that
  $\lambda_\phi(X)$ is a local minimum for $\lambda_\phi$ in $\Delta_X$.
  Suppose $X\notin\TT(\phi)$.

Then, for any open neighbourhood $U$ of $X$ in $\Delta_X$
  there is $Z\in U$, obtained from $X$ by folds directed by optimal maps, such that
  $\lambda_\phi(Z)=\lambda_\phi(X)$, and which admits a simple
  fold directed by an optimal map and in the tension graph, entering in a simplex $\Delta'$
  having $\Delta_X$ as a proper face. (See Figure~\ref{fig:lemmaX}.)

Moreover, by finitely many such folds we find an $X'$ s.t. $\Delta_X$ is a proper face of
$\Delta_{X'}$ and $\lambda_\phi(X')<\lambda_\phi(X)$. In particular $X$ is an exit point of $\Delta_X$.
\end{lem}
\setlength{\unitlength}{1ex}
\begin{figure}[htbp]
  \centering
  \begin{picture}(52,17)
    \put(10,2){\line(1,0){30}}
    \put(10,2){\line(1,1){15}}
    \put(40,2){\line(-1,1){15}}
    \put(30,2){\line(0,1){3}}
    \put(25,2){\makebox(0,0){$\bullet$}}
    \put(30,2){\makebox(0,0){$\bullet$}}
    \put(30,5){\makebox(0,0){$\bullet$}}
    \put(18,2){\makebox(0,0){$($}}
    \put(32,2){\makebox(0,0){$)$}}
    \put(21,0){\makebox(0,0){$U$}}
    \put(25,0){\makebox(0,0){$X$}}
    \put(30,0){\makebox(0,0){$Z$}}
    \put(30,7){\makebox(0,0){$X'$}}
    \put(22,10){\makebox(0,0){$\Delta_{X'}$}}
    \put(8,2){\makebox(0,0){$\Delta_{X}$}}
  \end{picture}
  \caption{Graphical statement of Lemma~\ref{LemmaX}}
  \label{fig:lemmaX}
\end{figure}
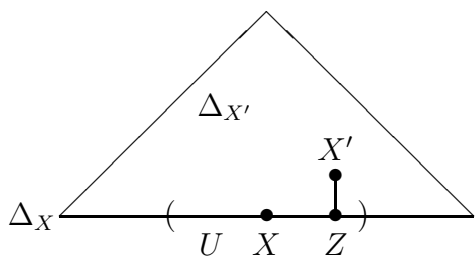
\proof Let's prove the first claim. Since $X\notin \TT(\phi)$, by Theorem~\ref{Theoremtt} there is a neighbourhood of $X$
in $\Delta_X$ which is contained in the complement of $\TTo(\phi)$. Without loss generality we may assume
that $U$ is contained in such neighbourhood.

Let $f:X\to X$ be an optimal map representing $\phi$. If there is a non-trivalent foldable
vertex in $X_{\max}$ then we set $Z=X$ and we are done. Otherwise, consider $Z\in U$ obtained
from $X$ by a fold directed by $f$ (we still denote by $f:Z\to X$ the map induced by $f$). We have $\lambda(Z)\leq \lambda(X)$. Since
$\lambda(X)$ is a local minimum, we must have $\lambda(Z)=\lambda(X)$. Let $Y\subset Z_{\max}$ be
an $f$-invariant sub-graph given by Lemma~\ref{Lemmatt3}. Since $Z\notin \TTo(\phi)$, the
restriction $f|_Y$ is not a strict train-track. It follows that by using folds directed by
optimal maps we can either
\begin{itemize}
\item[$a)$] reduce the tension graph; or
\item[$b)$] increase the number of foldable vertices; or
\item[$c)$] create a non-trivalent foldable vertex.
\end{itemize}
So far $Z$ is generic. We choose $Z\in U$ so that, in order:
\begin{enumerate}
\item it locally minimizes the tension graph;
\item it locally maximizes the number of foldable vertices among points satisfying $(1)$.
\end{enumerate}

For such a $Z$ the only possibility that remains in the above list of alternatives is $c)$, and we are done.

The last claim follows from the fact that the dimension of $\O(\Gamma)$ is bounded.\qed

\section{Behaviour of $\lambda$ at bordification points}
For the rest of the section we fix $G$ and
$\Gamma=\sqcup_i\Gamma_i$ as in Notation~\ref{not:gamma}. We also fix $\phi\in\Aut(\Gamma)$ and
we understand that $\lambda=\lambda_\phi$. In this section we discuss the behaviour of $\lambda$
when we reach points in $\partial_\infty(\O(\Gamma))$. We will see that the function $\lambda$
is not continuous and we will provide conditions that assure continuity along particular
sequences. We will also focus on the behaviour of $\lambda$ on
horospheres. In this section we will mainly think to points of $\O(\Gamma)$ as graphs.

Points near the boundary at infinity have some sub-graph that is almost collapsed. This is
usually referred to as the ``thin'' part of outer space. We will introduce a more fine notion
of ``thinness''.

\begin{defn}
  Let $\e>0$. A point $X\in\O(\Gamma)$ is {\em $\e$-thin} if there is a loop $\gamma$ in $X$
  such that $L_X(\gamma)<\e\vol(X)$.
\end{defn}

\begin{defn}
  Let $M,\e>0$. A point $X\in\O(\Gamma)$ is {\em $(M,\e)$-collapsed} if there is a loop
  $\gamma$ in $X$ such that $L_X(\gamma)<\e\vol(X)$ and for any other loop $\eta$ such that
  $L_X(\eta)\geq\e\vol(X)$ we have $L_X(\eta)>M\vol(X)$.
\end{defn}

\begin{defn}
  Let $\e>0$. For any $X\in\O(\Gamma)$ we define $X_\e$ the $\e$-thin part of $X$  as the core
  graph formed by the axes of elements $\gamma$ with $L_X(\gamma)<\e\vol(X)$.
\end{defn}

\begin{defn}
  Let $X\in\O(\Gamma)$. A sub-graph $A\subset X$ is called {\em $\phi$-invariant} if there is a
  $\PL$-map $f:X\to X$ representing $\phi$ such that $f(A)\subseteq A$.
\end{defn}

We now state some easy facts, the first of which can be found in \cite{BestvinaBers}.

\begin{prop}
  For any $C>\lambda(\phi)$ there is $\e>0$ such that for any $X\in\O(\Gamma)$, if $\lambda(X)<C$ and
  $X_\e\neq\emptyset$ then $X$
  contains a $\phi$-invariant subgraph.
\end{prop}
For a proof in the case $\Gamma$ is connected see~\cite[Section 8]{FM13} (connectedness plays
in fact no role).

However, we will need a slightly more precise statement, in order to be able to determine a particular invariant subgraph.

\begin{prop}\label{propfinvariance}
  Let $C \geq 1$ and $M>0$. 

\noindent  
  Let $\e = 1/2 \min \{ M/CD, 1/D \}$, where $D$ is the maximal number of (orbits of) edges for any graph in  $\O(\Gamma)$. Then, for  $X\in\O(\Gamma)$, if $\lambda_\phi(X)<C$ and $X$
  is $(M,\e)$-collapsed, then $X_\e$ is not the whole $X$ and it is $\phi$-invariant.
\end{prop}
\proof By definition any edge in $X_\e$ is shorter than $\e\vol(X)$. Thus we have $\vol(X_\e)<\e\vol(X)D$. In particular, since $\e D<1$ then
$X_\e\neq X$ (and thus there exists a loop $\eta$ with $L_X(\eta)>\e\vol(X)$, whence $L_X(\eta)>M\vol(X)$), since $X$ is $(M,\e)$-collapsed).

\medskip

Let $f:X\to X$ be an optimal $\PL$-map representing $\phi$. 

By picking a maximal tree in the quotient, we may find a generating set of the fundamental group of (each component of) $X_\e$ whose elements have length at most $2 \vol(X_\e)$. For any such generator, $\gamma$, we have that $L_X(f(\gamma))/L_X(\gamma) \leq C$ and hence, $L_X(f(\gamma)) \leq C L_X(\gamma) \leq 2C \vol(X_\e) < 2C D \e \vol(X) \leq M \vol(X)$. But since $X$ is $(M,\e)$-collapsed, we get that $L_X(f(\gamma)) < \e \vol(X)$. Hence $f(\gamma)$ is homotopic to a loop in $X_\e$. 

Varying $\gamma$ we deduce that $X_\e$ is $\phi$-invariant. 
\qed

\begin{prop}\label{oss}
  Let $X\in\O(\Gamma)$ and $\phi\in\Aut(\Gamma)$. Suppose that $A\subset X$ is a
  $\phi$-invariant core graph. Then $\lambda_{\phi|_A}(A)\leq \lambda_\phi(X)$.
\end{prop}
\proof Let $f:X\to X$ be a $\PL$-map representing $\phi$. Since $A$ is $\phi$-invariant,
$f(A)\subset A$ up to homotopy. By passing to the universal covering we see that $f|_A:A\to X$
retracts to a map $f_A:A\to A$ representing $\phi$ with $\Lip(f_A)\leq\Lip(f)$, hence
$\lambda_{\phi|_A}(A)\leq\Lip(f_A)\leq \Lip(f)=\lambda_\phi(X)$.\qed

\begin{thm}[Lower semicontinuity of $\lambda$]\label{fatto1}
  Fix $\phi\in\Aut(\Gamma)$ and $X\in\O(\Gamma)$. Let $(X_i)_{i\in\N}\subset \Delta_X$ be a
  sequence such that for any $i$,  $\lambda_\phi(X_i)<C$ for some $C$. Suppose that $X_i\to
  X_\infty\in\partial_\infty\Delta_X$ which is  obtained from $X$  by collapsing a sub-graph
  $A\subset X$. Then $\phi$ induces an element of $\Aut(X/A)$, still denoted by
  $\phi$.

  Moreover $\lambda_\phi(X_\infty)\leq \liminf_{i\to\infty} \lambda_\phi(X_i)$, and
  if strict inequality holds, then there is a sequence of
  minimal optimal maps $f_i:X_i\to X_i$ representing $\phi$ such that eventually on $i$ we have
  $(X_i)_{\max}\subseteq \core(A)$.
\end{thm}
\proof Let $M$ be the ``systole'' of $X_\infty$, that is to say the shortest length of simple
non-trivial loops in $X_\infty$. For any $M/\vol(X)>\e>0$, eventually on $i$, $X_i$ is
$(M/2\vol(X),\e)$-collapsed and $(X_i)_\e=\core(A)$. By Proposition~\ref{propfinvariance} $A$ is
$\phi$-invariant, thus $\phi\in\Aut(X/A)$.

For any loop $\gamma$ the lengths $L_{X_i}(\gamma)$ and $L_{X_i}(\phi(\gamma))$ converge to
$L_{X_\infty}(\gamma)$ and $L_{X_\infty}(\phi(\gamma))$ respectively.  Therefore, if $\gamma$
is a candidate in $X_\infty$ that realizes $\lambda_\phi(X_\infty)$, we have that
$\lambda_\phi(X_i)\geq L_{X_i}(\phi(\gamma))/L_{X_i}(\gamma)\to \lambda_\phi(X_\infty)$ whence
the lower semicontinuity of $\lambda$.

On the other hand, by Theorem~\ref{thmminopt} for any $i$ there is a minimal optimal map
$f_i:X_i\to X_i$ representing $\phi$. Let $\gamma_i$ be a candidate that realizes
$\lambda_\phi(X_i)$, i.e. a $f_i$-legal candidate in $(X_i)_{\max}$. Since $X$ is
combinatorically finite, we may assume w.l.o.g. that $\gamma_i=\gamma$ is the same loop for any
$i$. We have
$$\lambda_\phi(X_i)=\frac{L_{X_i}(\phi(\gamma))}{L_{X_i}(\gamma)}\to \frac{L_{X_\infty}(\phi(\gamma))}{L_{X_\infty}(\gamma)}$$

Thus if $L_{X_\infty}(\gamma)\neq 0$ we have
$\lambda_\phi(X_\infty)=\liminf\lambda_\phi(X_i)$. It follows that if there is a jump in $\lambda$
at $X_\infty$, then any legal candidate is contained in $A$. Since $f_i$ is minimal this
implies that $\core(A)$ contains the whole tension graph.\qed

\begin{rem}
 A comment on Theorem~\ref{fatto1} is required. To avoid cumbersome notation, we have decided to denote by
 $\phi$ both the element of $\Aut(X)$ and the one induced in $\Aut(X/A)$. So when we write
 $\lambda_\phi(X_\infty)$ we mean $\Lambda(X_\infty,\phi X_\infty)$ as elements in $\O(X/A)$. In
 particular, $\lambda_\phi=\inf_X\lambda_\phi(X)$ can be different if computed in $\O(X)$ or in
 $\O(X/A)$. When this will be crucial we will specify in which space we take the infimum.

Moreover, if $\phi|_A$ is the restriction of $\phi$ to $A$, then $\lambda_{\phi|_A}$ is
 calculated in the space $\O(A)$. While the simplex $\Delta_{X_\infty}$ is a simplicial face
 of $\Delta_X$, $\Delta_A\in\O(A)$ has not that meaning. One could argue that $\Delta_A$ is the
 simplex ``opposite'' to $\Delta_{X_\infty}$ in $\Delta$, but $\phi$ does not necessarily
 produces an element of $\Aut(X/(X\setminus A))$ as the complement of $A$ may be not invariant.
\end{rem}

Clearly, if $A\subset X$ is $\phi$-invariant then $\lambda_\phi(X/A)<\infty$. On the other
hand, if $A$ is not $\phi$-invariant, its collapse makes $\lambda$ explode. Thus we can extend
the function $\lambda$ as follows.

\begin{defn}\label{dli}
  Let $X_\infty\in\partial_\infty\O(\Gamma)$. We say that $\lambda_\phi(X_\infty)=\infty$ if
  $X_\infty$ is obtained from a $\Gamma$-graph $X$ by collapsing a core sub-graph $A\subset X$
  which is not $\phi$-invariant.
\end{defn}

In general, the function $\lambda$ is not uniformly continuous with respect to the Euclidean metric,
even in region where it is bounded, and so we cannot extend it to the simplicial closure of
simplices. However we see now that the behaviour of $\lambda$ is controlled on segments.

We recall the description of horospheres given in~\ref{sechor}. Suppose that $X_\infty$ is
obtained from a $\Gamma$-graph $X$ by
collapsing a $\phi$-invariant core sub-graph $A=\cup_i A_i$. Let $k_i$ be the number of
germs of edges incidents to $A_i$ in $X\setminus A$. Then $\Hor(X_\infty)$
is a product of outer spaces with marked points $\O(A_i,k_i)$. 
\begin{nt}
We denote
$\pi:\Hor(X_\infty)\to\mathbb P \O(A)$ the projection that forgets marked points.  
\end{nt}
 Note that we
chosen $X_\infty$ not projectivized and $\mathbb P\O(A)$ projectivized.
For any
$Y\in\mathbb P\O(A)$ if
$Z\in\pi^{-1}(Y)$, then there is a scaled copy of $Y$ in $Z$. We denote by $\vol_A(Y)$ the
volume of $Y$ in $Z$. With this notation in place, we can now prove a key
regeneration lemma.

\begin{lem}[Regeneration of optimal maps]\label{lemma9}
  Fix $\phi\in\Aut(\Gamma)$ ad $X\in\O(\Gamma)$. Let $X_\infty\in\partial_\infty\Delta_X$ be
  obtained from $X$ by collapsing a $\phi$-invariant core sub-graph $A$.
  Then, for any $\PL$-map $f_A:A\to A$ representing $\phi|_A$, and for any $\e>0$
  there is $X_\e\in\Delta_X$ such that $$\lambda_\phi(X_\e)\leq
  \max\{\lambda_\phi(X_\infty)+\e, \Lip(f_A)\}.$$

More precisely, for any $Y\in \mathbb P\O(A)$ and map $f_Y:Y\to Y$ representing $\phi|_A$, for
any map $f:X_\infty\to X_\infty$ representing $\phi$,
for any $\widehat X\in\Hor(X_\infty)\cap\pi^{-1}(Y)$, and for any $\e>0$; there is
$0<\delta=\delta(f,f_Y,X_\infty,\Delta_{\widehat X})$,
such that for any  $Z\in\Delta_{\widehat X}\cap\pi^{-1}(Y)$, if $\vol_Z(Y)<\delta$  there is a $\PL$-map $f_Z:Z\to Z$
representing $\phi$ such that $f_Z=f_Y$ on $Y$ and
$$\Lip(f_Z)\leq\max\{\lambda_\phi(X_\infty)+\e, \Lip(f_Y)\}$$ (and hence the optimal map
  $\opt(f_z)$ satisfies the same inequality\footnote{We notice that while $f_Z=f_Y$ on $Y$,
  this is no longer true for $\opt (f_Z)$}).
\end{lem}
\proof For this proof we will need to work with both graphs and tree, and we will the usual $\tilde\
$-notation for the universal coverings. We denote by $\sigma:X\to X_\infty$ the map
that collapses $A$. If $A_i$ is a component of $A$, we denote by $v_i$ the non-free vertex
$\sigma(A_i)$. We set $V_A=\{v_i\}$. Let $k_i$ be the valence of $v_i$ in $X_\infty$.
For any $v_i$ let $E_i^1,\dots,E_i^{k_i}$ be the half-edges incident to $v_i$ in $X_\infty$.
Let $Y_i$ be the components of $Y$. Points in $\Delta_{\widehat X}$ are built by inserting a
scaled copy of each $Y_i$  at the $v_i$ as follows. (Now we need to pass to the universal coverings.)

For every half-edge $E_i^j$ of $X_\infty$ we choose a lift in $\wt X_\infty$.
The tree $\wt{\widehat X}$ is given by attaching $\wt E_i^j$ to a point $\wt y_i^j$ of $\wt
Y_i$, and  then equivariantly attaching any other lift of the $E_i^j$. At the level of graphs
this is equivalent to choose $y_i^j\in Y_i$. Two different choices at the level of universal
coverings differ, at the level of graphs by,  closed paths $\gamma_i^j$ in $Y_i$ and based
at $y_i^j$. The choice of the simplex $\Delta_{\widehat X}$ fixes such ambiguity. Moreover for
any two graphs in $\pi^{-1}(Y)\cap \Delta_{\widehat X}$ the points $y_i^j$ are attached to the
the same edge of $Y_i$. Let $Z\in\pi^{-1}(Y)\cap\Delta_{\widehat X}$.

Given $f_Y:Y\to Y$, consider its lift to $\wt Y$ and set $\wt z_i^j=\wt f_Y(\wt y^j_i)$. There is a
unique embedded arc $\tilde\gamma_i^j$ from $\wt z_i^j$ to $\wt y_i^j$. Let $L_i$ be the
simplicial length of $\gamma_i^j$. $L_i$ depends only on $f_Y$ and the choices of $\wt y_i^j$,
hence it depends only on $f_Y$ and $\Delta_{\widehat X}$.

Now, given $f:X_\infty\to X_\infty$, there exists a
continuous map  $g:Z\to Z$ representing $\phi$, which agrees with $f_Y$ on
$Y$ and which is obtained by a perturbation of $f$ on edges of $X_\infty$. Namely on $E_i^j$ we
need to attach $\gamma_i^j$ to $f(E_i^j)$, and in each point of $f^{-1}(V_A)$ we need to insert
a small segment whose image is a suitable path in $Y$. We refer the reader to the Appendix (section~\ref{appendix})
for an accurate and detailed discussion on the properties of such a map.
For the present purpose it is sufficient to note that there is a constant $C$ such that $g$ can
be obtained so that $\Lip(\PL(g))<\Lip(f)+C\vol(Y)$. Moreover the constant $C$ depends only on
the $L_i$'s, the paths added in $f^{-1}(V_A)$, and the edge-lengths of $X_\infty$. Hence it depends
only on $f_Y,\Delta_{\widehat X},X_\infty$.

The thesis follows by setting $\delta <\e/C$ and $f_Z=\PL(g)$.\qed

\medskip

\begin{defn}
  Fix $\phi\in\Aut(\Gamma)$. Let $X_\infty\in\partial_\infty\Delta\subset\O(\Gamma)$. We say
  that $X_\infty$ has {\em not jumped in $\Delta$} if there is a sequence of points $X_i\in\Delta$
  such that $\lambda_{\phi}(X_\infty)=\lim_i\lambda_\phi(X_i)$.  We say that
  $X_\infty\in\partial_\infty\O(\Gamma)$ {\em has not jumped} if there is a simplex
  $\Delta\subset\Hor(\Delta_{X_\infty})$ such that $X_\infty$ has not jumped in $\Delta$.
\end{defn}
The above definition is for points in $\partial_\infty\O(\Gamma)$. We decide to say that $X$ has not jumped for any $X\in\O(\Gamma)$.

Notice that even if $X_\infty$ has not jumped, there may exist a simplex
$\Delta\in\Hor(\Delta_{X_\infty})$ such that $X_\infty$ has jumped in $\Delta$. This is because
if $A$ is the collapsed part and $\phi|_A$ does not have polynomial growth, then we can choose a point in
$\O(A)$ with arbitrarily high $\lambda_{\phi|_A}$. Moreover, even if $X_\infty$ has not jumped in
$\Delta$ it may happen that $X_\infty$ is not a continuity point of $\lambda$. For example if
the collapsed part $A$ has a sub-graph $B$ which is not invariant, then the collapse of $B$
forces $\lambda$ to increase due to Proposition~\ref{propfinvariance}, and thus we can approach $X_\infty$ with arbitrarily high $\lambda$.

\begin{thm}\label{newjump}
  Let $\phi\in\Aut(\Gamma)$. Let $X\in\O(\Gamma)$ containing an invariant
  sub-graph $A$. Let $X_\infty=X/A$ and $C=\core(A)$. Then
  $$\lambda_{\phi|_C}(\Delta_C)\leq\lambda_\phi(\Delta_X).$$ Moreover $X_\infty$ has not jumped in
  $\Delta_X$ if and only if $$\lambda_\phi(X_\infty)\geq\lambda_{\phi|C}(\Delta_C)$$
if and only if
$$\lambda_\phi(X_\infty)\geq\lambda_{\phi}(\Delta_X).$$
In particular there is gap as $\lambda_\phi(X_\infty)$ cannot belong to the interval
$(\lambda_{\phi|_C}(\Delta_C),\lambda_\phi(\Delta_X))$ (if non-empty).
\end{thm}
\proof The first claim is a direct consequence of Proposition~\ref{oss}.
 Let $X_i\in\Delta_X$ with $X_i\to X_\infty$ without jump, and let $A_i$ be the metric
version of $A$ in $X_i$. Then by Proposition~\ref{oss} we have
$$\lambda_{\phi|_C}(\Delta_C)\leq\lambda_{\phi|_C}(\core(A_i))\leq\lambda_\phi(X_i)\to\lambda_\phi(X_\infty).$$
Conversely, suppose $\lambda_\phi(X_\infty)\geq\lambda_{\phi|C}(\Delta_C)$. For any $\e>0$
there is $C_\e\in\Delta_C$ and a $\PL$-map $f_{C_\e}:C_\e\to C_\e$ representing $\phi|_C$ such
that $\Lip(f_{C_\e})<\lambda_{\phi|_C}(\Delta_C)+\e$. By Lemma~\ref{lemma9} there is a point
$X_\e\in X$ and a map $f_\e:X_\e\to X_\e$ representing $\phi$ such that $X_\e\to X_\infty$ as
$\e\to 0$ and $\Lip(f_e)\leq \lambda_\phi(X_\infty)$. Then
$\lambda_\phi(X_\e)\to\lambda_\phi(X_\infty)$. The second claim is proved.
Now, if $\lambda_\phi(X_\infty)\geq\lambda(\Delta_X)$, by the first two claims it has not
jumped in $\Delta_X$. And if $X_\infty$ as not
jumped $$\lambda(\Delta_X)\leq\lambda_\phi(X_i)\to \lambda_\phi(X_\infty).$$
\qed

\begin{lem}\label{lexplode}
  Let $\phi\in\Aut(\Gamma)$ and let $X\in\O(\Gamma)$. If $\lambda(\phi)>1$, then $\lambda_\phi$
  is not bounded on $\Delta_X$.
\end{lem}
\proof If there is a loop which is not $\phi$-invariant, then by collapsing that loop we force
$\lambda_\phi$ to explode. If any loop is $\phi$-invariant then by
Theorem~\ref{sausagelemma} we get  $\lambda(\phi)=1$.
\qed

\begin{lem}\label{lconst}
  Let $\phi\in\Aut(\Gamma)$. Let $X\in\O(\Gamma)$ containing an invariant
  sub-graph $A$. Let $X_\infty=X/A$ and let $C=\core(A)$. Let $$X_t=(1-t)X_\infty+tX$$ and let
  $C_t$ be the metric
  version of $C$ in $X_t$. If $\lambda_\phi(X_\infty)<\liminf \lambda_\phi(X_t)$ then
  then for $t>0$ small enough $\lambda_\phi(X_t)$ is locally constant, more precisely we have $$\lambda_\phi(X_t)=\lambda_{\phi|_C}(\core(A_1)).$$

  In particular, this is the case if $X_\infty$ has jumped in $\Delta$ along the segment $XX_\infty$.
\end{lem}
\proof
By Lemma~\ref{fatto1} for $t$ small enough there is an optimal map $f_t:X_t\to X_t$ whose
tension graph of $X_t$ is contained in $C_t$.
Since $C_t$ is $\phi$-invariant, $f_t(C_t)\subset C_t$ up to homotopy. Since the vertices of
$(X_{t})_{\max}$ are at least two gated, $f((X_t)_{\max})\subset C_t$.
Therefore $\lambda_{\phi|_C}(C_t)=\Lip(f_t)$ and
$\lambda_\phi(X_t)=\Lip(f_t)=\lambda_{\phi|_C}(C_t)=\lambda_{\phi|_C}(C_1)$ (where the last
equality follows from the fact that $[C_t]=[C_1]\in\mathbb P\O(C)$).

The last claim follows because by Theorem~\ref{newjump}, and since $X_\infty$ has jumped in
$\Delta$, we have $$\lambda_\phi(X_\infty)<\lambda_{\phi|_C}(\Delta_C)\leq\lambda_\phi(\Delta)\leq\lambda_\phi(X_t)$$
hence $\lambda_\phi(X_\infty)<\liminf_t\lambda_\phi(X_t)$.
\qed

\begin{cor}\label{fatto2}
  Let $\phi\in\Aut(\Gamma)$. Let $\Delta$ be a simplex of $\O(\Gamma)$. If there a point in
  $\overline{\Delta}^\infty$ which jumped in $\Delta$, then there is a min-point in $\Delta$ that
  realizes $\lambda(\Delta)$.

  In particular, there is always a min-point $X_{\min}$ in
  $\overline{\Delta}^\infty$ which has no jumped and such that
  $$\lambda_\phi(\Delta_X)=\lambda_\phi(X_{\min})=\lambda_\phi(\Delta_{X_{\min}}).$$
\end{cor}
\proof Let $\Delta=\Delta_X$ for some $\Gamma$-graph $X$. Suppose that $A$ is a
$\phi$-invariant subgraph so that $X_\infty=X/A$ has a jump in $\Delta$. Let $C=\core(A)$.
By Theorem~\ref{newjump} $$1\leq\lambda(X_\infty)<\lambda(\Delta_C)\leq\lambda(\Delta)$$ hence by
Lemma~\ref{lexplode} there  $C_o\in\Delta_{C}$ be such that $\lambda(C_o)=\lambda(\Delta)$.
Thus there is $X_o\in\Delta$ obtained by inserting on $X_\infty$ a copy of a metric graph
isomorphic to $A$ and with core $C_0$. By Lemma~\ref{lconst}, for small enough $t$, we have
$\lambda(tX+(1-t)X_\infty)=\lambda(\Delta)$. The last claim follows clearly from the first.
\qed

\section{Convexity properties of the displacement function}\label{section_conv}

We recall that we are using
the terminology ``simplex'' in a wide sense, as $\Delta_X$ is a standard simplex if we work in
$\mathbb P\O(\Gamma)$ and the cone over it if we work in $\O(\Gamma)$. (Remember we use
Notation~\ref{not:gamma} for $\Gamma$.)

The function $\lambda$ is scale invariant on $\O(\Gamma)$ so it descends to a function on
$\mathbb P\O(\Gamma)$. In order to control the value of $\lambda$ on segments in terms of
its value on vertices, we would like to say that $\lambda$ is convex on segments.
A little issue appears with projectivization. If $\Delta$ is a simplex of
$\O(\Gamma)$, its euclidean segments are well defined, and their projections on $\mathbb
P(\O(\Gamma))$ are euclidean segments in the image of $\Delta$. However, the linear
parametrization is not a projective invariant (given $X,Y$, the points $(X+Y)/2$ and $(5X+Y)/2$
are in different projective classes).

It follows that convexity of a scale invariant function is not
well-defined. In fact if $\sigma$ is a segment in $\Delta$, $\pi:\Delta\to\mathbb
P\Delta$ is the projection, and $f$ is a convex function on
$\sigma$, then $f\circ\pi^{-1}$ may be not convex. It is convex only up to reparametrization of
the segment $\pi(\sigma)$.  Such functions are called quasi-convex, and this notion will be
enough for our purposes.

\begin{defn}
  A function $f:[A,B]\to \mathbb R$ is called {\em quasi-convex} if for all $[a,b]\subseteq[A,B]$
$$\forall t\in[a,b]\qquad f(t)\leq\max\{f(a),f(b)\}.$$
\end{defn}
Note that quasi-convexity is scale invariant.

\begin{lem}\label{lconvexity}
  For any $\phi\in\Aut(\Gamma)$ and for any open simplex $\Delta$ in $\O(\Gamma)$ the function
  $\lambda$ is quasi-convex on segments of $\Delta$. Moreover, if $\lambda(A)>\lambda(B)$ then
  $\lambda$ is not locally constant near $A$.
\end{lem}
\proof
Let $X$ be a $\Gamma$-graph such that $\Delta=\Delta_X$.
We use the Euclidean coordinates of $\Delta$ labelled with edges of $X$, namely a point $P$ in
$\Delta$ is given by a vector whose $e^{th}$ entry is the length of the $e$ in $P$. In the same way,
to any reduced loop $\eta$ in $X$ we associate its occurrence vector, whose $e^{th}$ entry is
the number of times that $\eta$ passes through the edge $e$. We will denote by $\eta$ both the
loop and its occurrence vector. With this notation, the length function is bilinear:
$$L_X(\gamma)=\langle X,\gamma\rangle$$
(where $\langle,\rangle$ denotes the standard scalar product on $\mathbb R^k$.)

Let $\sigma$ be a segment in $\Delta$ with endpoints $A,B$.
Let $\gamma$ be a candidate. We consider both $\gamma$ and $\phi\gamma$ as loops in $X$.
Up to switching $A$ and $B$, we may assume that
$$\frac{\langle A,\phi\gamma\rangle}{\langle A,\gamma\rangle}\geq
\frac{\langle B,\phi\gamma\rangle}{\langle B,\gamma\rangle}.$$

Such a condition is scale invariant, and since $\lambda$ is scale invariant, up to rescaling
$B$ we may assume that $\langle B,\gamma\rangle>\langle A,\gamma\rangle$. Now, we parametrize
$\sigma$ as usual with $[0,1]$
$$\sigma(t)=A_t=Bt+(1-t)A$$
We are interested in the function

$$F_\gamma(t)=\frac{\langle A_t,\phi\gamma\rangle}{\langle A_t,\gamma\rangle}
=\frac{\langle Bt+(1-t)A,\phi\gamma\rangle}{\langle Bt+(1-t)A,\gamma\rangle}
=\frac{\langle A,\phi\gamma\rangle+t\langle
  B-A,\phi\gamma\rangle}{\langle A,\gamma\rangle+t\langle B-A ,\gamma\rangle} $$

A direct calculation shows that the second derivative of a function of the type
$f(t)=(a+tb)/(c+td)$ is given by $2(ad-bc)d/(c+td)^3$.

So the sign of $F_\gamma''(t)$ is given by $$\big(\langle A, \phi\gamma\rangle \langle
B,\gamma\rangle-\langle B,\phi\gamma\rangle\langle A,\gamma\rangle\big)
\big(\langle B-A,\gamma\rangle\big)$$
which is positive by our assumption on $A,B$. Hence $F_\gamma(t)$ is convex and therefore
quasi-convex:
$$F_\gamma(t)\leq\max\{F_\gamma(A),F_\gamma(B)\}.$$

Now, by the Sausage Lemma~\ref{sausagelemma} we have:

\begin{eqnarray*}
  \lambda_\phi(A_t)&=&\max_\gamma F_\gamma(t)\leq\max\{\max_\gamma F_\gamma(A),\max_\gamma F_\gamma(B)\}\\&=&\max\{\lambda_\phi(A),\lambda_\phi(B)\}.
\end{eqnarray*}

Finally, since lengths of candidates
are finitely many, there is a candidate $\gamma_o$ such that for $t$ sufficiently small
we have $\lambda_\phi(A_t)=F_{\gamma_o}(t)$. So, if $\lambda_\phi$ is locally constant
near $A$, then we must have $F_{\gamma_o}''(t)=0$ hence

$$\lambda_\phi(A)=\lambda_\phi(A_0)=\frac{\langle A,\phi\gamma_o\rangle}{\langle A,\gamma_o\rangle}=
\frac{\langle B,\phi\gamma_o\rangle}{\langle B,\gamma_o\rangle}\leq\lambda_\phi(B).$$

\qed

\begin{lem}\label{lconv2}
  Let $\phi\in\Aut(\Gamma)$ and let $\Delta$ be a simplex in $\O(\Gamma)$. Let
  $A,B\in \overline{\Delta}^\infty$ be two points that have not jumped in $\Delta$. Then for any
  $P\in\overline{AB}$ $$\lambda(P)\leq\max\{A,B\}$$

  Moreover, if $\lambda(A)\geq \lambda(B)$, then $\lambda|_{\overline{AB}}$ is continuous at $A$.
\end{lem}
\proof Let $X$ be a graph of groups so that $\overline{AB}=\Delta_X$. By
Lemma~\ref{lconvexity}, the function $\lambda$ is quasi-convex on the interior of $\overline{AB}$ as
a segment in $\O(X)$. Let $\{A_i\}$ and $\{B_i\}$ sequences in $\Delta$ such that
 $A=\lim A_i$ and $B=\lim B_i$ with $\lim\lambda(A_i)=\lambda(A)$ and $\lim
 \lambda(B_i)=\lambda(B)$. Such sequences exists because of the non jumping hypothesis.
For all point $P$ in the
segment $\overline{AB}$, there is a sequence of points $P_i$ in the segment $\overline{A_iB_i}$
such that $P_i\to P$. By Lemma~\ref{lconvexity} we know
$$\lambda(P_i)\leq\max\{\lambda(A_i),\lambda(B_i)\},$$
and by lower semicontinuity
(Theorem~\ref{fatto1}) of
$\lambda$ and the non jumping assumption, such inequality passes to the limit. In particular,
if $\lambda(A)\geq \lambda(B)$, then $\lambda(P)\leq\lambda(A)$ for any $P\in\overline{AB}$.

Now suppose that $P^j\to A$ is a sequence in the segment $\overline{AB}$. Let
$P^j=\lim_iP_i^j$. Then by lower semicontinuity Theorem~\ref{fatto1}
applied to the space  $\O(P)$, on the segment  $\overline{AB}$ we have
$$\lambda(A)\geq\lim_j\lambda(P_j)\geq\lambda(A).$$
\qed

We end this section with an estimate of the derivative of functions like the $F_\gamma(t)$ defined as
in Lemma~\ref{lconvexity}, which will be used in the sequel. As above, we use the formalism
$\langle X,\gamma\rangle=L_X(\gamma)$.

\begin{lem}\label{derivative}
  Let $\Delta=\Delta_X$ be a simplex of $\O(\Gamma)$ and $A,B\in\overline{\Delta}^\infty$. Let
  $\gamma$ be a loop in $X$ which is not collapsed neither in $A$ nor in $B$ and set  $$C=\max\{\frac{L_A(\gamma)}{L_B(\gamma)},\frac{L_B(\gamma)}{L_A(\gamma)}\}$$
   Let $\phi$ be any  automorphism of $\Gamma$. Suppose that 
$\frac{\langle B,\phi\gamma\rangle}{\langle B,\gamma\rangle} \geq \frac{\langle
  A,\phi\gamma\rangle}{\langle A,\gamma\rangle}$. Let $A_t=tB+(1-t)A$ be the linear
parametrization of the segment $AB$ in $\Delta$ and define
$F_\gamma(t)=
\frac{\langle A_t,\phi\gamma\rangle}{\langle A_t,\gamma\rangle}$.
Then 
$$0\leq F'_\gamma(t)\leq C \frac{\langle B,\phi\gamma\rangle}{\langle B,\gamma\rangle}$$

In particular, for any point $P$ in the segment $AB$ we have
$$\lambda_\phi(P)\geq\frac{\langle P,\phi\gamma\rangle}{\langle P,\gamma\rangle}\geq \frac{\langle
  B,\phi\gamma\rangle}{\langle B,\gamma\rangle}-C\lambda_\phi(B)\frac{||P-B||}{||A-B||}$$
where $||X-Y||$ denotes the standard Euclidean metric on $\Delta$. 
\end{lem}
Before the proof, a brief comment on the statement is desirable. First, note that the constant
$C$ does not depend on $\phi$. Moreover, by taking the supremum where $\gamma$ runs over all candidates given
by the Sausage Lemma~\ref{sausagelemma}, then $C$ does not even depend on $\gamma$. 
Finally if $\gamma$ is a candidate that realizes $\lambda_\phi(B)$, then we
get a bound of the steepness of $F_\gamma$ which does not depend on $\phi$ nor on $\gamma$ but
just on $\lambda_\phi(B)$ and $||A-B||$.
\proof We have
$$F_\gamma(t)=
\frac{\langle A_t,\phi\gamma\rangle}{\langle A_t,\gamma\rangle}
=\frac{\langle Bt+(1-t)A,\phi\gamma\rangle}{\langle Bt+(1-t)A,\gamma\rangle}
=\frac{\langle A,\phi\gamma\rangle+t\langle
  B-A,\phi\gamma\rangle}{\langle A,\gamma\rangle+t\langle B-A ,\gamma\rangle}$$
and a direct calculation show that
\begin{equation}
  F'_\gamma(t)=\frac{\langle B,\gamma\rangle\langle A,\gamma\rangle}{(\langle
    A_t,\gamma\rangle)^2}\left(\frac{\langle B,\phi\gamma\rangle}{\langle B,\gamma\rangle}
-\frac{\langle A,\phi\gamma\rangle}{\langle A,\gamma\rangle}\right)
\end{equation}
The first consequence of this equation is that the sign of $F'_\gamma$ does not depend on
$t$, and since 
$\frac{\langle B,\phi\gamma\rangle}{\langle B,\gamma\rangle} \geq \frac{\langle
  A,\phi\gamma\rangle}{\langle A,\gamma\rangle}$, then $F'_\gamma\geq 0$. Moreover, 
since $\langle A_t,\gamma\rangle$ is linear on $t$, 
we have $\frac{\langle B,\gamma\rangle\langle A,\gamma\rangle}{(\langle
    A_t,\gamma\rangle)^2}\leq C$. Therefore we get 
$$F'_\gamma(t)\leq C
\frac{\langle B,\phi\gamma\rangle}{\langle B,\gamma\rangle}
$$
and the first claim is proved.
For the second claim, note that the parameter $t$ is nothing but $||A-A_t||/||A-B||$ and thus
$$F_\gamma(1)-F_\gamma(t)\leq 
(1-t)C\frac{\langle B,\phi\gamma\rangle}{\langle B,\gamma\rangle}=
\frac{||B-A_t||}{||B-A||}C\frac{\langle B,\phi\gamma\rangle}{\langle B,\gamma\rangle}.
$$
If $P=A_t$, we have $F_\gamma(1)=\frac{\langle B,\phi\gamma\rangle}{\langle B,\gamma\rangle}$
and $F_\gamma(t)=\frac{\langle P,\phi\gamma\rangle}{\langle P,\gamma\rangle}$.
 By taking in account $\lambda_\phi(B)\geq\frac{\langle B,\phi\gamma\rangle}{\langle
  B,\gamma\rangle}$ and $\lambda_\phi(P)\geq\frac{\langle P,\phi\gamma\rangle}{\langle
  P,\gamma\rangle}$ we get the result.\qed

\section{Existence of minimal displaced points and train tracks at the bordification}

The existence of points that minimizes the displacement is proved in~\cite[Theorem 8.4]{FM13} for irreducible automorphisms, but in
fact the philosophy of the proof works in the general case if we are allowed to pass to the boundary
at infinity, and by
taking in account possible jumps. A part jumps, the problem is that one cannot use compactness for
claiming that a minimizing sequence has a limit as the bordification of $\O(\Gamma)$ is not
even locally compact. The trick used in~\cite{FM13} is to use Sausage Lemma.
We use Notation~\ref{not:gamma} for $G$ and $\Gamma$.

\begin{lem}\label{trick}
  For any $\Gamma$, for any $X\in\overline{\O(\Gamma)}^\infty$ the set
  $\{\lambda_\phi(X):\phi\in\Out(\Gamma)\}$ is discrete. In other words, given $X$, all possible
  displacements of $X$ with respect to all automorphisms (and hence markings) run over a
  discrete set. 
\end{lem}
\proof
This proof is similar to that of~\cite[Theorem 8.4]{FM13}, we include it by
completeness. If $\phi\notin\Aut(X)$ (i.e. if $X$ has a collapsed part which is not
$\phi$-invariant) than $\lambda_\phi(X)=\infty$ and there is nothing to
prove. Otherwise, 
by Sausage Lemma~\ref{sausagelemma}, $\lambda_\phi(X)=\Lambda(X,\phi X)$ is
computed by the quotient of translation lengths of candidates. The possible values of
$L_X(\phi\gamma)$ (with $\gamma$ any loop) are a discrete set just because $X$ has
finitely many edges. 
 Candidates are in general
infinitely many, but there are only finitely many lengths of them. Thus the possible values of
$\Lambda(X,\phi X)$ runs over a discrete subset of $\mathbb R$.\qed

\begin{thm}\label{conj}
For any $\Gamma$ the global simplex-displacement spectrum
$$\operatorname{spec}(\Gamma)= \Big\{\lambda_\phi(\Delta):\Delta\text{ a simplex of }
\overline{\O(\Gamma)}^\infty, \phi\in\Out(\Gamma)\}$$ is well-ordered as a subset of $\mathbb R$.
In particular, for any $\phi\in\Out(\Gamma)$ the spectrum of possible minimal displacements $$\operatorname{spec}(\phi)= \Big\{\lambda_\phi(\Delta):\Delta\text{ a simplex of }
\overline{\O(\Gamma)}^\infty\}$$ is well-ordered as a subset of $\mathbb R$.
\end{thm}
\proof Recall that
we defined $\lambda_\phi(\Delta)$ as $\inf_{X\in\Delta}\lambda_\phi(X)$. For this proof we
fix the volume-one normalization, and in any simplex we use the standard Euclidean norm,
denoted by $||\cdot||$.

We argue by induction on the rank of $\Gamma$ (See Definition~\ref{pr4_rank}). Clearly if the rank of $\Gamma$ is one there is nothing to
prove. We now assume the claim true for any $\Gamma'$ of rank smaller than $\Gamma$.

We will show than any
monotonically decreasing sequence in $\operatorname{spec}(\Gamma)$ has a (non trivial)
sub-sequence which is constant, whence the original sequence 
is eventually constant itself. This
implies  that $\operatorname{spec}(\Gamma)$ is well-ordered. For the second claim, since
$\operatorname{spec}(\phi)$ is a subset of a well-ordered set, it is  well-ordered. 

We follow the line of reasoning of~\cite[Theorem 8.4]{FM13}.  Let $\lambda_i\in\operatorname{spec}(\Gamma)$
be a monotonically decreasing sequence. Note that displacements are non-negative so $\lambda_i$ converges.
For any $i$ we chose $\phi_i$ and a point $X_i\in \overline{\O(\Gamma)}^\infty$ such that
$\lambda_{\phi_i}(X_i)=\lambda_{\phi_i}(\Delta_{X_i})=\lambda_i$ and let $\Delta_i=\Delta_{X_i}$.
Up to possibly passing to sub-sequences we may assume
that there is $\psi_i\in\Out(\Gamma)$ such that $\psi_iX_i$ belongs to a fixed simplex
$\Delta$. Therefore, by replacing $\phi_i$ with $\psi_i\phi_i\psi_i^{-1}$ we may assume that
the $X_i$ all belong to the same simplex $\Delta$. Let $X$ be the graph of groups corresponding to $\Delta$,
i.e. $\Delta=\Delta_X$.\footnote{Note that $X$ may be a boundary point of $\O(\Gamma)$ and that
  we have made no assumption about jumps, so $X$ may jump.}
 
 Up to sub-sequences, $X_i$ converges to a point  $X_\infty$ in the simplicial
closure of $\Delta$. By Lemma~\ref{trick} up to possibly passing to a subsequence we may assume
that  $\lambda_{\phi_i}(X_\infty)$ is a constant $L$. (The only issue here is that $X_{\infty}$ is displaced a finite amount by each $\phi_i$ - we show this is true for all but finitely many $i$.)

Note that if $X_\infty$ is in
$\partial_\infty\overline{\Delta}$, then there exist $M,\varepsilon > 0$ such that $X_i$ is eventually
$(M,\varepsilon)$-collapsed. Namely, assuming $X_\infty$ has volume 1, take $M$ to be the length of the shortest loop in $X_\infty$, and take $\epsilon$ to be a constant small enough to satisfy the hypotheses of Proposition~\ref{propfinvariance}. Then the thin part, $(X_i)_\varepsilon$, is the core of the the sub-graph
of $X$ which collapsed to obtain $X_\infty$. By Proposition~\ref{propfinvariance}
$(X_i)_\varepsilon$ is $\phi_i$-invariant and so $\lambda_{\phi_i}(X_\infty)<\infty$. Hence $L<\infty$.

Since $X_i$ is a
 min-point for the function $\lambda_{\phi_i}$, by Lemma~\ref{lconvexity} the function
 $\lambda_{\phi_i}$ either is constant on the segment $X_iX_\infty$ or it is not locally
 constant near $X_\infty$. By Lemma~\ref{lconst} in the latter case $X_\infty$ has not jumped
 w.r.t. $\lambda_{\phi_i}$ along the segment $X_iX_\infty$. 

Therefore we have the following three cases, and 
 up to subsequences we may assume that we are in the same case for any $i$:
 \begin{enumerate}
 \item $\lambda_{\phi_i}$ is constant and continuous on $X_iX_\infty$;
 \item $\lambda_{\phi_i}$ is constant on the interior of $X_iX_\infty$ and there is a jump at
   $X_\infty$, hence  $\lambda_{\phi_i}(X_\infty)<\lambda_{\phi_i}(X_i)$ by lower
   semicontinuity Theorem~\ref{fatto1};
 \item $\lambda_{\phi_i}$ is monotone increasing near $X_\infty$ and continuous at $X_\infty$.
 \end{enumerate}

In the first case $\lambda_i=\lambda_{\phi_i}(X_i)=L$ and we are done. 
In the second case we use the inductive hypothesis. Namely, by Lemma~\ref{lconst} there is a
core $\phi_i$-invariant sub-graph $A_i$ of $X_i$ such that
$\lambda_{\phi_i}(X)=\lambda_{\phi_i|_{A_i}}(A_i)$ for any $X$ in the interior of the segment
$X_iX_\infty$. Moreover, up to sub-sequences we may assume that $A_i$ is topologically the same
graph for any $i$. If $A_i$ does not minimizes locally $\lambda_{\phi_i|_{A_i}}$ in its
simplex, the we could perturb a little $A_i$ and strictly decrease $\lambda_{\phi_i}(X_i)$
contradicting the minimality of $X_i$. By quasi-convexity Lemma~\ref{lconvexity}, in any
simplex local minima are global minima and thus $\lambda_{\phi_i|_{A_i}}(A_i)=\lambda_{\phi_i|_{A_i}}(\Delta_{A_i})$.
By induction the global simplex-displacement spectrum of $A_i$ is well ordered, hence the
decreasing sequence $\lambda_i=\lambda_{\phi_i}(X_i)=\lambda_{\phi_i|_{A_i}}(\Delta_{A_i})$ is
eventually constant.

It remains case $(3)$. In this case $$\lambda_{\phi_i}(X_i)<L=\lambda_{\phi_i}(X_\infty).$$
 Let $R>0$ be such that for any face $\Delta'$ of $\Delta$ such that $X_\infty\notin
 \overline{\Delta'}^\infty$, the ball $B(X_\infty,2R)$ is disjoint from $\Delta'$. In other
 words, if $P\in B(X_\infty,2R)$ is obtained form $X$ by collapsing a sub-graph $P_0$, then
 $P_0$ is collapsed also in $X_\infty$. Eventually on $i$, $X_i\in B(X_\infty,R)$. Let $Y_i$ be
 the point on the Euclidean line trough $X_i,X_\infty$ at distance exactly $R$ from $X_\infty$.

Let $\gamma_i$ be a candidate in $X_\infty$ that realizes $\lambda_{\phi_i}(X_\infty)$ and such
that the stretching factor $\frac{L_X(\phi_i(\gamma_i))}{L_X(\gamma_i)}$ of $\gamma_i$ locally decreases toward $X_i$. Such a $\gamma_i$ exists
because  $\lambda_{\phi_i}(X_i)<\lambda_{\phi_i}(X_\infty)$. 

By Lemma~\ref{derivative} applied with $A=Y_i$ and $B=X_\infty$ we have 
 $$\lambda_{\phi_i}(X_i)\geq \lambda_{\phi_i}(X_\infty)\left(1-C\frac{||X_i-X_\infty||}{R}\right).$$
where
$C=\max\{\frac{L_{Y_i}(\gamma_i)}{L_{X_\infty}(\gamma_i)},\frac{L_{X_\infty}(\gamma_i)}{L_{Y_i}(\gamma_i)}\}$.
Since there are finitely many lengths of candidates and by our choice of $R$, the
constant $C$ is uniformly bounded independently on $i$. Since $X_i\to X_\infty$ we have
$\varepsilon_i=||X_i-X_\infty||\to 0$ and thus
$$L(1-C\varepsilon_i)\leq \lambda_{\phi_i}(X_i)\leq L.$$
Thus $\lambda_i\to L$ and since it is a monotonically decreasing sequence bounded above by its limit,
it must be constant.\qed

\medskip
We suspect that $\operatorname{spec}(\phi)$ is not only well-ordered but in fact
discrete. However, Theorem~\ref{conj} will be enough for our purposes.

\begin{thm}[Existence of minpoints]\label{thmminptE}
  Let $\phi$ any element in $\Aut(\Gamma)$. Then there exists $X\in\overline{\O(\Gamma)}^\infty$
  that has not jumped and such that $$\lambda_\phi(X)=\lambda(\phi).$$
\end{thm}
\proof Let $X_i\in\O(\Gamma)$ be a minimizing sequence for $\lambda_\phi$. Without loss of
generality we may assume that the sequence $\lambda_\phi(\Delta_{X_i})$ is monotone
decreasing and Theorem~\ref{conj} implies that it is eventually constant. Therefore $X_i$ can be
chosen in a fixed simplex $\Delta$. Corollary~\ref{fatto2} concludes.\qed

\medskip

An interesting corollary of Theorem~\ref{thmminptE} is that we can characterize (global) jumps
extending the first statement of Theorem~\ref{newjump} from a local to a global statement.

\begin{thm}\label{thmjump}
  Fix $\phi\in\Aut(\Gamma)$. Let $X\in\O(\Gamma)$ and let $X_\infty\in\partial_\infty\Delta_X$ be obtained from $X$ by
  collapsing a $\phi$-invariant core graph $A$. Then $X_\infty$ has not jumped if and only
  if $$\lambda(\phi|_A)\leq \lambda_\phi(X_\infty).$$
\end{thm}
\proof
Suppose that $X_\infty$ has not jumped. Then there
is a simplex $\Delta$ where $X_\infty$ has not jumped, and the claim follows from
Theorem~\ref{newjump} because $\lambda_{\phi|_A}\leq\lambda_{\phi|_A}(\Delta_A)$.

On the other hand, suppose $\lambda(\phi|_A)\leq\lambda_\phi(X_\infty)$. By
Theorem~\ref{thmminptE} there is a simplex in $\O(A)$ containing a minimizing sequence for $\phi|_A$.
Let $A_\e$ be an element in that simplex so that $\lambda_{\phi|_A}(A_\e)<\lambda(\phi|_A)+\e$, and let
$f_A:A_\e\to A_\e$ be an optimal map representing $\phi|_A$. Note that $A_\e$ and $A$ may be not
homeomorphic. Let $\widehat X$ be a $\Gamma$-graph obtained by inserting a copy of $A_\e$ in
$X_\infty$. (We notice that since $A_\e$ may be not homeomorphic to $A$, we can have $\Delta_{\widehat
  X}\neq\Delta_X$. We also notice that such $\Delta_{\widehat X}$ is not unique as we have plenty
of freedom of attaching the edges of $X_\infty$ to $A_\e$.)
By Lemma~\ref{lemma9}, for any $\e>0$ there is an element $X_\e\in\Delta_{\widehat X}$ and an
optimal map $f_\e:X_\e\to X_\e$ representing $\phi$ so that $X_\e\to X_\infty$ and
 $\Lip(f_\e)\leq \lambda_\phi(X_\infty)+\e$, hence $\lambda_\phi(X_\e)\leq
(X_\infty)+\e$. Thus $X_\infty$ has no jump in $\Delta_{\widehat X}$, and therefore has not jumped.\qed

\medskip

By Theorem~\ref{Theoremtt} we know that minimal displaced points and train tracks coincide. But
some care is needed here, as that theorem is stated for point of $\O(\Gamma)$, and not for points at infinity. In fact, given $\phi\in\Aut(\Gamma)$,
$X\in\O(\Gamma)$ and $A\subset X$ a $\phi$-invariant sub-graph, a priori
it may happen that $\lambda(\phi)$ is different if we consider $\phi$ as an element of
$\Aut(X)$ or of $\Aut(X/A)$. That is to say, we may have $X_\infty=X/A$ such that
$\lambda_\phi(X_\infty)=\lambda(\phi)$ but $X_\infty$ is not a train track point in $\O(X/A)$.

For instance, consider the case where $X=A\cup B$, with both $A$ and $B$ invariant. Suppose
that $\lambda(\phi)=\lambda(\phi|_A)>\lambda(\phi|_B)$. Now suppose that
$\lambda(\phi)=\lambda_{\phi}(X)=\lambda_{\phi|_A}(A)=\lambda_{\phi|_B}(B)$. Collapse $A$. Then the resulting
point $X_\infty$ is a min point for $\phi$ in $\overline{\O(\Gamma)}^\infty$ which has not
jumped, but since $\lambda(\phi|_B)<\lambda(\phi)$, it is not a min point for $\phi$ on
$\O(X/A)$.

We want to avoid such a pathology. Here we need to make a difference between $\lambda(\phi)$
computed in different spaces, so we will specify the space over which we take the infimum.

\begin{lem}[Existence of train tracks]\label{lemmattE}
  Let $\phi\in\Out(\Gamma)$. Let $X_\infty\in\overline{\O(\Gamma)}^\infty$ be such that:
  \begin{itemize}
  \item There is $X\in\O(\Gamma)$ such that $X_\infty$ is obtained from $X$ by collapsing a
    (possibly empty) core sub-graph $A$ in $X$;
  \item $\lambda_\phi(X_\infty)=\inf_{Y\in\O(\Gamma)}\lambda_\phi(Y)$;
  \item it has no jump in $\Delta_X$.
  \end{itemize}
Suppose moreover that $X_\infty$  maximizes the dimension of $\Delta_{X_\infty}$ among the set
of elements in $\overline{\O(\Gamma)}^\infty$ satisfying such conditions (such a set is not empty
by Theorem~\ref{thmminptE}). Then
  $\lambda_\phi(X_\infty)=\inf_{Y\in\O(X/A)}\lambda_\phi(Y)$. (Hence it is in
  $\TT(\phi)\subset\O(X/A)$.)
\end{lem}
\proof If $A$ is empty this is an instance of Theorem~\ref{Theoremtt}. Otherwise,
suppose $X_\infty$ is not a train track point of $\O(X_\infty)$. We claim that near $X_\infty$
there is a point $X'_\infty\in\O(X_\infty)$ such that $\lambda_\phi(X_\infty')<\lambda_\phi(X_\infty)$.
Indeed, if $X_\infty$ is not a local
min point in $\Delta_{X_\infty}\subset\O(X_\infty)$, then we can find $X'_\infty$ just near
$X_\infty$ in $\Delta_{X_\infty}$. Otherwise, by Lemma~\ref{LemmaX}
there is a point $X'_\infty$ obtained form $X_\infty$ by folds directed by optimal maps (and such
that $\dim(\Delta_{X'_\infty})>\dim(\Delta_{X_\infty})$)  such that
$\lambda_\phi(X_\infty')<\lambda_\phi(X_\infty)$. 

Let $\e=(\lambda_\phi(X_\infty)-\lambda_\phi(X_\infty'))/2$.

Since $X_\infty$ has not jumped, by Theorem~\ref{thmjump} we have
$\lambda(\phi|_A)\leq \lambda_\phi(X_\infty)$. If $\lambda(\phi|_A)< \lambda_\phi(X_\infty)$,
let $A'\in\O(A)$ be a point such that $\lambda_{\phi|_A}(A')<\lambda_\phi(X_\infty)$. Now
Lemma~\ref{lemma9} provides an element of $\O(\Gamma)$ which is displaced less or equal than
$\max\{\lambda_{\phi|_A}(A'),\lambda_\phi(X_\infty')+\e\}$, contradicting the fact that
$X_\infty$ is a minpoint for $\lambda$. Therefore $\lambda(\phi|_A)=\lambda(X_\infty)$.

By Theorem~\ref{thmminptE} there is
$A_\infty\in\overline{\O(A)}^\infty$ such that
$\lambda_{\phi|_A}(A_\infty)=\lambda(\phi|_A)$ and which has not jumped in $\O(A)$.
Thus $A_\infty$ is obtained, without jumps, from a point
$A'_\infty\in\O(A)$ by collapsing a (possibly empty) invariant core sub-graph $B$. So $A_\infty\in\O(A'_\infty/B)$.

Let $Y$ be a $\Gamma$-graph obtained by inserting a copy of $A_\infty'$ in $X'_\infty$. Let
$Y'$ be the graph obtained collapsing $B$. $Y'$ belongs to the simplicial boundary of
$\Delta_Y$ and, since $A_\infty$ has no jump, then so does $Y'$. Now, observe that
$Y'\in\O(Y/B)$ and $A_\infty$ is a $\phi$-invariant subgraph of $Y'$ so that $Y'/A_\infty=X_\infty'$.
Lemma~\ref{lemma9}
provides an element in $Y'_\infty\in\O(Y/B)$, in the same simplex of $Y'$ which is displaced no more than
$\lambda_{\phi|_A}(A_\infty)$ (because
$\lambda_\phi(X_\infty')<\lambda_\phi(X_\infty)=\lambda_{\phi|_A}(A_\infty)$).
Now, $Y'_\infty$ is a new minpoint for $\lambda$ with
$\dim(\Delta_{Y'_\infty})>\dim(\Delta_{X_\infty})$ contradicting the maximality hypothesis on
$X_\infty$. It follows that $X_\infty$ is a train track point in $\O(X_\infty)$ as desired.\qed

So we have seen that, even if non-jumping min-points are not necessarily train tracks, some of
them are. Conversely, we see now non-jumping train tracks at the bordification are always min-point for $\lambda_\phi$.

\begin{lem}\label{triangin}
Let $\phi\in\Aut(\Gamma)$ and let $X\in\O(\Gamma)$. If there is $k$ so that there is a constant
$A>0$ such that for any $n>>1$
$$Ak^n\leq \Lambda(X,\phi^n X)$$ then $k\leq \lambda(\phi)$.
\end{lem}
\proof This follows from the multiplicative triangular inequality.
For any $Y\in\O(\Gamma)$ we have $\Lambda(Y,\phi^nY)\leq \Lambda(Y,\phi Y)^n$. Define a constant
$C=\Lambda(X,Y)\Lambda(Y,X)$ and notice that we also have $C=\Lambda(X,Y)\Lambda(\phi
Y,\phi X)$. Then,
$$Ak^n\leq\Lambda(X,\phi^n X)\leq\Lambda(X,Y)\Lambda(Y,\phi^n Y)
\Lambda(\phi^nY,\phi^n X)\leq C\Lambda(Y,\phi Y)^n$$
whence, for any $n$ $$\left(\frac{k}{\Lambda(Y,\phi Y)}\right)^n\leq \frac{C}{A}.$$
This implies $k\leq \Lambda(Y,\phi(Y))$. By choosing a minimizing sequence of points
$Y_i$ we get $k\leq \lambda(\phi)$.\qed

\begin{lem}\label{dom5}
  Let $\phi\in\Aut(\Gamma)$. Let $X_\infty\in\overline{\O(\Gamma)}$ which has not jumped.
  Suppose that there is a loop $\gamma \in X_\infty$ and $k>0$ such that
  $L_{X_\infty}(\phi^n)(\gamma)\geq k^n L_{X_\infty}(\gamma)$. Then
  $$k\leq\lambda(\phi).$$

  In particular, if $X_\infty$ is a train track for $\phi$ as an element of $\Aut(X_\infty)$,
  then it is a minpoint  for $\phi$ as an element of $\Aut(\Gamma)$.
\end{lem}
\proof Let $X\in\O(\Gamma)$ so that $X_\infty$ is obtained from $X$ by collapsing a core sub-graph
$A\subset X$. Let $X_\e$ be a point of $X$ where $\vol(A)<\e$.
Let $\gamma$ be as in the hypothesis.
For $\e$ small enough we have
$L_{X_\e}(\gamma)\leq 10 L_{X_\infty}(\gamma)$, and therefore
\begin{eqnarray*}
\Lambda(X_\e,\phi^nX_\e) \geq  \frac{L_{X_\e}(\phi^n\gamma)}{L_{X_\e}(\gamma)}\geq
\frac{L_{X_\infty}(\phi^n\gamma)}{10L_{X_\infty}(\gamma)}\geq
\frac{k^nL_{X_\infty}(\gamma)}{10L_{X_\infty}(\gamma)}=
\frac{k^n}{10}.
\end{eqnarray*}
By Lemma~\ref{triangin} we have $\lambda(\phi)\geq k$.

For the second claim it suffice to choose let $\gamma$ a legal candidate that
realizes $\Lambda(X_\infty,\phi X_\infty)$. So
$L_{X_\infty}(\phi^n(\gamma))=\lambda_\phi(X_\infty)^nL_{X_\infty}(\gamma)$.

Hence $\lambda(\phi)\geq \lambda_\phi(X_\infty)$ and since $X_\infty$
has not jumped $\lambda(\phi)\leq\lambda_\phi(X_\infty)$.\qed

We are now in position of extending the second claim of Theorem~\ref{newjump}.

\begin{cor}\label{corlalx}
  Let $\phi\in\Aut(\Gamma)$. Let $X\in\O(\Gamma)$ and $X_\infty$ be obtained from $X$ by
  collapsing a $\phi$-invariant core sub-graph $A$. Then $$\lambda(\phi|_A)\leq \lambda(\phi).$$

Moreover, if $\lambda(\phi|_A)=\lambda_\phi(X_\infty)$, then $$\lambda(\phi)=\lambda(\phi|_A).$$

In particular $X_\infty$ has not jumped
if and only if $$\lambda(\phi)\leq \lambda(X_\infty).$$
\end{cor}
\proof Let $\lambda=\lambda(\phi|_A)$.
By  Lemma~\ref{lemmattE} and Theorem~\ref{Theoremtt}, there is $\bar A\in\overline{\O(A)}^\infty$ which is a min-point for
$\phi|_A$, which has not  jumped in $\O(A)$, and which is a train track for $\phi|_A$ as an
element of $\Aut(\bar A)$. Let $f_A$ be a
train track map $f_A:\bar A\to \bar A$ representing $\phi|_A$.
Therefore, there is a legal loop $\gamma$ in $\bar A_{\max}$ whit legal images in $\bar A_{\max}$ and
stretched exactly by $\lambda$. Let now $\widehat X$ be a metric $\Gamma$-graph obtained by
inserting a copy of $\bar A$ in $X_\infty$. Since $\overline A$ has not jumped in $\O(A)$, then
$\widehat X$ has not jumped in $\O(\Gamma)$.

Let $f:\widehat X\to\widehat X$ be any $\PL$-map
representing $\phi$ so that $f|_A=f_A$. Therefore $f_A^n(\gamma)$ is immersed for any $n$ and the length of $f_A^n(\gamma)$ is $\lambda^n$
times the length of $\gamma$. It follows that $L_{\widehat X}((\phi^n)\gamma)=\lambda^n(L_{\widehat
  X}(\gamma))$.

By Lemma~\ref{dom5}
$\lambda(\phi|_A)=\lambda\leq \lambda(\phi)$, and the first claim is proved. Moreover, if
$\lambda(\phi|_A)=\lambda_\phi(X_\infty)$, then
$$\lambda(\phi)\leq\lambda(X_\infty)=\lambda(\phi|_A)=\lambda\leq\lambda(\phi)$$ and therefore
all inequalities are equalities.
Finally, if $X$ has no jumped then $\lambda(X)\geq \lambda(\phi)$ just because this inequality is true by definition for points in
$\O(\Gamma)$ and clearly passes to limits of non-jumping sequences, and the converse inequality
follows from the second claim and Theorem~\ref{thmjump}.
\qed

\medskip
Note that Corollary~\ref{corlalx} implies that {\em a posteriori} we can remove the non-jumping
requirement from Theorem~\ref{thmminptE} and Lemma~\ref{lemmattE}.

\begin{cor}[Min-points don't jump]
  Let $\phi$ any element in $\Aut(\Gamma)$. If $X\in\overline{\O(\Gamma)}^\infty$
  is such that $\lambda_\phi(X)=\lambda(\phi)$, then it has not jumped.
\end{cor}
\proof This is a direct consequence of Corollary~\ref{corlalx}.\qed

\medskip

We introduce the notion of train track at infinity.

\begin{defn}[Train track at infinity]\label{dttinfty}
  Let $\phi\in\Aut(\Gamma)$. The set $\TT^\infty(\phi)$ is defined as the set of points
  $X\in\overline{\O(\Gamma)}^\infty$ such that  $X$ has not jumped, and $X$ is a train track
  point for $\phi$ in $\O(X)$. (Hence $\lambda_\phi(X)=\lambda(\phi)$ by Lemma~\ref{dom5}.)
\end{defn}

Note that $\TT(\phi)\subset\TT^{\infty}(\phi)$. The main differences are that $\TT(\phi)$ may
be empty (if $\phi$ is thin) while any $\phi$ has a train track in $\TT^\infty(\phi)$. On the
other side, $\TT(\phi)$ coincides with the set of minimally displaced points, while
$\TT^\infty(\phi)$ may be strictly contained in the set of minimally displaced points.

With this definition  we can collect some of the above results
in the following simple statement, which is a straightforward
consequence of Theorems~\ref{Theoremtt},~\ref{thmminptE} and Lemmas~\ref{lemmattE},~\ref{dom5}.
\begin{thm}\label{corttE}
  For any $\phi \in\Out(\Gamma)$, $TT^\infty(\phi)\neq\emptyset$. For any
  $X\in\TT^{\infty}(\phi)$, $\lambda_\phi(X)=\lambda(\phi)$.
\end{thm}

The following corollary shows that if $\phi$ is reducible then there is a train track showing reducibility.

\begin{cor}[Detecting reducibility]\label{corred}
  Let $\phi\in\Aut(\Gamma)$ be reducible. Then there is $T\in\TT^\infty(\phi)$ such that either
  $T\in\partial_\infty\O(\Gamma)$ or there is an optimal
  map $f_T:T\to T$ representing $\phi$ such that there is a proper sub-graph of $T$ which is $f_T$-invariant.
\end{cor}
\proof Since $\phi$ is reducible there is  $X\in\O(\Gamma)$, a $\PL$-map $f:X\to X$
representing $\phi$ and a proper sub-graph $A\subset X$ such that $f(A)=A$. We can therefore
collapse $A$ and $\lambda$ won't explode. By Theorem~\ref{corttE} there is a train track $Z$ for
$\phi$ in $\overline{\O(X/A)}^\infty$ and a train track $Y$ for $\phi|_A$ in
$\overline{\O(A)}^\infty$.  If
$\lambda_{\phi|_A}(Y)\leq\lambda_\phi(Z)$, then
$Z\in\TT^\infty(\phi)\cap\partial_\infty\O(\Gamma)$ and  we are done.
Otherwise, since $Z$ has not jumped (as a point of $\partial_\infty\O(X/A)$),
we can regenerate it to a point $Z'\in\O(X/A)$ with
$\lambda_\phi(Z')<\lambda_{\phi|_A}(Y)$. We now apply regeneration Lemma~\ref{lemma9} to $Y$ and
$Z'$. If $Y\in\partial_\infty\O(A)$, then we get a train track for $\phi$ in
$\partial_\infty\O(\Gamma)$. If $Y\in\O(A)$ we get a train track for $\phi$ in $\O(\Gamma)$
admitting $Y$ as an invariant sub-graph.\qed

\medskip

In fact, the proof of Corollary~\ref{corred} proves more: that train tracks detect any
invariant free factor.
\begin{cor}[Strong reformulation of Corollary~\ref{corred}]\label{strongcorred}
  Let $\phi\in\Aut(\Gamma)$. Let $X$ be a $\Gamma$-graph having a $\phi$-invariant core
  sub-graph $A$. Then there is $Z\in\O(X/A)$ and $W\in\Hor(Z)$ such that the simplex $\Delta_W$
  contains a minimizing sequence for $\lambda$. Moreover if $Y\in\O(A)$ is the graph used to
  regenerate $W$ from $Z$, then the minimizing sequence can be chosen with $\PL$-maps $f_i$
  such that $f_i(Y)=Y$ and $\Lip(f_i)\to\lambda(\phi)$.
\end{cor}
\proof Follows from the proof of Corollary~\ref{corred} (and Lemma~\ref{lemma9}).\qed

\medskip

Finally, as in the case of irreducible automorphisms, the existence of train tracks gives the
following fact.
\begin{cor}
  For any $\phi\in\Aut(\Gamma)$ we have $\lambda(\phi^n)=\lambda(\phi)^n.$
\end{cor}
\proof It follows from Theorem~\ref{corttE} and Lemma~\ref{Lemmatt2}.\qed
 
\section{Statement of main theorem and regeneration of paths in the bordification}
We use Notation~\ref{not:gamma}, that we recall here for the benefit of the reader.
\begin{itemize}
\item $G$ will always mean a group with  a splitting $\mathcal G:G=G_1*\dots*G_p*F_n$;
\item  $\Gamma=\sqcup \Gamma_i$ will always mean that $\Gamma$ is a finite disjoint
union of finite graphs of groups $\Gamma_i$, each with trivial edge-groups and non-trivial fundamental group $H_i=\pi_i(\Gamma_i)$, each
$H_i$ being equipped with the splitting given by the vertex-groups.
\end{itemize}
Also, we recall that we defined the rank of $\Gamma$ in Definition~\ref{pr4_rank}.
Finally, we recall the notation for $\lambda$ (Definition~\ref{defdispl}).
  For any automorphism $\phi \in\Out(\Gamma)$ we define the function
$$\lambda_\phi:\O(\Gamma)\to\R\qquad \lambda_\phi(X)=\Lambda(X,\phi X)$$
If $\Delta$ is a simplex of $\O(\Gamma)$ we define
$$\lambda_\phi(\Delta)=\inf_{X\in\Delta}\lambda_\phi(X)$$
If there is no ambiguity we write simply $\lambda$ instead of $\lambda_\phi$.
Finally, we set
$$\lambda(\phi)=\inf_{X\in\O(\Gamma)}\lambda_\phi(X)$$
We agree that we extend the function $\lambda$ to points in
$X_\infty\in\partialì_\infty(\O(\Gamma))$ for which there is a sequence of points
$X_i\in\O(\Gamma)$ such that $X_i\to X_\infty$ with $\lambda(X_i)$ bounded above, and we set
$\lambda=\infty$ on other points. (Definition~\ref{dli}.)

\begin{defn}\label{def:sp}
  Let $X,Y\in\overline{\O(\Gamma)}^\infty$. A {\em simplicial path} between $X,Y$ is given by:
  \begin{enumerate}
  \item A finite sequence of points $X=X_0,X_1,\dots,X_k=Y$, called vertices, such that
    $\forall i=1,\dots, k$,  there is a minimal simplex $\Delta_i$ such that
    $\Delta_{X_{i-1}}$ and $\Delta_{X_i}$ are  both simplicial  faces of $\Delta_i$ (we allow
    one of them or even both to coincide with $\Delta_i$).
  \item Euclidean segments $\overline{X_{i-1}X_i}\subset \Delta_i$, called edges.
 \end{enumerate}
\end{defn}

\begin{defn}
  We say that a set $\chi$ is {\em connected by simplicial paths} if for any $x,y\in\chi$ there is a
  simplicial path between $x$ and $y$ which is entirely  contained in $\chi$.
\end{defn}

\begin{thm}[Level sets are connected]\label{tconnected}
Let $\phi \in\Out(\Gamma)$. For any $\e > 0 $ the set
$$\{X\in\O(\Gamma):\lambda_\phi(X)\leq \lambda(\phi) + \e \}$$ is connected in $\O(\Gamma)$ by simplicial paths.
The set
$$\{X\in\overline{\O(\Gamma)}^\infty:\lambda_\phi(X)=\lambda(\phi)\}$$
is connected by simplicial paths in $\overline{\O(\Gamma)}^\infty$.
\end{thm}
The remaining goal of the paper is devoted to the proof of Theorem~\ref{tconnected}. The rough strategy is
to prove the second claim and then prove that paths in the bordification can regenerate to
paths in $\O(\Gamma)$ without increasing $\lambda$ too much. The proof goes by induction on the
rank of $\Gamma$ (see definition below).

\begin{rem}
  Theorem~\ref{tconnected} is trivially true if $\rank(\Gamma)=1$, because in that case
  either $\O(\Gamma)$ or $\mathbb P\O(\Gamma)$ is a single point.
\end{rem}

\begin{lem}[Regeneration of segments]\label{regseg}
Fix $[\phi]\in\Out(\Gamma)$. Let $X_\infty,Y_\infty\in\overline{\O(\Gamma)}^\infty$ such that
$\Delta_{Y_\infty}$ is a (not necessarily proper) simplicial face of
$\Delta_{X_\infty}$. Suppose that  $\lambda(X_\infty)\geq\lambda(\phi)$.
Then there is an open simplex $\Delta$ of $\O(\Gamma)$ such
that for any $\e>0$ there is  $Y\in \Hor(Y_\infty)\cap \overline\Delta$ and
$X\in\Hor(X_\infty)\cap\Delta$
 such that $$\lambda_\phi(Y),\lambda_\phi(X)<\max\{\lambda_\phi(Y_\infty),\lambda_\phi(X_\infty)\}+\e.$$
Moreover, such inequality holds on the whole segments $\overline{XX_\infty}$ and $\overline{YY_\infty}$.
\end{lem}
\proof Let $X_\infty$ be obtained by collapsing a $\phi$-invariant core-subgraph $A$ from a
$\Gamma$-graph $\widehat X$. Since $\lambda_\phi(X_\infty)\geq\lambda(\phi)$, by Corollary~\ref{corlalx} $\lambda(\phi|_A)\leq\lambda_\phi(X_\infty)$.
By Theorem~\ref{thmminptE} there is a simplex in $\O(A)$ that contains a minimizing sequence
for $\lambda(\phi|_A)$. Let $A_\e$ be a point in that simplex such that
$\lambda(A_\e)<\lambda(\phi|_A)+\e$. The required simplex $\Delta$ is
obtained by inserting a copy of $A_\e$ at the place of $A$ in $X_\infty$. We notice that such
a $\Delta$ is not unique.
 By Lemma~\ref{lemma9} there is a
point $X\in\Delta\cap\Hor(X_\infty)$ such that
$\lambda_\phi(X)\leq\lambda_\phi(X_\infty)+\e$.

Let's now see what happens to the points in $\overline\Delta\cap \Hor(Y_\infty)$.
By hypothesis there is a $\phi$-invariant $B\subseteq X_\infty$ such that
as a graph (i.e. forgetting the metrics), $Y_\infty$ is obtained from $X_\infty$ by collapsing
$B$. $B$  has a pre-image in $X$ still denoted by $B$. Let $T$ be the forest $(A\cup
B)\setminus \core(A\cup B)$. If $Y'=X/T$, as a graph, $Y_\infty=X/(A\cup B)=Y'/\core(A\cup B)$.

Thus the finitary face $\Delta_{Y'}$ of $\Delta$ obtained by the collapse of $T$ intersects
$\Hor(Y_\infty)$.

Let $f:X\to X$ be an optimal map representing $\phi$. Since $\core(A\cup B)$ is $\phi$-invariant,
$f(\core (A\cup B))\subset \core(A\cup B)$ up to homotopy. It follows that there is
a $\PL$-map $g:\core(A\cup B)\to \core(A\cup B)$ representing $\phi|_{A\cup B}$ such that
$\Lip(g)\leq\lambda_\phi(X)\leq\lambda_\phi(X_\infty)+\e$.
By Lemma~\ref{lemma9} there is a point $Y\in\Hor(Y_\infty)\cap \Delta_{Y'}$ such that
$\lambda_\phi(Y)\leq\max\{\lambda_\phi(Y_\infty)+\e,\Lip(g)\}\leq\max\{\lambda_\phi(Y_\infty)+\e,\lambda_\phi(X_\infty)+\e\}$. The
last claim also follows by Lemma~\ref{lemma9}. \qed

\medskip

Now we can plug in the inductive hypothesis in the proof of Theorem~\ref{tconnected}. Recall
that if $X=T/S$ as graphs of groups, then we  denote by $\pi:\Hor(X)\to \mathbb P\O(S)$ the
projection that associates to a point in $\Hor(X)$ its collapsed part (see section~\ref{sechor}).

\begin{lem}[Regeneration of horospheres]\label{inhor}
   Suppose that Theorem~\ref{tconnected} is true in any rank less than $\rank(\Gamma)$. Let
  $\phi \in\Out(\Gamma)$. Let $T\in\O(\Gamma)$ be a $\Gamma$-graph having a proper
  $\phi$-invariant core sub-graph
  $S$. Let $X\in\partial_\infty\O(\Gamma)$ be the graph obtained from $T$ by collapsing $S$, and
  let $A,B\in\Hor(X)\subset\O(\Gamma)$. Let $m_A$
  and $m_B$ be the supremum of $\lambda_\phi$ on the Euclidean segments $\overline{AX}$ and
  $\overline{BX}$ respectively. Then, for any $\e>0$ there is a simplicial path $\gamma$
  between $A$ and $B$, and  in $\Hor(X)$, such that for any vertex $Z$ of $\gamma$ we have $$\lambda_\phi(Z)<\max\{m_A,m_B\}+\e.$$
\end{lem}
\proof Let $L=\max\{m_A,m_B\}$. Since $S$ is $\phi$-invariant, by Lemma~\ref{fatto1} we have
that $\lambda_\phi(X)$ is finite and by Lemma~\ref{lemma9} both $m_A$ and $m_B$ are finite.

For any $Y\in\Hor(X)$, Theorem~\ref{sausagelemma} implies $\lambda_\phi(\pi(Y))\leq\lambda_\phi(Y)$ so
$$\lambda_\phi(\pi(A))\leq\lambda_\phi(A)\quad
\lambda_\phi(\pi(B))\leq\lambda_\phi(B)
$$
hence, $\lambda_\phi(\pi(A)),\lambda_\phi(\pi(B))\leq L$. The rank of $S$ is strictly
smaller than $\rank(\Gamma)$ because it is a proper sub-graph of $T$. Hence Theorem~\ref{tconnected} holds for $\O(S)$. So there
is a finite simplicial path $(Y_i)\in\O(S)$ between $\pi(A)$ and $\pi(B)$ such that $\lambda_\phi(Y_i)<L+\e$. Then, there is
a finite  simplicial path in $\Hor(X)$ between $A$ and $B$
whose vertices are points $\widehat T_j$ such that for
any $j$ there is $i$ such that $\pi(\widehat T_j)=Y_i$. By Lemma~\ref{lemma9} there is a
simplicial path in $\Hor(X)$ whose vertices are points $Z_j\in\Delta_{\widehat T_j}$ such that $\pi(Z_j)=\pi(\widehat
T_j)=Y_i$ and $\lambda_\phi(Z_j)<L+\e$.\qed

\medskip

We recall that we are using the notation of Definition~\ref{def:sp}.
\begin{thm}[Regeneration of paths]\label{tregge}
  Suppose that Theorem~\ref{tconnected} is true in any rank less than $\rank(\Gamma)$. Let
  $\phi \in\Out(\Gamma)$. Let
  $\gamma=(X_i)$ be a simplicial path in $\overline{\O(\Gamma)}^\infty$ such that
  for every $i$ either  $\Delta_{X_{i-1}}$ is a simplicial face of $\Delta_{X_i}$ vice versa,

  Suppose that there is $L$ so that for any point $X_i$ we have
  $$\lambda(\phi)\leq\lambda_\phi(X_i)\leq L.$$

  Then, for any $\e>0$ there exists a simplicial path $\eta$ in $\O(\Gamma)$, contained in the
  level set $\lambda_\phi^{-1}(L+\e)$, and such that each
  vertex of $\eta$ belongs to the horosphere of some $X_j$.
\end{thm}
\proof By Lemma~\ref{lconvexity} it suffices to define the vertices of the path $\eta$.
By Lemma~\ref{regseg} For any $i$ there
are points $A_i,B_i\in\Hor(X_i)$ such that $\lambda_\phi(A_i),\lambda_\phi(B_i)\leq L+\e$ and
such that $B_i,A_{i+1}$ are in the same closed simplex of $\O(\Gamma)$. By Lemma~\ref{inhor}
there is a simplicial path $Y_{ij}$ between $A_i$ and $B_i$ such that $Y_{ij}\in\Hor(X_i)$ and
$\lambda_{\phi}(Y_{ij})\leq L+\e$. The path $\eta$ is now defined by the concatenation of such
paths and the segments $\overline{B_iA_{i+1}}$. \qed

\medskip
The proof of Theorem~\ref{tconnected} now continues by an argument of peak reduction among
simplicial paths connecting two points in the same level set. In next section we prove
the results that will allow to reduce peaks.

\section{Preparation to peak reduction}

We keep Notation~\ref{not:gamma}. We also recall that for $\phi\in\Aut(\Gamma)$ and a simplex
$\Delta\in\overline{\O(\Gamma)}^\infty$ we are using the notation $$\lambda(\Delta)=\lambda_\phi(\Delta)=\inf_{X\in\Delta}\lambda_\phi(X).$$

For the remaining of the section we fix $\phi\in\Aut(\Gamma)$.
Recall that we are using the notation of Definition~\ref{def:sp} for simplicial paths. In
Theorem~\ref{tregge} we required that given two consecutive points
  $X_{i-1},X_i$ then one of $\Delta_{X_{i-1}},\Delta_{X_i}$ is a face of the other. In general
  such a  condition is easy to obtain by adding a middle point, but we need to do it in such a
  way to control the function $\lambda$, which is not in general continuous on
  $\overline{\O(\Gamma)}^\infty$.

We describe now a procedure for locally minimizing $\lambda$ on simplicial path in
$\O(\Gamma)$.

Let $(X_i)_{i=0}^{k}$ be a simplicial path such that:
\begin{itemize}
\item $X_2,\dots,X_{k-1}\in\O(\Gamma)$;
\item If $X_0\notin\O(\Gamma)$ then $X_0\in\partial_\infty\Delta_{X_1}$ and has no jump in
  $\Delta_{X_1}$.
\item If $X_k\notin\O(\Gamma)$ then $X_k\in\partial_\infty\Delta_{X_{k-1}}$ and has no jump in $\Delta_{X_{k-1}}$.
\end{itemize}

Then, we define a new simplicial path by doing the following steps:
\begin{enumerate}
\item For any $i$, if $X_{i-1}$ and $X_{i}$ are both proper faces of $\Delta_i$, then we add to
  the path a   point $\widehat X_i\in\Delta_i$.
\item  We renumber the sequence of vertices, still
  denoted by $X_i$. So now the sequence is $(X_i)_{i=0}^m$ for some $m\geq k$.
\item We set $Y_0=X_0$ and $Y_m=X_m$.
\item For any any $i$, we chose $Y_i\in\overline{\Delta_{X_i}}^\infty$ so that
  $\lambda(Y_i)=\lambda(\Delta_{X_i})=\lambda(\Delta_{Y_i})$, moreover we require that
   $Y_i$ maximizes the
  dimension of $\Delta_{Y_i}$ among such points. Such a $Y_i$ exists and has not jumped in
  $\Delta_{X_i}$ by  Corollary~\ref{fatto2}.
\item If it happens that $Y_i=Y_{i+1}$ for some $i$, then we identity such points and we
  renumber the sequence accordingly.
\end{enumerate}

\begin{defn}
  We say that the path $(Y_i)$ as above is obtained by {\em minimizing} the path $(X_i)$. A
  path in $\overline{\O(\Gamma)}^\infty$ is said {\em minimized} if it is the optimization of a
  simplicial path $(X_i)_{i=0}^k$ contained in $\O(\Gamma)$ except at most at its endpoints $X_0,X_k$,
  which have no jump in $\Delta_{X_0}$ and $\Delta_{X_{k-1}}$ respectively.
\end{defn}

Note that if $(Y_i)$ has the same endpoints of $(X_i)$, and any other vertex minimizes
$\lambda$ in its simplex. In particular
$$\sup_i\lambda(Y_i)\leq\sup_i\lambda(X_i).$$

\begin{lem}\label{lkey}
 Let $A,B$ two consecutive vertices of a minimized simplicial path $\Sigma$.
 \begin{itemize}
 \item Por any point $P$ of $\overline{AB}$ we have
 $\lambda(P)\geq\lambda(\phi)$;
\item if  $\lambda(A)=\lambda(B)$, then
$\lambda$ is constant on the segment $\overline{AB}$;
\item if $\lambda(A)>\lambda(B)$ then $\lambda$ is continuous and strictly monotone near $A$. \end{itemize}
 \end{lem}
\proof
By definition of minimized path, there is a sequence $A_i\to A$ and $B_i\to B$,
both without jump and in the same closed simplex $\overline{\Delta}$, and one of them is in the
open simplex $\Delta$ of $\O(\Gamma)$. Without loss of generality we may  assume
$\Delta=\Delta_{B_i}$. Clearly $\lambda(A),\lambda(B)\geq\lambda(\Delta)\geq\lambda(\phi)$
because they did not jumped in $\Delta$. If
$B\in\Delta$ then $P\in\Delta$, and $$\lambda(P)\geq\lambda(\Delta)\geq\lambda(\phi).$$

If $B\notin \Delta$, by
Corollary~\ref{fatto2} no point of $\overline{\Delta}^\infty$ has jumped in $\Delta$, so $P$
has not jumped and again $\lambda(P)\geq\lambda(\Delta)\geq\lambda(\phi).$

Suppose now that $\lambda(A)=\lambda(B)=L$. Since $B$ minimizes $\lambda$ on $\Delta$ we have
$L=\lambda(\Delta)$. By quasi convexity
$$\lambda(\Delta)\leq\lambda(P)\leq\lambda(B)=\lambda(\Delta).$$
The last claim follows directly from Lemma~\ref{lconvexity} and Lemma~\ref{lconv2}.
\qed

\medskip
Note that in particular, Lemma~\ref{lkey} implies that
$\lambda$ is bounded from below by $\lambda(\phi)$ on any minimized simplicial path.

\begin{defn}
  A simplicial path $\Sigma$ is said {\em $L$-calibrated} if $\lambda$ is
  continuous on $\Sigma$ and for any point $P$ of $\Sigma$ we have
  $$\lambda(\phi)\leq \lambda(P)\leq L,$$
and if max-points minimize $\lambda$ on their simplices. (That is to say, if $X$ is such that $\lambda(X)=\max_P\in\Sigma\lambda(P)$, then $\lambda(X)=\lambda_{\Delta_X}$.)
\end{defn}
Note that by Lemma~\ref{lconvexity} if $A,B$ are consecutive vertices of a calibrated path
such that $\lambda(A)>\lambda(B)$ then $\lambda$ is strictly monotone near $A$.

\begin{lem}\label{calibration}
  Let $\Sigma$ be a minimized simplicial path and let $$L=\max_{P\in\Sigma}\lambda(P).$$
Then there is a $L$-calibrated simplicial path obtained by
  adding some extra vertices to $\Sigma$.
\end{lem}
\proof By Lemma~\ref{lkey} we have to care only about continuity.
Let $A,B$ be two consecutive points of $\Sigma$. By Lemma~\ref{lkey} $\lambda$ can be
not continuous only at the endpoint with lower $\lambda$. Suppose $\lambda(A)>\lambda(B)$.
By definition of minimized path there is a sequence $A_i\to A$ and $B_i\to B$,
both without jump and in the same closed simplex $\overline{\Delta}$, and one of them is in the
open simplex $\Delta$ of $\O(\Gamma)$. Moreover, since $\lambda(B)=\lambda(\Delta_{B_i})$,
$\lambda(A)=\lambda(\Delta_{A_i})$ and  one of them equals $\lambda(\Delta)$,
$\lambda(B)<\lambda(A)$ forces $\Delta=\Delta_{B_i}$. If $B\in\Delta$ then $\lambda$ is
continuous at $B$ because $\lambda$ is continuous on any open simplex $\Delta$ of
$\O(\Gamma)$. If $B\notin\Delta$, then $B$ has not jumped in $\Delta$. Therefore (by
Theorem~\ref{newjump} and Lemma~\ref{lemma9}), there is a point
$\widehat B\in\Delta$ such that

\begin{itemize}
\item $\lambda(\widehat B)<\lambda(A)$;
\item $\lambda$ is continuous on the segment $\overline{\widehat B B}$.
\end{itemize}
by Lemma~\ref{lconv2} $\lambda$ is continuous on the segment $\overline{A\widehat B}$.
Since we did not modified $\Sigma$ at its max-points, also the second condition of calibration
is assured because it is already satisfied for minimized paths.
\qed

\medskip

We prove now a (technical) fact that  can be informally
phrased as follows\footnote{We recall
  that by definition $\overline{\O(\Gamma)}^\infty=\overline{\O(\Gamma)}$ and that the symbol
  $\infty$ is just to put emphasis on the fact that we are considering the simplicial
  bordification of the outer space obtained by adding all simplices at infinity.}:

\centerline{\parbox{0.8\textwidth}{Given $X\in\overline{\O(\Gamma)}^\infty$ and $f:X\to X$
    an optimal map representing $\phi$, if $Y$ is sufficiently close to $X$ for the Euclidean
    metric, then any fold in $X$ directed by $f$ closely reads in $Y$.}}

\begin{thm}\label{thmZnew}
Let $X,Y\in\overline{\O(\Gamma)}$. Suppose that $\Delta_X$ is a simplicial face of
$\Delta_Y$. Thus as graphs, $Y$ is obtained by collapsing a sub-graph $A$. Suppose that
$\core(A)$ is $\phi$-invariant. For $t\in[0,1]$ let $Y_t=(1-t)X +t Y$ be a
parametrization of the Euclidean segment from $X$ to $Y$. Let $\sigma_t:Y_t\to X$ be the map
obtained by collapsing  $A$ and by linearly rescaling the edges in $Y\setminus A$.

Let $f:X\to X$ be an optimal map representing $\phi$. Then for any $\e >0$ there is $t_\e>0$ such that $\forall 0\leq t<t_\e$ there is an optimal map
$g_t:Y_t\to Y_t$ representing $\phi$ such that $$d_\infty(\sigma_t\circ g_t,f\circ \sigma_t)<\e.$$
\end{thm}
\proof The proof of this theorem relies on accurate (but boring) estimates. For the happiness of
the reader we postpone the proof to the appendix.\qed

\begin{cor}\label{corpeak}
In the hypotheses, and with notation  of Theorem~\ref{thmZnew}, let $\tau$ be an $f$-illegal
turn of $X$ and let
$\Delta^\tau$ be the simplex obtained from $\Delta_X$ by fold a little that turn. Then, for any
$\e>0$, there is $t_\e$ such that $\forall t<t_\e$, there is a finite  simplicial path
$\Sigma_t$ in $\O(Y)$ with vertices
$Z_0^t=Y_t,Z_1^t,\dots,Z_m^t$ such that $\Delta_{Z_i^t}$ has $\Delta_X$ as a simplicial
face for $i\neq m$, $\Delta_{Z_m^t}$ has $\Delta^\tau$ as a simplicial
face, and such that for any point $Z$ of $\Sigma_t$
we have
$$\lambda(X)-\e<\lambda(Z)\leq\lambda(Y_t).$$
Moreover, for $s\in[0,t]$, $s\mapsto Z_i^s$ parametrizes the segment from $X$ to $Z_i^t$, and $Z_m^s$
that from $X^\tau$ to $Z_m^t$.
\end{cor}
\proof  For this proof we will work entirely with trees. So $Y$ will denote a $\Gamma$-trees,
$A$ an equivariant family of sub-trees, and so on.

We denote by $A_t$ the metric copy of $A$ in $Y_t$.
By hypothesis there are two different segments $\alpha_\tau,\beta_\tau$ incident at the same vertex $v$ in
$X$ such that $f$ overlaps $\alpha_\tau$ and $\beta_\tau$. If $v\notin\sigma_t(A_t)$ then, for any small enough
$\e$ and $t<t_\e$, also $g_t$ must overlap $\alpha=\sigma_t^{-1}(\alpha_\tau)$ and
$\beta=\sigma_t^{-1}(\beta_\tau)$, and the claim follows by (equivariantly) performing the
corresponding simple fold directed by $g_t$. The inequality ``$\leq\lambda(Y_t)$'' follows
because the fold is directed by an optimal map, the inequality ``$>\lambda(X)-\e$'' follows by
lower semicontinuity of $\lambda$.

Otherwise, $\alpha$ and $\beta$ are segments incident to the
same component of $A_t$. If $\alpha$ and $\beta$ are incident to the same point, then we
proceed as above, so we can suppose that they are incident to different points of $A$.

For small enough $\e$ and $t<t_\e$ we
have that $g_t$ overlaps some open sub-segments of $\alpha$ and $\beta$. Let $a\in\alpha$
and $b\in\beta$ such that $g_t(a)=g_t(b)$ and such that $a$ is the closest possible to $A$.

Let $\gamma$ be the shortest path from $\alpha$ and $\beta$ in $A_t$.
It turns out that  $\gamma$ is a simple simplicial path. On $\gamma$ we put an extra simplicial
structure given by the pull-back via $g_t$: we declare new vertices of $\gamma$ the points
whose $g_t$-image is a vertex of $Y_t$.
$g_t(\gamma)$ is a tree because $Y_t$ is. Moreover, since $g_t (a)=g_t(b)$,
the restriction of $g_t$ to $\gamma$ cannot be injective. In particular, if $x\in\gamma$ is a
point such that $d_{Y_t}(g_t(x),g_t(a))$ is maximal, then $x$ is a vertex of $\gamma$, and the two
sub-segments of $\gamma$ incident to $x$ are completely overlapped.

Let $Z_1^t$ be the tree obtained by equivariantly identify such segments. Clearly, $g_t$ induces
a map $g_t^1:Z_1^t\to Z_1^t$. Such map is continuous and not necessarily $\PL$. However,
$$\Lip(g_t^1)\leq\Lip(g_t)$$ and $\PL(g_t^1)$ still represents $\phi$. Since
$\Lip(\PL(g_t^1))\leq\Lip(g_t^1)$ he have
$$\lambda(Z_1^t)\leq\lambda(Y_t).$$
Note also that $\Delta_{Z_1^t}$ has $\Delta_X$ as a simplicial face because our identification
occurred in $A_t$. Also, since $Y_t$ parametrizes the segment from $X$ to $Y$, as $t$ varies
$Z_1^t$ parametrizes the segment from $X$ to $Z_1^t$.

Note that a priori we may have $\Delta_{Z_1^t}=\Delta_Y$, but in any case
$\Delta_{Z_1^t}$ is either a (non necessarily proper)
simplicial face of $\Delta_Y$ or vice versa.

In $Z_1^t$ we have a simple path $\gamma_1$ resulting from $\gamma$ by the cancellation  of the
two identified segments at $x$. By construction $g_t^1$ is simplicial and not injective on
$\gamma_1$. Therefore we can iterate the above procedure and define points $Z_i^t$
with $$\lambda(Z_i^t)\leq \Lip(g_t)=\lambda(Y_t)$$
and such that $\Delta_{Z_i^t}$ has $\Delta_X$ as a simplicial face. Moreover either
$\Delta_{Z_i^t}$ has $\Delta_{Z_{i-1}ìt}$ as a simplicial face or vice versa.
Since $\gamma$ has a finite number of vertices, we must stop, and we do when $\gamma_i$ is a
single point. At this stage, $\alpha$ and $\beta$ are incident to the same point and we are
reduced to the initial case. Note that any $Z_i^t\to X$ as $t\to 0$, thus so does any point
in segment from $Z_i^t$ to $Z_{i+1}^t$. Therefore by lower
semicontinuity of $\lambda$ for any $\e>0$, since we have
finitely many points, for sufficiently small $t$ we have that for any $i$
$$\lambda(X)-\e<\lambda(Z_i^t)$$
and the same inequality holds for points in the segments from $Z_i^t$ to $Z_{i+1}^t$.
\qed
\begin{cor}\label{cordai}
  In the hypothesis of Theorem~\ref{thmZnew}, suppose that $X$ is an exit point for
  $\Delta_X$\footnote{See Definition~\ref{exitp}}, and let $X_E$ as
  Definition~\ref{exitp}. Then, for any $\e>0$, there is $t_\e$ such that for any
  $t<t_\e$, there is a finite  simplicial path $\Sigma_t$ in $\O(Y)$ with vertices
  $Y_t=Z_0^t,Z_1^t,\dots,Z_k^t$ such that $\Delta_{Z_i^t}$ has $\Delta_X$ as a simplicial
face for $i\neq k$, $\Delta_{Z_k^t}$ has $\Delta_{X_E}$ as a simplicial
face, and such that for any point $Z$ of $\Sigma_t$
we have
$$\lambda(X)-\e<\lambda(Z)\leq\lambda(Y_t).$$
Moreover,
for $s\in[0,t]$, $s\mapsto Z_k^s$ parametrizes the segment from from $X_E$ to $Z_k^t$.
\end{cor}
\proof It follows by recursively apply Corollary~\ref{corpeak}. The uniform estimate on $t$
follows because the path from $X$ to $X_E$ is finite.\qed

\section{The end of the proof of Theorem~\ref{tconnected}: peak reduction on simplicial paths}
We fix $\Gamma$ as in Notation~\ref{not:gamma} and $\phi\in\Aut(\Gamma)$. Let
$\lambda=\lambda_\phi$.

We will prove that for any $\e > 0$, the set
$$\{X\in\overline{\O(\Gamma)}^\infty:\lambda(\phi)\leq\lambda_\phi(X)\leq \lambda(\phi) + \e \}$$
is connected by $\lambda(\phi) + \e$-calibrated simplicial paths. This in particular gives the second claim of
Theorem~\ref{tconnected}.

Moreover, if $\Sigma$ is calibrated, then by possibly adding some
extra vertices to $\Sigma$ we obtain a path in the same level set that satisfies
the hypotheses of Theorem~\ref{tregge} and therefore can be regenerated to $\O(\Gamma)$.
Therefore, this proves also the first claim of Theorem~\ref{tconnected}.

\medskip

From now on we fix $A,B\in\overline{\O(\Gamma)}$ such that $\lambda(A),\lambda(B)\geq\lambda(\phi)$.
Let $L \geq \max\{\lambda(A),\lambda(B)\}$.

Let $\Sigma_L(A,B)$ be the set of $L$-calibrated simplicial paths from $A$ to $B$.

\begin{lem}
For some $L$, $\Sigma_L(A,B)\neq\emptyset$.
\end{lem}
\proof Since $\lambda(A),\lambda(B)\geq\lambda(\phi)$, they have not jumped.
Let $A'\in\Hor(A)$ and $B'\in\Hor(B)$. Since $A',B'\in\O(\Gamma)$, which is connected,
there is a simplicial path in $\O(\Gamma)$ between $A',B'$. After minimizing such path, by Lemma~\ref{calibration} we obtain
an element of $\Sigma_L$ (where the $L$ is the maximum displacement along such a path).\qed

\begin{defn}
  For any simplicial path $\Sigma=(X_i)$ we define
  $\max(\Sigma)=\max_{X_i}\lambda(X_i)$, and we say that $X_i$ is a {\em peak} if
  $\lambda(X_i)=\max(\Sigma)$. A pair of two consecutive peaks $X_{i-1},X_{i}$ is called a {\em flat
    peak}. A peak is {\em strict} if it is not part of a flat peak.
\end{defn}

Let $\Sigma_0=(X_i)\in\Sigma(A,B)$ such that among all elements $\sigma\in\Sigma$ it minimizes, in order
\begin{enumerate}
\item $\max(\Sigma)$
\item the number peaks;
\item the number of flat peaks.
\end{enumerate}

\begin{lem}
  Such a $\Sigma_0$ exists.
\end{lem}
\proof By Theorem~\ref{conj} we are minimizing over a well-ordered set.\qed

\medskip

Note that if $X$ is a strict peak then $\lambda$ it is strictly monotone on both sides of $X$.
(By Lemma~\ref{lconvexity}.)

Once again, we need the inductive hypothesis.

\begin{lem}
  Suppose that Theorem~\ref{tconnected} is true in any rank less than $\rank(\Gamma)$.
  Then $\Sigma_0$ has no strict peaks in its interior.
\end{lem}

\proof
Suppose that $\lambda(X_{i-1})<\lambda(X_i)>\lambda(X_{i+1})$. In particular we have
$\lambda(\phi)<\lambda(X_i)$; a strict inequality. By calibration $\Delta_{X_i}$ minimizes $\lambda$ in its simplex,
hence $\Delta_{X_i}$ is a proper face of both $\Delta_i$
and $\Delta_{i+1}$.
Thus for any $Y$ and $Z$, respectively in $\overline{X_{i-1}X_i}$ and
$\overline{X_iX_{i-1}}$ we have
$$\lambda(\phi)<\lambda(Y),\lambda(Z)<\lambda(X_i).$$

We set $X=X_i$. If $C$ is the collapsed part of $X$, then by
Corollary~\ref{corlalx}
$$\lambda(\phi|_C)<\lambda(X).$$
As this is an open condition, it is preserved in an open neighbourhood $U$ of $X$ in
$\O(X)$.\footnote{Note that $\O(X)$ may be different from $\O(\Gamma)$.}

Since  $X$ is not a $\phi$-minimally displaced point, by Lemma~\ref{dom5}
$X\notin\TT(\phi)\subset\O(X)$. By
Lemma~\ref{LemmaX}, $X$ is an exit point. Let $X_E$ as in Definition~\ref{exitp}.

Now we invoke Corollary~\ref{cordai}.
With the terminology of Corollary~\ref{cordai}, let $Y_t$ parametrize the segment from $X$ to
$Y$. By Corollary~\ref{cordai} there exists a simplicial path $\Sigma_Y$ in $\O(Y)$ connecting
$Y_t$ to a point in $Y_E\in\Hor(X_E)$ (as a subset of $\O(Y)$).

Moreover, Corollary~\ref{corpeak} gives the estimate
$$\lambda(X)-\e<\lambda(P)\leq\lambda(Y_t).$$
for any point of $\Sigma_Y$. In particular, if $\e$ is mall enough we have
$$\lambda(\phi)<\lambda(P)<\lambda(X)$$
Strict inequalities. Similarly, there is a simplicial path $\Sigma_Z$ connecting $Z_t$ to a
point $Z_E\in\Hor(X_E)$ with the same estimate above.

Since $\Sigma_Y$ is in $O(Y)$ and
$\lambda$ is continuous on $\O(Y)$, then $\lambda$ is continuous on $\Sigma_Y$.
The same for $\Sigma_Z$.

Let $\wt Y_E,\wt Z_E\in\O(\Gamma)$ be points respectively in $\Hor(Y_E)\cap\Hor(X_E)$ and
$\Hor(Z_E)\cap\Hor(X)$ such that
$\lambda(\wt Y_E),\lambda(\wt Z_E)\leq\lambda(X)$
and such that $\lambda$ is continuous on $\overline{Y_E\wt Y_E}$ and $\overline{Z_E\wt Z_E}$
(such points exist by Lemma~\ref{lemma9}).

By Lemma~\ref{inhor} there is a  path $\Theta_E$ in $\O(\Gamma)$ and in the level set
 $\{\lambda(P)<\lambda(X)\}$ connecting $\wt Y_E$ to $\wt Z_E$. We add $Y_E$ and $Z_E$ to such
 points, we minimize and then apply Lemma~\ref{calibration}. The result is a calibrated
 path $\theta_E$ connecting $Y_t$ to $Z_t$ in the level set $\{\lambda < X\}$.

The path obtained by following $\Sigma$ till $Y_t$, then $\theta_E$ till $Z_t$ and then
$\Sigma$ again, has either a lower maximum than $\Sigma$ or one peak less than
$\Sigma$.\qed

\begin{lem}
  $\Sigma_0$ has no flat peaks unless $\lambda$ is constant on $\Sigma$ and $\lambda(\Sigma)=\lambda(\phi)$.
\end{lem}
\proof If $\lambda$ is not constantly $\lambda(\phi)$ on $\Sigma$, in particular $\lambda$ is
strictly bigger than $\lambda(\phi)$ on peaks.

Suppose that there is $Y,X$ two consecutive vertices of $\Sigma_0$ with
$\lambda(X)=\lambda(Y)=\max(\Sigma)>\lambda(\phi)$. The idea is to find a third point $Z$ to
add between $Y$ and $X$ in order to destroy the flat peak.

If there is a point $Z$ in the interior of the segment $YX$, with $\lambda(\phi)\leq\lambda(Z)<\lambda(X)=\lambda(Y)$, then we
add it.

Otherwise, $\lambda$ is constant on $\overline{XY}$. Let $W$ be a point in the interior of the
segment $\overline{XY}$. If $W$ is not a local minimum for $\lambda$ in $\Delta_W$, then near
$W$ we find $Z$ with the above properties. We add it.

If $W$ is a local minimum for $\lambda$ in $\Delta_W$ then, by
Lemma~\ref{LemmaX}, near $W$ in $\O(W)$  there is a point $Z$ with the above properties and such that $\Delta_W$
is a finitary face of $\Delta_Z$ in $\O(W)$. We add $Z$.

Since $\lambda$ is continuous on closed simplices of $\O(W)$, then by adding $Z$ in each of the
above case we obtain a new $L$-calibrated path $\Sigma_1$ which has the same number of peaks
and exactly one flat peak less then $\Sigma_0$. This is impossible by the minimality of $\Sigma_0$.
\qed

\medskip
To finish the proof of Theorem~\ref{tconnected}, simply observe that we have shown that we can
connect any two points in $\{X\in\overline{\O(\Gamma)}^\infty:\lambda_\phi(X)=\lambda(\phi)\}$ by a calibrated simplicial path with no peaks, either strict or  flat. This immediately implies that the displacement is constant along the path. \qed

\section{Applications}\label{appl}

In this section we show how the connectedness of the level sets gives a solution to the conjugacy problem. 

We start with some technical results. Recall that a point, $X$, of $CV_n$ is called $\epsilon$-thin if there is a homotopically non-trivial loop in $X$ of length at most $\epsilon$. Conversely, $X$ is called $\epsilon$-thick if it is not $\epsilon$-thin.

\begin{prop}[\cite{BestvinaBers}, Proposition 10] \label{thick}
Let $X \in CV_n$ (that is, $X$ is a marked metric graph) and $f:X \to X$ a PL-map representing some automorphism of $F_n$. Let $\lambda  = Lip(f)$, let $N$ equal maximal chain of topological subgraphs of any graph in $CV_n$ (this is clearly a finite number) and let $\mu$ be any real number greater than $\lambda$. Then if $X$ is $1/((3n-3)\mu^{(N+1)})$-thin, the automorphism represented by $f$ is reducible. For instance, one can take $N=3n-3$.
\end{prop}

\begin{defn}
Let $X \in CV_n$. Then we call $R$ an adjacent uniform rose if it obtained by collapsing a maximal tree in $X$ and then rescaling so that all edges in $R$ have the same length (that is, $1/n$, as we will work with volume 1). 
\end{defn}

\begin{prop}
Let $X \in CV_n$ be a point which is $\epsilon$-thick and let $R$ be any adjacent uniform rose (both of volume 1). 
Then, $\Lambda(X, R) \leq 1 / \epsilon$ and $\Lambda(R, X) \leq n$. 
\end{prop}
\proof By Theorem~\ref{sausagelemma}, we can look at candidates that realise the stretching factor. Since, topologically, one passes from $X$ to $R$ by collapsing a maximal tree, we get that a candidate in $X$, when mapped to $R$, crosses every edge at most twice. In fact the candidate crosses every edge of $R$ at most once in the case of an embedded simple loop or an infinity loop. This gives the first inequality, on taking into account that $X$ is $\epsilon$-thick and that barbells have length at least $2 \epsilon$. 

For the second inequality note that a embedded loop in $R_0$ is a edge and has length $1/n$ and lifts to an embedded loop in $X$, of length at most $1$. An infinity loop in $R_0$ consists of two distinct edges, has length $2/n$ and lifts to a loop in $X$ which goes through every edge at most twice. (Barbells are not present in $R_0$). 
\qed

\begin{cor} \label{roseapprox}
Let $X \in CV_n$ be $\epsilon$-thick and let $R$ be an adjacent uniform rose. Consider $\Phi \in Out(F_n)$. Then $\Lambda(R, \phi R) \leq \frac{ n }{\epsilon}\Lambda(X, \phi X)$.  
\end{cor}

\begin{prop} \label{roseconnect}
Let $R, R_{\infty}$ be two points in $CV_n$ which are both uniform roses (graphs with exactly one vertex and so that every edge has the same length). Let $\phi \in Out(F_n)$ be irreducible and suppose that $\mu$ is any real number greater than $\max( \Lambda(R, \phi R), \Lambda(R_{\infty}, \phi R_{\infty}))$.

Then there exist $R_0=R, R_1, R_2, \ldots, R_k=R_{\infty}$, which are all uniform roses in $CV_n$ such that: 
\begin{itemize}
\item For each $i$, there exists a simplex $\Delta_i$ such that $\Delta_{R_i}$ is a rose face of both  $\Delta_i$ and $\Delta_{i+1}$. 
\item $\Lambda(R_i, \phi R_i) \leq \frac{ n }{\epsilon} \mu$, where $\epsilon = 1/((3n-3)\mu^{(N+1)})$. 
\end{itemize}
\end{prop}
\proof
This follows from Theorem~\ref{tconnected}, using Definition~\ref{def:sp}, since each pair $\Delta_i$ and $\Delta_{i+1}$ have a (at least one) common rose face; just take any uniform adjacent rose in any  common rose face. The remaining point follows from Corollary~\ref{roseapprox} and Proposition~\ref{thick}. 
\qed

\medskip

\noindent
{\em Proof of Theorem~\ref{conjirred}}: We clearly have an algorithm which terminates, and it is apparent that if $\psi \in S_{\phi}$ then these automorphisms are conjugate. It remains to show the converse; that if they are conjugate, then $\psi \in S_{\phi}$.

Let $R$ be the uniform rose corresponding to the basis $B$. If $\psi$ were conjugate to $\phi$, then there would be a conjugator, some $\tau \in Out(F_n)$ such that $\psi = \tau^{-1} \phi \tau$. Let $R_{\infty} = \tau R$. Now use Proposition~\ref{roseconnect} to find a sequence $R=R_0, R_1, \ldots, R_k=R_{\infty}$, such that each consecutive pair are incident to a common simplex and $\Lambda(R_i, \phi R_i) \leq n(3n-3) \mu^{3n-1} = K$.

We let $\zeta_i$ be an automorphism which sends $R_i$ to $R_{i+1}$; the fact that these roses are both incident to a common simplex implies that each $\zeta_i$ is a CMT automorphism. Inductively, we may define, $\tau_i = \zeta_0 \ldots \zeta_{i-1}$, and note that $\tau_i R = R_i$. We make these choices so that $\tau=\tau_k$. (This possible since regardless of the choices made, we always have that $\tau_k^{-1} \tau$ fixes $R$ and is therefore a CMT automorphism, therefore by possibly adding a single repetition of roses at the start we may assume that $\tau=\tau_k$.)

Now let $\phi_i =  \tau_i^{-1} \phi \tau_i$. 

Since $\phi_{i+1} = \zeta_i^{-1} \phi_i \zeta$, to finish the proof we just need that $||\phi_i || \leq K$. This follows since, 
$$
\Lambda( R_i, \phi R_i ) = \Lambda(\tau_i R, \phi \tau_i R) = \Lambda(R, \phi_i R) = ||\phi_i||_B. 
$$

\qed

We conclude the paper by proving Theorem~\ref{detectirred}. First a lemma, 

\begin{lem}\label{reducerose}
Let $X$ be a core graph and $f$ a homotopy equivalence on $X$, having a proper, homotopically non-trivial subgraph $X_0$ such that $f(X_0) = X_0$.  Then there is a maximal tree, $T$, such that the automorphism induced by $f$ on the rose $X/T$ is visibly reducible.  
\end{lem}
\proof
Choose $X_0$ to be minimal. Therefore it will have components, $X_1, \ldots, X_k$ such that $f(X_i)=X_{i+1}$ with subscripts taken modulo $k$. Take a maximal tree for each $X_i$ and extend this to a maximal tree, $T$, for $X$. It is then clear that if we take $B_i$ to be the set of edges in $X/T$ coming from $X_i$, that $\psi$ will be visibly reducible as witnessed by $B_1, \ldots, B_k$. (Note each subgroups generated by each $B_i$ are only permuted/preserved up to conjugacy, since the $X_i$ are disjoint and so one cannot choose a common basepoint).  
\qed

\medskip

\noindent
{\em Proof of Theorem~\ref{detectirred}}: We proceed much as in the proof of Theorem~\ref{conjirred}, but here we do not know that the points in $CV_n$ we encounter will remain uniformly thick. 

The algorithm clearly terminates, and if there is a $\psi$ in $S^+$ which is visibly reducible, then $\phi$ is reducible. It remains, therefore, to show that if $\phi$ is reducible, then there is some $\psi \in S^+$ which is visibly reducible. 

Let $R$ be the uniform rose corresponding to the basis $B$. By Corollary~\ref{strongcorred}, there exists an $X \in CV_n$ with a core invariant subgraph and such that $\Lambda(X, \phi(X)) \leq \mu$. 

By  Theorem~\ref{tconnected}, there exist points, $X_0=R, X_1, \ldots, X_k=X$, such that $\Lambda(X, \phi(X_i)) \leq \mu$. Choose the maximal index, $M$, such that $X_0, X_1, \ldots, X_M$ are all $\epsilon$-thick, where $\epsilon = 1/((3n-3)\mu^{(N+1)})$ as in Lemma~\ref{thick}. Now for each $i \leq M$, choose an $R_i$ which is a adjacent uniform rose to both $X_i$ and $X_{i+1}$ (choose $R_0=R$ and if $M=k$, let $R_k$ be any uniform rose adjacent to $X_k$). 

If $M=k$, we set $R_{M+1}=R_M$. Otherwise, by Lemma~\ref{thick}, we have that $X_{M+1}$ has an optimal $PL$-representative for $\phi$ which admits an invariant subgraph. So by Lemma~\ref{reducerose}, we may find an adjacent uniform rose face, $R_{M+1}$ so that the representative, $\psi$, of $\phi$ at $R_{M+1}$ is visibly reducible. 

As above,  we $\tau \in Out(F_n)$ such that $\psi = \tau^{-1} \phi \tau$. Then let $\zeta_i$ be an automorphism which sends $R_i$ to $R_{i+1}$; each $\zeta_i$ is a CMT automorphism. Inductively, we may define, $\tau_i = \zeta_0 \ldots \zeta_{i-1}$, and note that $\tau_i R = R_i$. We make these choices so that $\tau=\tau_{M+1}$.

Now let $\phi_i =  \tau_i^{-1} \phi \tau_i$, so that $\phi_0=\phi$ and $\phi_{M+1} = \psi$. Since each $X_0, \ldots, X_M$ is $\epsilon$-thick we get, by Corollary~\ref{roseapprox}  that each $\phi_i \in S_i$ for $i \leq M$. Hence $\psi \in S^+$ and is visibly reducible.
\qed

\newpage
\section{Appendix: proof of Theorem~\ref{thmZnew}}\label{appendix}
In this section we give the proof of Theorem~\ref{thmZnew}, which we recall
\begin{thm*}[Theorem~\ref{thmZnew}]
Let $X,Y\in\overline{\O(\Gamma)}$. Suppose that $\Delta_X$ is a simplicial face of
$\Delta_Y$. Thus as graphs, $Y$ is obtained by collapsing a sub-graph $A$. Suppose that
$\core(A)$ is $\phi$-invariant. For $t\in[0,1]$ let $Y_t=(1-t)X +t Y$ be a
parametrization of the Euclidean segment from $X$ to $Y$. Let $\sigma_t:Y_t\to X$ be the map
obtained by collapsing  $A$ and by linearly rescaling the edges in $Y\setminus A$.

Let $f:X\to X$ be an optimal map representing $\phi$. Then for any $\e >0$ there is $t_\e>0$ such that $\forall 0\leq t<t_\e$ there is an optimal map
$g_t:Y_t\to Y_t$ representing $\phi$ such that $$d_\infty(\sigma_t\circ g_t,f\circ \sigma_t)<\e.$$
\end{thm*}
\proof We split the proof in two sub-cases. First when $A$ is itself a core-graph, and then the
case when $\core(A)$ is empty. Clearly the disjoint union of the two cases implies the mixed case.

We will work at once with graphs and
trees, by using the usual $\wt\ $ notation: if $X$ is a graph, $\wt X$ is its universal
covering an for any object $o$ (a point, a sub-set, a map, $\dots$) $\tilde o$ is one of its
lift to the universal coverings. Vice versa, given an object $\tilde o$, we understand that
$\tilde o$ descends (by equivariance if for example $o$ is a map) to an object $o$ at the level of graphs.

\begin{lem}\label{lemmaZnew}
Let $X,Y\in\overline{\O(\Gamma)}$. Suppose that as graphs of groups, $X$ is obtained from $Y$
by collapsing a $\phi$-invariant core sub-graph $A=\sqcup A_i$.
For $t\in[0,1]$ let $Y_t=(1-t)X +t Y$ be a
parametrization of the Euclidean segment from $X$ to $Y$. Let $\sigma_t:Y_t\to X$ be the map
obtained by collapsing  $A$ and by linearly rescaling the edges in $Y\setminus A$.

Let $f:X\to X$ be an optimal map representing $\phi$. Then for any $\e >0$ there is $t_\e>0$ such that $\forall 0\leq t<t_\e$ there is an optimal map
$g_t:Y_t\to Y_t$ representing $\phi$ such that $$d_\infty(\sigma_t\circ g_t,f\circ \sigma_t)<\e.$$
\end{lem}
\proof  We begin by fixing  some notation. Firs of all, we will use the symbol $\lambda$ to
denote any displacement functions of $\phi$ (i.e. $\lambda_\phi,\lambda_{\phi|_A},\dots$)
If $x$ is a point in a metric space, we denote by $B_r(x)$ the open metric ball centered at $x$
and radius $r$.
For any $i$, we denote by $v_i$ the non-free vertex of $X$
obtained by collapsing $A_i$. For any $t$ we denote by $A^t$ the metric copy of $A$ in
$Y_t$. Note that $A$ is uniformly collapsed in $Y_t$, that is to say, $[A^t]\in\mathbb P\O(A)$
is the same element for any $0<t\leq 1$, and we have $\vol(A^t)=t\vol(A^1)$.

By lower
semicontinuity of $\lambda$ (Theorem~\ref{fatto1}) we have that
\begin{equation}
  \label{eq:0}
  \forall \e_0>0\exists t_{\e_0}>0 \text{ such that } \forall t<t_{\e_0} \text{ we have }
\lambda (Y_t)>\frac{\lambda (X)}{1+\e_0}.
\end{equation}

A priori $f$ may collapse some edge, in any case   $\forall \e_1>0\exists f_1:X\to X$
a $\PL$- map representing $\phi$ such that $f_1$ does not collapse any edge, and
\begin{equation}
  \label{eq:1}
  d_\infty(f,f_1)<\e_1\quad \text{and}\quad \Lip(f_1)<\Lip(f)(1+\e_1)=\lambda(X)(1+\e_1).
\end{equation}

Moreover $\exists 0<\rho_0=\rho_0(X,f_1)$ such that $\forall \rho<\rho_0$
\begin{itemize}
\item $B_\rho(x)$ is star-shaped for any $x\in X$ (i.e. it contains at most one vertex);
\item for any $i$, each connected component of $f_1^{-1}(B_\rho(v_i))$ is star-shaped
  and contains exactly one pre-image of $v_i$;
\item for any $i,j$ the connected components of
  $f_1^{-1}(B_\rho(v_i))$ and those of $f_1^{-1}(B_\rho(v_j))$ are pairwise disjoint.
\end{itemize}

We fix an optimal map $\varphi:A^1\to A^1$ representing $\phi|_A$. Since $[A^t]\in\mathbb
P\O(A)$ does not depend on $t$, $\varphi:A^t\to A^t$ is an optimal map for any $t\in(0,1]$ and the
Lipschitz constant does not change. Clearly (by Sausage
Lemma~\ref{sausagelemma})
\begin{equation}
  \label{eq:fi}
   \Lip(\varphi)\leq \lambda(Y_t)\quad \text{ for any } t.
\end{equation}
The natural option is to define $g_t$ by using $\sigma_t^{-1}\circ
f_1\circ\sigma_t$. Hence, we need to deal with places where $\sigma_t^{-1}$ is not
defined. First we fix a lift $\wt\varphi$ of $\phi$.

Each germ of edge $\alpha$ at $v_i$ in $X$ corresponds to a germ $\alpha_Y(=\sigma^{-1}_t(\alpha))$ in $Y$ incident to
$A_i$ at a point that we denote by $p_\alpha$. For any such $\alpha$ we {\bf choose} a lift
$\wt\alpha$, that corresponds to a germ $\wt\alpha_Y$ incident to $\wt p_\alpha\in\wt
A_i$. (See Figure~\ref{pippo}.)
\setlength{\unitlength}{1ex}
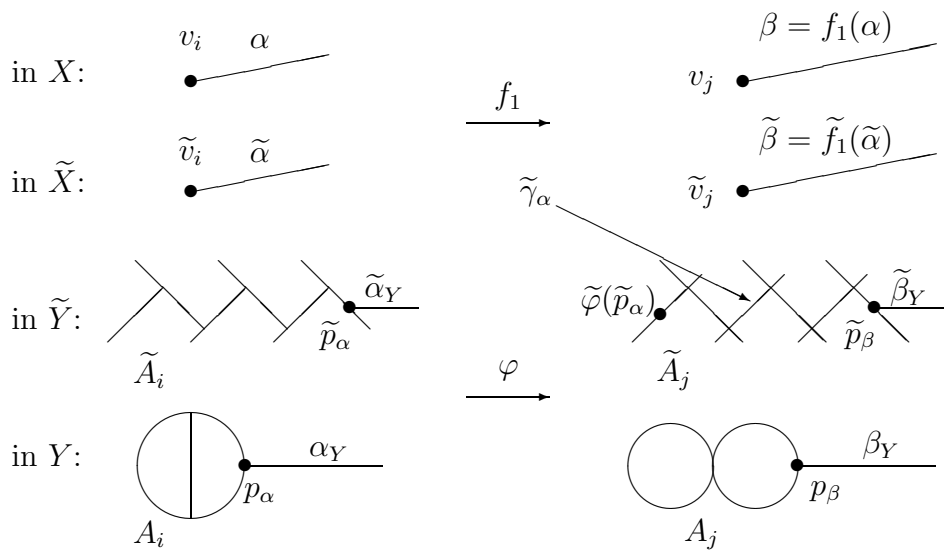
\begin{figure}[htbp]
  \centering
  \begin{picture}(70,38)
\put(8,0){
    \put(5,5){\circle{10}}
    \put(5,1.1){\line(0,1){7.8}}
    \put(8.9,5){\line(1,0){10}}
    \put(8.9,5){\makebox(0,0){$\bullet$}}
    \put(2,0){\makebox(0,0){$A_i$}}
    \put(10,3){\makebox(0,0){$p_\alpha$}}
    \put(15,6){\makebox(0,0){$\alpha_Y$}}

    \put(40,0){
      \put(-0.25,5){\circle{5.85}}
      \put(5.9,5){\circle{5.85}}
      \put(8.9,5){\line(1,0){10}}
      \put(8.9,5){\makebox(0,0){$\bullet$}}
      \put(2,0){\makebox(0,0){$A_j$}}
      \put(11,3){\makebox(0,0){$p_\beta$}}
      \put(15,6.5){\makebox(0,0){$\beta_Y$}}
      \put(-2,0){\multiput(0,0)(6,0){3}{
                  \put(-1,14){\line(1,1){5}}
                  \put(1,20){\line(1,-1){6}}}
                \put(16.5,16.5){\makebox(0,0){$\bullet$}}
                \put(16.5,16.5){\line(1,0){5}}
                \put(15.5,14.5){\makebox(0,0){$\wt p_\beta$}}
                \put(19,18){\makebox(0,0){$\wt \beta_Y$}}
                \put(1,16){\makebox(0,0){$\bullet$}}
                \put(-2,17){\makebox(0,0){$\wt\varphi(\wt p_\alpha)$}}
      \put(2,12){\makebox(0,0){$\wt A_j$}}
                }
    \put(5,25){\line(5,1){14}}
    \put(5,25){\makebox(0,0){$\bullet$}}
    \put(11,29){\makebox(0,0){$\wt \beta=\wt f_1(\wt \alpha)$}}
    \put(2,25){\makebox(0,0){$\wt v_j$}}

    \put(0,8){
      \put(5,25){\line(5,1){14}}
      \put(5,25){\makebox(0,0){$\bullet$}}
      \put(11,29){\makebox(0,0){$\beta=f_1(\alpha)$}}
      \put(2,25){\makebox(0,0){$v_j$}}
               }

      }

    \multiput(0,0)(6,0){3}{
       \put(-1,14){\line(1,1){4}}
       \put(1,20){\line(1,-1){5}}}
      \put(16.5,16.5){\makebox(0,0){$\bullet$}}
      \put(16.5,16.5){\line(1,0){5}}
      \put(15.5,14.5){\makebox(0,0){$\wt p_\alpha$}}
      \put(19,18){\makebox(0,0){$\wt \alpha_Y$}}
      \put(2,12){\makebox(0,0){$\wt A_i$}}

    \put(5,25){\line(5,1){10}}
    \put(5,25){\makebox(0,0){$\bullet$}}
    \put(10,28){\makebox(0,0){$\wt\alpha$}}
    \put(5,28){\makebox(0,0){$\wt v_i$}}

    \put(0,8){
      \put(5,25){\line(5,1){10}}
      \put(5,25){\makebox(0,0){$\bullet$}}
      \put(10,28){\makebox(0,0){$\alpha$}}
      \put(5,28){\makebox(0,0){$v_i$}}
               }
\put(30,25){\makebox(0,0){$\wt\gamma_\alpha$}}
\put(31.5,24){\vector(2,-1){14}}
}
\put(0,5){in $Y$:}
\put(0,15){in $\wt Y$:}
\put(0,25){in $\wt X$:}
\put(0,33){in $X$:}
\put(33,30){\vector(1,0){6}}
\put(36,32){\makebox(0,0){$f_1$}}
\put(33,10){\vector(1,0){6}}
\put(36,12){\makebox(0,0){$\varphi$}}

  \end{picture}
  \caption{How to  choose the paths $\wt\gamma_\alpha$}\label{pippo}
\end{figure}

Suppose $f_1(v_i)=v_j$, and let $\beta=f_1(\alpha)$. Then $\wt\beta$ is a germ at $\wt v_j$ and
corresponds to a germ $\wt\beta_Y$ incident to $\wt A_j$ at a point $\wt p_\beta$.

Let $\wt \gamma_\alpha$ be the unique path in $\wt A_j$ connecting $\wt\varphi(p_\alpha)$ to $\wt
p_\beta$.
\begin{rem}
  We choose a path $\wt \gamma_\alpha$ for any germ $\alpha$ in $X$, which is a finite graph.
Therefore we have only finitely many such $\wt \gamma_\alpha$'s. We can then complete that family of
paths by equivariance.
\end{rem}

Now we do a similar construction for other pre-images of the $v_i$'s. For any $x\in X$ such
that $f_1(x)=v_i$ for some $i$, but $x\notin\{v_j\}$, we choose a base-point $\wt x_i\in \wt
A_i$. Any germ of edge $\alpha$ at $x$ correspond to an edge $\alpha_Y$ is $Y$
(note that $x$ is not necessarily a vertex of $X$). For any such $\alpha$ we {\bf choose} a lift $\wt
\alpha$. Since $f_1$ does not collapse edges, $\wt f_1(\wt \alpha)$ is  a germ of edge $\wt
\beta$ at $\wt v_i$, and
corresponds to a germ $\wt \beta_Y$ at $\wt A_i$ in $\wt Y$. Let $\wt \gamma_\alpha$ be the unique
path in $\wt A_i$ connecting $\wt x_i$ and $\wt \beta_Y$.

\begin{rem}
As above  we choose only finitely many such $\wt\gamma_\alpha$'s and we complete the choices
equivariantly.
\end{rem}
Note that, as germs, $\alpha_Y=\sigma_t^{-1}(\alpha)$
and $\beta_Y=\sigma_t^{-1}(\beta)=\sigma_t^{-1}(f_1(\alpha))$.
Now we have a path $\gamma_\alpha\subset A$ for any pre-image of germs at the $v_i$'s, chosen
independently on $t$. Let
$t\in(0,1]$. We define a map $$g: Y_t\to Y_t$$ representing $\phi$ as follows:
\begin{itemize}
\item in $\sigma_t^{-1}\big(X\setminus f_1^{-1}(\sqcup_i B_\rho(v_i))\big)$ we just set
  $g=\sigma_t^{-1}\circ f_1\circ\sigma_t$;
\item in $\sigma_t^{-1}\big(f_1^{-1}(\sqcup_i B_\rho(v_i))\big)\setminus A^t$ we use the paths
  $\gamma_\alpha$. More precisely, let $N$ be a connected component of
  $f_1^{-1}(B_\rho(v_i))$  and let $x\in N$ such that $f_1(x)=v_i$. For any
  edge $\alpha\in N$ emanating from $x$ we define $g(\sigma_t^{-1}(\alpha))$ by mapping
  linearly\footnote{I.e. at constant speed} $\sigma_t^{-1}(\alpha)$ to the path given by the
  concatenation of
  $\beta_Y=\sigma_t^{-1}(f_1(\alpha))$ and $\gamma_\alpha$. Note that
  $g|_{\sigma_t^{-1}(\alpha)}=\PL(g|_{\sigma_t^{-1}(\alpha)})$.
\item in $A^t$ we set $g=\varphi$;

\end{itemize}
finally, we set $$g_t=\opt(\PL(g))$$
where $\PL$-ization and optimization are made with respect to the metric structure of $Y_t$.
We now estimate the Lipschitz constant of $g$. Clearly $$\lambda(Y_t)=\Lip(g_t)\leq \Lip(g).$$

In $\sigma_t^{-1}\big(X\setminus f_1^{-1}(\sqcup_i B_\rho(v_i))\big)$ we have
  $g=\sigma_t^{-1}\circ f_1\circ\sigma_t$. Then
$$\Lip(g)\leq \Lip(\sigma_t^{-1})\Lip(f_1)\Lip(\sigma_t).$$

Since on edges of $Y\setminus A$ the map $\sigma_t$ is just a rescaling of edge-lengths, for any $\e_2>0$ there is $t_{\e_2}>0$ such that $\forall t<t_{\e_2}$
\begin{equation}
  \label{eq:2}
  \Lip(\sigma_t)<1+\e_2 \qquad \Lip(\sigma_t^{-1})<1+\e_2
\end{equation}
hence, by~\eqref{eq:1}, and by setting $(1+\e_2)^2(1+\e_1)=1+\e_3$ we have
\begin{equation}
  \label{eq:3}
  \Lip(g)\leq(1+\e_2)^2\lambda(X)(1+\e_1)=(1+\e_3)\lambda(X).
\end{equation}

Now, let $N$ be a connected component of $f_1^{-1}(\sqcup_i B_\rho(v_i))$. Let $x\in N$ such
that $f_1(x)=v_i$ and let $\alpha$ be an edge of $N$ emanating from $x$.
By definition $g$ is linear on $\sigma_t^{-1}(\alpha)$, thus in order to estimate its Lipschitz
constant we need to know only the lengths of $\sigma_t^{-1}(\alpha)$ and  its image.
We have $L_X(f_1(\alpha))=\rho$ and therefore
$$\rho\leq \Lip(f_1) L_X(\alpha)
\qquad
L_X(\alpha)=L_X(\sigma_t(\sigma_t^{-1}(\alpha)))
\leq
\Lip(\sigma_t) L_{Y_t}(\sigma_t^{-1}(\alpha))
$$ whence, by~\eqref{eq:2} and~\eqref{eq:3}, we obtain
$$L_{Y_t}(\sigma_t^{-1}(\alpha))\geq\frac{L_X(\alpha)}{\Lip(\sigma_t)}
>\frac{\rho}{(1+\e_2)\Lip (f_1)}>
\frac{\rho}{\lambda(X)(1+\e_1)(1+\e_2)}.$$
Since $\gamma_\alpha$ is the same loop in $A$ for every $t$, its length in $A^t$ depends
linearly on $t$, namely here is a constant $C_\alpha$ such that $$L_{Y_t}(\gamma_\alpha)= C_\alpha t$$
whence, setting $C=\max_\alpha C_\alpha$,
\begin{eqnarray*}
  &~&\Lip(g|_{\sigma_t^{-1}(\alpha)})
\leq\frac{L_{Y_t}(\sigma_t^{-1}(f_1(\alpha))+L_{Y_t}(\gamma_\alpha)}{L_{Y_t}(\sigma_t^{-1}(\alpha))}
\leq \frac{\Lip(\sigma_t^{-1})\rho+tC_\alpha}{L_{Y_t}(\sigma_t^{-1}(\alpha))}\\
&<&((1+\e_2)\rho+tC_\alpha)\frac{\lambda(X)(1+\e_1)(1+\e_2)}{\rho}\\&=&
\lambda(x)\big[(1+\e_3)+\frac{(1+\e_1)(1+\e_2)+tC_\alpha}{\rho}\big]\\
&<&\lambda(x)\big[(1+\e_3)+\frac{(1+\e_3)+tC}{\rho}\big]
\end{eqnarray*}
Therefore $\forall\e_4>0\exists t_{\e_4}>0$ such that $\forall t<t_{\e_4}$, for any $\alpha$ as
above we have

\begin{eqnarray}
  \label{eq:4}
  \Lip(g|_{\sigma_t^{-1}(\alpha)})<\lambda(X)(1+\e_4).
\end{eqnarray}

Finally,  on $A^t$ we have $g=\varphi$ and so $\Lip(g|_{A^t})=\Lip(\varphi)$.
Since by~\eqref{eq:0} $\lambda(X)\leq\lambda(Y_t)(1+\e_0)$, by putting
together~\eqref{eq:fi},~\eqref{eq:3}, and~\eqref{eq:4} we have that for any $\e_5>0$ there is
$t_{\e_5}>0$ such that for any $t<t_{\e_5}$ we have

$$\Lip(g)\leq\lambda(Y_t)(1+\e_5)$$

Since $g_t$ is optimal $\Lip(g_t)=\lambda(Y_t)$ and by Theorem~\ref{Lemma_opt}
$$d_\infty(g_t,g)<\vol(Y_t)(\Lip(g)-\lambda(Y_t))<\vol(Y_t)\lambda(Y_t)\e_5.$$

\medskip

We now estimate $$d_\infty(\sigma_t\circ g,f_1\circ \sigma_t).$$
 In $\sigma_t^{-1}\big(X\setminus f_1^{-1}(\sqcup_i B_\rho(v_i))\big)$ we have
  $g=\sigma_t^{-1}\circ f_1\circ\sigma_t$ so here the distance is zero. On $A^t$, since $g(A)=A$,
  for any $i$ there is $j$ such that  we have $\sigma_t(g(A_i))=\sigma_t(A_j)=v_j=f(v_i)$,
  hence also in $A^t$ the distance is zero.
Finally, let $N$ be a connected component of $f_1^{-1}(\sqcup_i B_\rho(v_i))$. Let $x\in N$ such
that $f_1(x)=v_i$ and let $\alpha$ be an edge of $N$ emanating from
$x$.
The path $g(\sigma_t^{-1}(\alpha))$ is given by the concatenation of $\sigma_t^{-1}(f_1(\alpha))$ with
$\gamma_\alpha$. The latter is collapsed by $\sigma_t$, and the image of the former is just
$f_1(\alpha)=f_1\circ\sigma_t(\sigma_t^{-1}(\alpha))$. Since the length of $\gamma_\alpha$ in
$A^t$  is bounded by $tC$ we have that
$$d_\infty(\sigma_t\circ g,f_1\circ \sigma_t)\to 0\qquad\text{ as }t\to 0.$$
In particular $\forall \e_6\exists t_{\e_6}$ such that $\forall t<t_{\e_6}$ we have
$$d_\infty(\sigma_t\circ g,f_1\circ \sigma_t)<\e_6.$$

Finally,
\begin{eqnarray*}
&~&d_\infty(\sigma_t\circ g_t,f\circ \sigma_t)\\&\leq&
d_\infty(\sigma_t\circ g_t,\sigma_t\circ g)+
d_\infty(\sigma_t\circ g,f_1\circ \sigma_t)+
d_\infty(f_1\circ\sigma_t,f\circ \sigma_t)\\&\leq&
\Lip(\sigma_t)d_\infty(g_t,\circ g)+\e_6+d_\infty(f_1,f)\\&<&(1+\e_2)\vol(Y_t)\lambda(Y_t)\e_5+\e_6+\e_1
\end{eqnarray*}
which is arbitrarily small for $t\to 0$.\qed

\begin{lem}\label{lemmaZ}
Let $X,Y\in\overline{\O(\Gamma)}$. Suppose that as graphs of groups, $X$ is obtained from $Y$
by collapsing a sub-forest $T=\sqcup T_i$ whose tree $T_i$ each contains at most one non-free
vertex.
For $t\in[0,1]$ let $Y_t=(1-t)X +t Y$ be a
parametrization of the Euclidean segment from $X$ to $Y$. Let $\sigma_t:Y_t\to X$ be the map
obtained by collapsing  $T$ and by linearly rescaling the edges in $Y\setminus T$.

Let $f:X\to X$ be an optimal map representing $\phi$. Then for any $\e >0$ there is $t_\e>0$ such that $\forall 0\leq t<t_\e$ there is an optimal map
$g_t:Y_t\to Y_t$ representing $\phi$ such that $$d_\infty(\sigma_t\circ g_t,f\circ \sigma_t)<\e.$$
\end{lem}
\proof
The proof goes exactly as that of Lemma~\ref{lemmaZnew}, and it is even simpler. As above $T^t$
denote the scaled version of $T$.
Let $v_i$ be the vertex of $X$ resulting from the collapse of $T_i$.
The function $\lambda$ is now continuous
$$\lambda(Y_t)\to \lambda(X)$$
as above, if $f$ collapses some edge we fine $f_1:X\to X$ a $\PL$-map representing $\phi$ which
collapses no edge and with
$$d_\infty(f,f_1)<\e_1\quad \text{and}\quad \Lip(f_1)<\Lip(f)(1+\e_1)=\lambda(X)(1+\e_1).$$

We choose $\rho$ so that $B_\rho(v_i)$ is star-shaped, the components of
$f_1^{-1}(B_\rho(v_i))$ are star-shaped and contain a unique pre-image of $v_i$, and so that
the  components of $f_1^{-1}(B_\rho(v_i))$  and $f_1^{-1}(B_\rho(v_j))$ are pairwise
disjoint. Finally we chose $\rho$ small enough so that if $f(v_i)\notin\{v_j\}$,
then $f(v_i)\notin\cup_jB_\rho(v_j)$.

For any $i$ we choose a base vertex $x_i\in T_i$ which  the non-free vertex of $T_i$ if
any. For any $x\in X$ such that $f_1(x)=v_i$ and for any edge $\alpha$ in $f_1^{-1}(B_\rho(v_i))$
incident to $x$, let $\gamma_\alpha$ be the unique embedded path connecting
$\sigma_t^{-1}(f_1(\gamma_\alpha))$ to $x_i$. We define $g:Y_t\to Y_t$ as follows:
\begin{itemize}
\item in $\sigma_t^{-1}\big(X\setminus f_1^{-1}(\sqcup_i B_\rho(v_i))\big)$ we just set
  $g=\sigma_t^{-1}\circ f_1\circ\sigma_t$;
\item in $\sigma_t^{-1}\big(f_1^{-1}(\sqcup_i B_\rho(v_i))\big)\setminus T^t$ we use the paths
  $\gamma_\alpha$. More precisely, let $N$ be a connected component of
  $f_1^{-1}(B_\rho(v_i))$  and let $x\in N$ such that $f_1(x)=v_i$. For any
  edge $\alpha\in N$ emanating from $x$ we define $g(\sigma_t^{-1}(\alpha))$ by mapping
  linearly\footnote{I.e. at constant speed} $\sigma_t^{-1}(\alpha)$ to the path given by the
  concatenation of
  $\sigma_t^{-1}(f_1(\alpha))$ and $\gamma_\alpha$. Note that
  $g|_{\sigma_t^{-1}(\alpha)}=\PL(g|_{\sigma_t^{-1}(\alpha)})$.
\item in the components $T_i^t$ so that $f_1(v_i)=v_j$, we set $g(T_i^t)=x_j$;
\end{itemize}
finally we set $g_t=\opt(\PL(g))$. The estimates on Lipschitz constants and distances now
follow exactly as in the proof of Lemma~\ref{lemmaZnew}.\qed

\providecommand{\bysame}{\leavevmode\hbox to3em{\hrulefill}\thinspace}
\providecommand{\MR}{\relax\ifhmode\unskip\space\fi MR }
\providecommand{\MRhref}[2]{%
  \href{http://www.ams.org/mathscinet-getitem?mr=#1}{#2}
}
\providecommand{\href}[2]{#2}

\end{document}